\documentclass[sn-aps]{sn-jnl}

\usepackage{graphicx}%
\usepackage{multirow}%
\usepackage{amsmath,amssymb,amsfonts}%
\usepackage{amsthm}%
\usepackage{mathrsfs}%
\usepackage[title]{appendix}%
\usepackage{xcolor}%
\usepackage{textcomp}%
\usepackage{manyfoot}%
\usepackage{booktabs}%
\usepackage{algorithm}%
\usepackage{algorithmicx}%
\usepackage{algpseudocode}%
\usepackage{listings}%
\usepackage{bm}
\usepackage[caption=false,font=normalsize,labelfont=sf,textfont=sf]{subfig}
\usepackage{stfloats}

\numberwithin{equation}{section}

\newcommand{\vs}{{\bm s}}
\newcommand{\vn}{{\bm n}}
\newcommand{\vp}{{\bm p}}
\newcommand{\vr}{{\bm r}}

\newcommand{\vy}{{\bm y}}
\newcommand{\vz}{{\bm z}}

\newcommand{\vS}{{\bm S}}
\newcommand{\vJ}{{\bm J}}
\newcommand{\vU}{{\bm U}}

\newcommand{\vV}{{\bm V}}

\newcommand{\vmu}{{\bm \mu}}

\newcommand{\mua}{{\bm{\delta\mu}}_a}

\newcommand{\vphi}{{\bm{\phi}}}

\def\lp{\left(}
\def\rp{\right)}

\DeclareMathOperator*{\argmin}{argmin}

\usepackage{tikz-cd}

\usepackage[capitalize]{cleveref}

\theoremstyle{thmstyleone}%
\theoremstyle{thmstyletwo}%
\newtheorem{remark}{Remark}%

\theoremstyle{thmstylethree}%

\begin{document}

\title[Article Title]{A Modular Deep Learning-based Approach for Diffuse Optical
Tomography Reconstruction}


\author[1]{\fnm{Alessandro} \sur{Benfenati}}\email{alessandro.benfenati@unimi.it}

\author*[2]{\fnm{Paola} \sur{Causin}}\email{paola.causin@unimi.it}

\author[2]{\fnm{Martina} \sur{Quinteri}}\email{martina.quinteri@studenti.unimi.it}

\affil*[1]{\orgdiv{Department of 
  Environmental Science and Policy}, \orgname{University of Milano}, \orgaddress{\street{via Celoria 2}, \city{Milano}, \postcode{20133}, \country{Italy}}}

\affil[2]{\orgdiv{Department of Mathematics}, \orgname{University of Milano}, \orgaddress{\street{via Saldini 50}, \city{Milano}, \postcode{20133}, \country{Italy}}}


\abstract{Medical imaging is nowadays a pillar in 
diagnostics and therapeutic follow--up. Current
research tries to integrate established - but ionizing - 
tomographic techniques with technologies
offering reduced radiation exposure.  
Diffuse Optical Tomography (DOT) uses non--ionizing light in the Near-Infrared (NIR) window to reconstruct 
optical coefficients in living beings,
providing functional indications about 
the composition of the investigated organ/tissue. 
Due to predominant light scattering at NIR wavelengths, DOT reconstruction is, however, a severely ill--conditioned
inverse problem.  
Conventional reconstruction approaches  show severe weaknesses when dealing
also with mildly complex cases and/or are computationally very intensive.
In this work we explore deep learning techniques for DOT inversion. Namely, we propose a fully data--driven approach based on a modularity concept: 
first data and originating signal
are separately processed via autoencoders,
then the corresponding low--dimensional 
latent spaces are connected via
a bridging network which acts
at the same time as a learned regularizer.}




\maketitle

\section{Introduction}\label{sec1}

{D}{iffuse} Optical Tomography (DOT) 
is an emerging technology which employs light in the NIR spectral window to investigate biological tissues in living beings
for diagnostic and monitoring purposes. 
Its use has been explored in different medical fields, including imaging of brain, thyroid, prostate and 
breast cancer screening~\cite{hoshi2016overview,yamada2014diffuse,jiang2018diffuse}. 
Tomographic reconstruction in DOT aims at recovering the spatial distribution of the tissue 
optical properties, which can be related 
to chromophores concentration 
(mainly water, oxy and deoxy-hemoglobin and, in the case of adipose tissues as breast, lipids) and used to monitor functional changes in blood flow. 
 In DOT technologies, light emitted
from external sources at different positions is let to propagate throughout 
the investigated organ and the emerging photon flux is measured
on the tissue boundary~\cite{arridge1999optical}. DOT reconstruction is a notoriously severely ill-posed and ill-conditioned
inverse problem since at NIR wavelengths the outgoing light consists of a mixture of very few  coherent and quasi--coherent photons and a predominant component of incoherent (diffusive) photons
which experienced multiple scattering events while
crossing the tissue.  As such, DOT  
is intrinsically very sensitive to measurement 
noise~\cite{nielsen2006impact}.

\null

 Model-based image reconstruction approaches have represented in the past the standard strategy to perform DOT reconstruction~\cite{hoshi2016overview}. They are based on a mathematical model which describes the physics beyond the generation of the measurable quantity on the boundary, and possibly includes also the effect of the acquisition system and the noise affecting the data. Under suitable hypotheses, the reconstruction of the object of interest is achieved by the minimization of a \emph{fit-to-data} functional, which measures the discrepancy between the recorded data and the corresponding quantity computed by the model. The minimization problem is then solved by an optimization algorithm, for example the Gradient Descent or Alternating Direction Method of Multipliers~\cite{nocedal1999numerical}.

Due to the ill-conditioned nature of the reconstruction problem, DOT solvers must include
a form of regularization, which provides hard or soft prior information on the solution. 
As for this latter category, the $\ell_2$--norm (Tikhonov/ridge regression) penalization is the benchmark classic regularization approach and has been widely used in DOT literature~\cite{Okkyun2011,Huseyin2016}. 
The $\ell_1$-norm (lasso) penalization has been also used in this context to improve
sharp edges detection, obtained as a by--product 
of sparsity enhancement~\cite{Okawa:11}. 
Some of the authors of the present paper proposed in~\cite{Causin19,Causin20} the use of 
an Elastic--Net regularization term.  
This approach was investigated both in 2D and 3D domains and provided 
results superior to the use of pure $\ell_2$  or $\ell_1$ strategies: the (convex) combination of these norms by the Elastic Net indeed 
was able to give a satisfactorily stable solution with a preserved sparsity pattern in selected test cases. 
An alternative approach was proposed in~\cite{Benfenati2020}, based on 
Bregman iterations~\cite{Zhang201120}: this iterative technique consists in substituting the regularization function with its Bregman distance from the previous iterate.  This procedure, which
has been proved to provide a solution to the optimization problem and to possess a remarkable contrast enhancement ability \cite{Benfenati6,Benfenati14}, also
performs satisfactorily in selected cases in the DOT framework~\cite{Benfenati2020}.
Other regularization functionals have also been explored in literature, among which we cite here the combined use of a TV method
with $\ell_1$--norm in the {\tt TOAST++} software for DOT reconstruction~\cite{Schweiger2014toast++}.
Nevertheless, despite the amount of the above research work, poor reconstruction results 
are common outcomes of classical approaches, already in mildly complex situations and/or the reconstruction 
is too computationally intensive 
for routinal clinical use.
Specifically, a main critical point is the fine tuning of regularization parameters required
to obtain a reasonable solution: as a matter of fact, a strong 
dependence exists on the geometry, mesh-size, optical coefficients distribution and level 
of noise of each single case, so that a general, effective, strategy is very
hard to obtain. 

\null

\noindent

Recently, the outstanding
performance  on computer vision tasks of deep learning (DL) algorithms based on Neural Networks (NNs), and especially Convolutionary Neural Networks (CNNs), 
has motivated studies to explore their use also in DOT reconstruction,
with the aim to overcome the above discussed shortcomings. The application of DL techniques 
to DOT is still at its beginnings and very recent: 
in~\cite{Wu20}, the so--called SIMBA strategy was proposed: it uses denoising priors, an online extension of classical \emph{Regularization by Denoising} approaches~\cite{Romano17,Cohen21,cascarano2023constrained}; the authors in~\cite{Yoo2020} trained a NN for inverting in an end--to-end fashion the Lippman--Schwinger equation (a model for light propagation in the tissue), by firstly learning the pseudoinverse of the nonlinear mathematical operator,  and then applying a downstream network to remove the remaining artifacts. In~\cite{mozumder2021model}, the authors combined DL with a model-based approach using Gauss--Newton iterations where the update function was learned via a CNN.  We refer to the recent work~\cite{aspri2023mathematical} for a review of the 
other few applications of DL in DOT reconstruction and, more in general, for an overview of the problem,
both from the theoretical and computational 
viewpoints. 

\null

In this work, 
we exploit the nonlinear approximation 
capabilities of specific deep neural networks, called autoencoders, to perform
DOT reconstruction. We propose a modular
approach based on the idea 
of specializing learning tasks 
for the data and originating signals
and putting them together upon
a pre-training phase. 
Namely, first measured data and generating signal are mapped 
into corresponding feature spaces. Then,
these reduced dimensional spaces are connected via a ``bridge`` neural operator which 
connects the latent spaces and at the same time 
acts as a data--driven regularization
operator. 
This strategy shares conceptual similarities with the learned-SVD 
approach introduced in~\cite{Boink19}.
In the present context
we provide a thorough analysis of the
training strategy which makes 
this approach modular and efficient. 
We also draw a parallel  
with DL-based Reduced Order Modeling techniques for parametrized 
partial differential equations (see, {\em e.g.},~\cite{franco2023deep}), where
a latent space of the full order solution
is connected with the input parameters 
and with the neural-network based inversion strategy for electric impedance tomography presented in~\cite{seo2019learning}.
In addition, we remark the fact that the idea of learning the regularization via NNs is not new in the general field of inverse problem
solution: just to cite a few examples, in~\cite{Lunz18} a regularization term was learned using only unsupervised data;  
in~\cite{Kobler_2020_CVPR, Kobler20} a Total Deep Variation functional was learned, while the NETT framework  of~\cite{Li2020} employs as regularizer the encoder part of an encoder--decoder NN. However, at the best of our knowledge, there
are no applications of these concepts in DOT reconstruction.

\null 

The paper is organized as follows: in Sect.~\ref{sec:phys} we introduce the problem and the classical and DL-based approaches. Namely, first we recall the main physical phenomena occurring in light propagation and we introduce the optical coefficients that are the object of DOT reconstruction, then we present mathematical models of light distribution, we define the DOT inverse problem and we present the discrete setting for DOT reconstruction. In 
Sect.~\ref{sec:DL_DOT} we introduce our DL-based approach, discussing the proposed network architecture, its relations with existing approaches and the training strategy. 
In Sect.~\ref{sec:num} we present extensive numerical experiments and we discuss them; eventually, in 
Sect.~\ref{sec:discu} we draw the conclusions of the work.

\section{Problem Setting}
\label{sec:phys}

\subsection{Optical properties of biological tissues}
\begin{figure}[H]
\begin{center}
\includegraphics[scale=1.3]{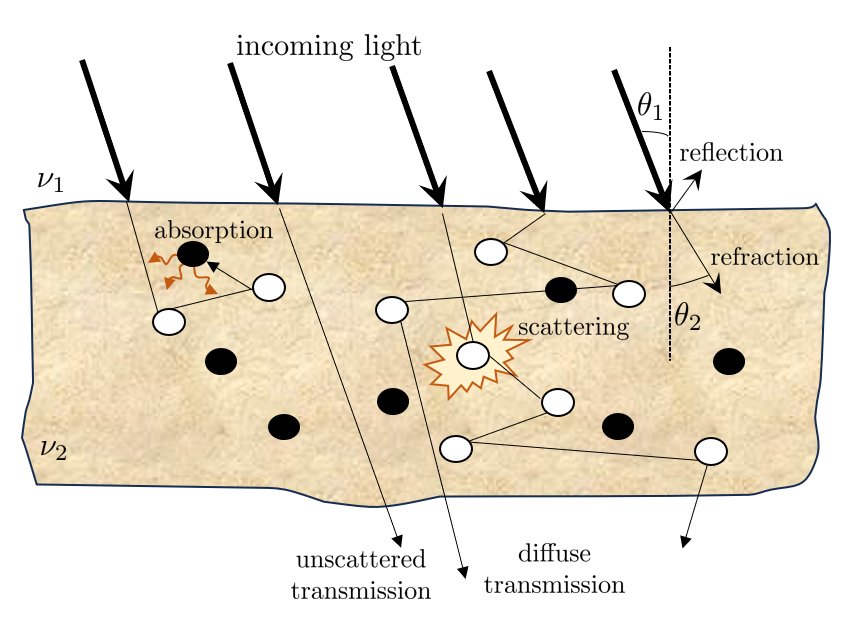}
\caption{Biological tissues are complex structures belonging to the class of turbid media with respect to light propagation. An incident light beam may undergo reflection, refraction, absorption and scattering. Of the penetrating photons, only a limited number of emerge from the sample via unscattered transmission, the most part of them being instead subjected to multi-scattering events.}
\label{fig:tissue}
\end{center}
\end{figure}

The optical properties of biological tissues are the result of their complex internal microstructure, which spans multiple scales and displays a non--homogeneous distribution of the refractive index~\cite{lorenzo2012principles}. 
In  Fig.~\ref{fig:tissue} we schematically depict the main phenomena occurring to a light beam incident on a slice of biological matter, that is, reflection, refraction, absorption and (multi)scattering. Such physical processes can be quantified by a series of parameters, which are presented as coefficients and include the absorption ($\mu_a$) and scattering ($\mu_s$)
coefficients  and the index of refraction ($\nu$) of the medium.
Absorption involves the extraction of energy from the photon by chromophores which undergo a transition process (electronic transition or vibrational state transition). Considering 
a biological medium comprising different chromophores, one can write 
an equivalent absorption coefficient 
at wavelength~$\lambda$, denoted by 
$\mu_a=\mu_a(\lambda)$, as
\begin{equation}\label{eq:mua:multiple}
\mu_{a}(\lambda) = \sum_i \varepsilon_{i}(\lambda) C_i,
\end{equation}
where $\varepsilon_{i}(\lambda)$ and $C_i$ are the absorption extinction coefficient at wavelength~$\lambda$ and concentration of the $i$-th chromophore, respectively. More specifically, since the main absorbing chromophores in biological tissues are oxy--genated or deoxy--genated haemoglobin, and  to a lesser extent water, lipids, melanin, myoglobin, and cytochromes~\cite{wang2012biomedical}, one can rewrite relation~(\ref{eq:mua:multiple}) in the form
\begin{equation}\label{eq:mua:NIR}
\mu_a(\lambda) =B \lp{\rm Sat} \cdot \mu_{a,oxy}(\lambda)+(1-{\rm Sat})\cdot \mu_{a,deoxy}(\lambda)\rp+ W\mu_{a,water}(\lambda) +
\sum_{i_{min}} \varepsilon_{i_{min}}(\lambda) C_{i_{min}},
\end{equation}
where Sat is blood oxygen saturation, $B$ blood volume fraction in the tissue, $W$ water content in the tissue and the sum accounts for other minor absorbers.
Diffuse optical imaging typically
works in the NIR spectrum since 
living tissues do not contain strong intrinsic chromophores that absorb radiation in correspondence to that wavelengths and this allows for deeper penetration.
In the following, we
will assume to collect data on a single
NIR frequency and we will 
write simply $\mu_{a}$ (cm$^{-1}$), thus omitting the
dependence on~$\lambda$. On the contrary, since we are
interested in the spatial distribution 
of $\mu_{a}$, we will assume 
that the distribution of chromophores 
in the tissue in not homogeneous, so that 
the absorption coefficient depends on the spatial position $\vr$, {\em i.e.} 
$\mu_{a}=\mu_{a}(\vr)$.

\null

\begin{figure}[htbp]
\centering
\includegraphics[width=.43\linewidth]{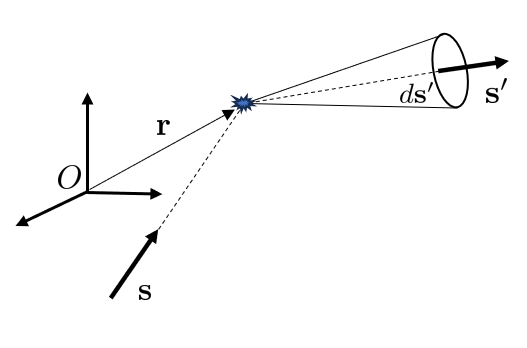}
\caption{Scattering in biological tissues causes a ray of light coming from direction ${\vs}$ to be scattered
by a particle in position $\vr$ in direction ${\vs}^{\, \prime}$ around a cone of angle $d{\vs}^{\, \prime}$.}
\label{fig:scatterevent}
\end{figure}

Scattering originates from the interaction of photons with structural heterogeneities present inside the biological matter, mainly cellular organelles such as mitochondria, thin fibrillar structures of connective tissues, melanin granules and red blood cells, and from the non-uniform distribution of the refractive index in the medium. The scattering interaction between a photon and a molecule results in the photon moving in a different direction.  The NIR window is characterized by the prevalence of elastic scattering, with strong directness upon the collision. Scattering from individual spheres can be calculated from Mie theory~\cite{bohren2008absorption}; for a medium comprising different types of scatterers one has the equivalent scattering coefficient $\mu_s=\mu_s(\lambda)$ given by
\begin{equation}
\mu_s(\lambda)=\sum_{i} \rho_i \sigma_i(\lambda),
\end{equation}
where $\rho_i$ is the number density of spheres of the $i$-th type and $\sigma_i(\lambda)$ the equivalent scattering  cross section
at wavelength $\lambda$. The scattering mean free path $l_s=1/\mu_s$ represents the average distance a photon travels between consecutive scattering events. The angular distribution of the scattered radiation about the outgoing direction~${\vs}^{\, \prime}$ is given by a cone of solid angle~$d {\vs}^{\, \prime}$ and
axis 
${\vs}^{\, \prime}$ originating at the scatterer at position $\vr$ from an incident ray of direction ${\vs}$  (see Fig.~\ref{fig:scatterevent}). 
Again, in the following we will consider a fixed $\lambda$ in the NIR spectrum and we will denote for short the scattering coefficient as $\mu_s$. 
In the present framework, we assume that
$\mu_s$ is a known constant throughout the domain.

\noindent {\em Aim of DOT.} Imaging by DOT aims at recovering the unknown distribution of optical coefficients in the
 investigated biological sample illuminated by NIR light sources. 
Specifically, knowledge of the absorption coefficient field in the tissue can be used to support screening and diagnosis of
 diseases as cancer or stroke, which are known to alter its physiological distribution due to the accumulation of specific 
chromophores~\cite{sandell2011review}. 
Due to relevant technological differences, which directly impact on the reconstruction problem, we specifically focus on DOT 
based on  Steady-State-Domain systems, which are among the most widely used in clinical settings for breast cancer screening,
 a main field of application of DOT. In this technology, the light is emitted into the tissue at a single frequency (continuous wave modality, CW)
  and the transmitted luminous signal is measured directly on the tissue boundary or imaged via a CCD camera. 

\subsection{Mathematical models of light propagation}
\label{sec:models}

The mathematical description of the process of light propagation in tissues can be performed via the Maxwell equations or using the theory of 
radiative transfer (RT) of energy. Models stemming from this latter theory are the most used in optical imaging modalities, since they avoid to 
deal with the complex statistical nature of the non--homogeneous medium required for the solution of Maxwell equations.

\noindent  We let~$\Omega\subset \mathbb{R}^d, d=2$ or $3$, be a connected open domain
representing the tissue to be investigated, and $S^{d-1}$ the unit sphere in $\mathbb{R}^d$. The main equation of RT reads (at steady--state)~\cite{peraiah2002introduction} 
\begin{equation}
\lp\,{\vs} \cdot 
\nabla +(\mu_a(\vr)+\mu_s)\,\rp I(\vr, {\vs})
+ \mu_s
\int_{S^{d-1}} I(\vr, {\vs}^{\, \prime})
p({\vs}, {\vs}^{\, \prime}) d {\vs}^{\, \prime}=f(\vr), \qquad \vr \in \Omega,
\label{eq:RT}
\end{equation}
where $I=I(\vr, {\vs})$ (W cm$^{-2}$sr$^{-1}$) is the radiance at point $\vr$ in direction ${\vs}$,  $f(\vr)$ an internal source term 
and the scattering phase function $p({\vs}, {\vs}^{\, \prime})$ represents the probability density function of photon scattering from 
direction ${\vs}$ to direction ${\vs}^{\, \prime}$. It is common assumption to let the scattering phase function depend only on the angle 
between the incident and scattered directions, hence to consider it as a function of only the scalar product $({\vs} \cdot {\vs}^{\, \prime})$. In this context, a constant of interest is the cosine-weighted average of the scattering 
\begin{equation}
g=\int_{S^{d-1}}p({\vs}, {\vs}^{\, \prime}) \,
{\vs}\cdot {\vs}^{\, \prime}
\, d{\vs}^{\, \prime},
\qquad g \in [-1,1].  
\label{eq:g}
\end{equation}
The above quantity is a measure of the forward scattering bias and quantifies how efficiently photons keep propagating in the forward 
direction despite scatter. When $g = 0$ one has fully isotropic scattering, while when $g \rightarrow 1$ the propagation is strongly 
forward biased, that is scattering vanishes, and for $g \rightarrow -1$ scattering becomes completely backward directed. For (in vitro) breast tissues
 at the visible and NIR wavelengths one has typically $0.65 < g < 0.95$~\cite{jacques87b}. 

\null

\noindent The RT model~(\ref{eq:RT}) is an integro--differential equation and its solution is computationally very  expensive. 
A physically reasonable and cost--effective model can be derived by performing an expansion of the RT equation in spherical harmonics 
(see, {\em e.g.},~\cite{Durduran2010DiffuseOF} for details and derivation). The following PDE model, known as the diffusion approximation (DA) model, is obtained truncating the expansion at first order:
\begin{equation}\label{eq:eq_diff}
 \lp-\nabla\cdot (D(\vr) \nabla) 
+ \mu_a(\vr)\rp U(\vr) 
 = f(\vr), \qquad \vr \in \Omega,
\end{equation}
where the field $U(\vr)=\int_{S^{d-1}}I(\vr, {\vs}^{\, \prime}) d {\vs}^{\, \prime}$ is the photon fluence due to the light source $f(\vr)$ and where $D=D(\vr)$ is the diffusion coefficient given by
\begin{equation}
D(\vr) = \frac{1}{3(\mu_a(\vr)+(1-g)\mu_s)}, 
\end{equation}
with $g$  defined by~\eqref{eq:g}. Model~\eqref{eq:eq_diff} is equipped with the Robin--type boundary condition 
\begin{equation}
\label{eq:BCs}
    U(\vr)+\frac{\zeta D(\vr)}{2c_d} \nabla U(\vr) \cdot {\vn} =0,  \qquad \vr \in \partial \Omega\\
\end{equation}
where $\zeta=(1+R)/(1-R)$, $R$ being a function of the refraction indices $\nu_1$ and $\nu_2$ of the two media concurring at the interface
 according to the Fresnel law~\cite{Mansuripur2002}, and $c_d$ is a constant depending on the space dimension ($c_d = 1/\pi$ if $d = 2$, $c_d = 1/4$ if $d = 3$).

\medskip

\noindent {\em The Rytov perturbation approach.}
A common approximation implemented in reconstruction software embedded in CW-DOT instrumentation is the Rytov linearization 
of the DA model.  A complete derivation of this methodology is provided 
in the general field of line-of-sight wave propagation through random media in Ishimaru~\cite[Vol. II, Ch. 17]{Ishimaru78} and for the specific  DOT application in~\cite{causin2020inverse}. In brief, one assumes the linearization $\mu_a(\vr)=\mu_{a,0}(\vr) + \delta \mu_a(\vr)$, where $0 < \mu_{a,0} \le \overline{\mu}_a $ is the  background value corresponding
to a condition with absence of regions with altered absorption and $\delta \mu_a$ a (small) perturbation term corresponding to a condition with the presence
of localized contrast regions with altered
- typically increased - absorption coefficient, the realization of the perturbed status depending on the specific technological setting.
An exponential change of variables is also introduced such that  
\begin{equation}
U(\vr)=U_0(\vr){e^{\psi_1(\vr)}},
\label{eq:expchange}
\end{equation}
where  $U_0(\vr)$ is the light fluence field in background conditions and $\psi_1=\log(U/U_0)$ is the so--called {\em logarithmic amplitude fluctuation} of the light due to the presence of contrast regions. Inserting~\eqref{eq:expchange} into~\eqref{eq:eq_diff}, performing algebraic manipulations and neglecting higher order terms, yields the following equation for  the combined quantity $(U_0 \psi_1)$ 
\begin{equation}
\begin{array}{ll}
\displaystyle\left[\Delta - \frac{\mu_{a,0}(\vr)}{D}\right](U_0 \psi _1)(\vr)=\frac{\delta \mu_a(\vr)}{D} U_0(\vr),& \quad \vr \in \Omega.
\end{array}
\label{eq:MH}
\end{equation} 
Notice that relation~(\ref{eq:MH}) is obtained under the
further assumption of
constant diffusion coefficient given by~$D=1/(3\mu_s)$,
a simplification which is quite common
in the CW-DOT context (see,~{\em e.g.},
\cite{boas2001simultaneous,konecky2008imaging}).

\subsection{DOT inverse problem}
\label{sec:inverse}
We let $Q$ be the function space of optical absorption coefficient given by
\begin{equation}\label{param:RTE:Q}
Q(\vr) = \left\{q(\vr)=\mu_a(\vr)\,:\, 0 < \mu_a(\vr) \le \overline{\mu}_a\right\},
\end{equation}
$\overline{\mu}_a$  being a biophysically plausible positive upper bound for the 
absorption coefficient. Moreover, we let 
$\mathbb{W}$ represent a (possibly infinite) set of source functions in~$\Omega$ and  $W$
represent the space spanned by the functions in $\mathbb{W}$. Then, for a particular applied light source
$f \in \mathbb{W}$, we denote by 
$\phi_{q}=\phi_{q}(\vr; f)$ the state variable
belonging to the vector space $\Phi_q$.
The field $\phi_{q}$  
corresponds to the light distribution field in the body due to the light source~$f$,
with optical absorption~$q$ and is solution of one of the above PDE-type models, of the generic form    
\begin{equation}
\begin{array}{ll}
\mathcal{L}\,\phi_q(\vr; f) = f(\vr), & \vr \in \Omega, \\[1mm]
\mathcal{B}\,\phi_q(\vr; f) = 0, & \vr \in \partial \Omega,    
\end{array}
\label{eq:PDE}
\end{equation}
where $\mathcal{L} : \Phi_q \times Q \rightarrow W$  is an (integro--) differential operator augmented by a set of boundary conditions  
represented by the operator $\mathcal{B}$. Hereafter, we drop the dependence on $\vr$ and $f$ to ease the readability. The DOT inverse 
problem consists in the identification of the set of parameters $q$, subject to Eq.~(\ref{eq:PDE}), which allow to obtain the best fit with 
respect to a set of observations of the state variable obtained on (a portion of) the domain boundary. In practice, one disposes of the  
noisy measurement vector $y_\delta \in Y$, $Y$ being (a possibly
infinite) vector space, for some noise level $\delta \ge  0$, obtained via the measurement map 
 $M: \Phi_q \times W  \rightarrow Y$ such that 
\begin{equation}
y_\delta = M(\phi_q; f)   
\end{equation}
represents the measured light field 
at the observation points. 
In addition, under suitable assumptions on $\mathcal{L}$, one can introduce the parameter-to-state map 
$S :  Q \rightarrow  \Phi_q$, i.e. $S(q) = \phi_q$, which is a solution to the boundary value problem~(\ref{eq:PDE}),
and define the {\em forward map} $\mathcal{F}= P \circ S: Q \rightarrow Y$, where $P$ is the observation operator 
providing the value of the computed solution field at the observations points. Then one aims at computing the best parameter set~$q^*$ such that 
\begin{equation}
q^* \in \argmin_{q \in Q}\,\mathcal{C}(\mathcal{F}(q),y_\delta), 
\label{eq:opt1}
\end{equation}
where $\mathcal{C}: Y \rightarrow [0,+\infty)$ is a
cost function. One may take $\mathcal{C}(\mathcal{F}(q),y_\delta)=||\mathcal{F}(q)-y_\delta||_{\ell_2}^2$ but also
other measures as the Kullback-Leibler divergence are possible
choices.
Observe that, when Rytov linearization is adopted, upon solving Eq.~(\ref{eq:MH}), one compares in~(\ref{eq:opt1}) 
the quantity $\psi_1$ evaluated at the measurement points with the ratio 
$\log(y_\delta/y_{0,\delta})$ at the same locations,
where light measurements $y_{0,\delta}$ and $y_{\delta}$ correspond to
baseline and perturbed
conditions, respectively.

As already mentioned, due to predominant scattering phenomena in light transmission, DOT reconstruction
(\ref{eq:opt1})
 is a severely ill-posed problem. Mathematically, this means that there is lack of continuous dependence of the solution 
 on the measurements (small errors in the data may lead to large errors in the model parameters).
 In these circumstances, it becomes essential the adoption of regularization techniques to solve the problem. One possible way is
to introduce a regularization term such that the problem becomes: find the optimal 
parameter set $q^*$ 
such that 
\begin{equation}\label{eq:opt_prob2}
q^*\in\argmin_{q\in Q}\,[\mathcal{C}(\mathcal{F}(q),y_\delta) 
+ \mathcal{R}(q)],
\end{equation}
where $\mathcal{R}: Q \rightarrow [0,+\infty]$ is a regularization functional 
which encodes soft or hard 
prior knowledge about the hidden structures to be investigated, favoring appropriate minimizers or penalizing those with undesired structures \cite{ArrMasOekSch19,benning_burger_2018,engl1996regularization}.

\subsection{Discrete DOT inverse problem}
\label{sec:discrete}
DOT reconstruction is inherently performed at the discrete level, since in practice one disposes only of a finite set 
of $n_d$ (possibly noisy) measures of light at the observation points (detectors), obtained when
each of the $n_s$ light sources is turned on, for a total of 
$M$ measures. We denote by 
${\bf y}^M_\delta \in \mathbb{R}^{M}$ such set. 
 In addition, also
 the absorption coefficient is reconstructed in a finite set of locations 
in the computational domain, typically
corresponding to a discretization of the domain
into $V$ voxels, with diameter size~$h$. The continuous 
set of optical coefficients is thus replaced by 
a finite dimensional set, whose elements
${\boldsymbol \mu}_{a}^h \in \mathbb{R}^{V}$
are the approximation of the absorption coefficient value in each voxel.

\subsubsection{PDE-based solvers}

Classic solvers are based on the solution of the optimization problem~(\ref{eq:opt1})
(or~(\ref{eq:opt_prob2})). This requires
to dispose of: {\em i)} 
a computationally efficient method to deal with the PDE system~(\ref{eq:PDE}) which provides the light field for a given set of optical parameters (forward problem); {\em ii)} an optimization procedure.

As for {\em i)}, the standard approach is to decompose the computational domain into triangular/trapezoidal areas (2D) or into cuboids of suitable shape (3D). Upon introducing a discretization of the domain, one then operates with functions belonging to finite dimensional subspaces parametrized by the mesh diameter~$h$. For example, using a finite elements method, one approximates the solution $\phi_q$ by a piecewise polynomial function $\phi_{q,h} \in \Phi_h$ given by
\begin{equation}
\phi_{q,h}=\sum_{j=1}^{N_h}\phi_j \, \xi_j(\vr),
\end{equation}
where $\Phi_h$ is a finite dimensional subspace 
(typically chosen such that $\Phi_{q,h}\subset \Phi_{q}$) spanned by the basis functions $\{\xi_j\}_{j =1,\dots,N_h}$ and where the nodal values $\phi_j$ are to be determined. Piecewise polynomial approximation to the optical
absorption coefficient, ${\boldsymbol \mu}_{a}^h$, is constructed in the same way. Applying a Galerkin approach transforms the continuous problem into the $N_h$-dimensional discrete problem of finding the nodal field values at all nodes, given the set of nodal parameters.

As for {\em ii)}, supposing that the forward map $\mathcal{F}$ is Fréchet differentiable, a Taylor expansion around a reference point $q_0$ reads
\begin{equation}
\mathcal{F}(q) \approx
\mathcal{F}(q_0)+
\mathcal{F}^\prime(q_0)\delta q,\,\, \delta q = q-q_0
\label{eq:linearizedF}
\end{equation}
where $\mathcal{F}^\prime$ is the first order Fr\'echet derivative of the forward operator with respect to $q$. Inserting~(\ref{eq:linearizedF}) into (\ref{eq:opt1}) and truncating to the first order yields, upon discretization of the fields and using
directly the increment~$\delta \vmu_a^h$ 
\begin{equation}
({\delta \vmu_a^h})^*\in\argmin_{\delta \vmu_a^h} [\mathcal{C}({\vJ}  \delta \vmu_a^h,  \vy_{\delta}^M)+\mathcal{R}(\delta \vmu_a^h)],
\label{eq:Ryt} 
\end{equation}
where $\vJ$ of size $M \times V$
denotes in the discrete setting the first order Fr\'echet derivative (Jacobian matrix) and where 
$\delta \vmu_a^h$ is the voxel-based
increment for the unknown optical field. 
The optimization problem (\ref{eq:Ryt})  can
be solved by one of the many available algorithms (see~\cite{aspri2023mathematical} for a review).

\section{Learning DOT reconstruction via Neural Networks}
\label{sec:DL_DOT}
We discuss in this section our original contribution. 
We aim to find an approximate  
solution to DOT reconstruction 
leveraging the approximation capabilities of
nonlinear neural networks. Namely, 
we introduce neural networks which produce mappings
of the type
\begin{equation}
\mathcal{N}: \mathbb{R}^{M} \rightarrow 
\mathbb{R}^{V}
\end{equation}
and perform data-driven reconstruction
of the optical coefficient from
observations. 
The architecture we propose is characterized
by the two autoencoder networks
\begin{equation}
\begin{array}{cc}
 \mbox{data-AE}:    & \Psi_y:\mathbb{R}^{M} \rightarrow 
\mathbb{R}^{M}, 
 \\[3mm]
\mbox{signal-AE}:     & 
\Psi_\mu:\mathbb{R}^{V} \rightarrow 
\mathbb{R}^{V}.
\end{array}
\end{equation}
The data-AE acts on observations ${\bf y}^M_\delta$
and is defined by $\Psi_y=\mathcal{D}_y \circ \mathcal{E}_y$, with 
\begin{equation}
\begin{array}{cc}
{\bm z}_y=\mathcal{E}_y({\bm y}^M_\delta)
     & 
\mathcal{E}_y:   \mathbb{R}^{M} 
\rightarrow \mathbb{R}^{n_y} 
\qquad (\text{encoder})  \\[3mm]
\tilde{{\bm y}}^M_\delta=\mathcal{D}_y({\bm z}_y),  & 
\mathcal{D}_\mu: \mathbb{R}^{n_y} \rightarrow  \mathbb{R}^{M}
\qquad (\text{decoder}),
\end{array}
\end{equation}
where ${\bf z}_y$ is the encoded feature (latent or compressed representation) of the data and $\tilde{\bf y}^M$ its decoded   
version. 
The signal-AE acts on the generating signal 
${\bm y}^M_\delta$
and is defined by $\Psi_\mu=\mathcal{D}_\mu \circ \mathcal{E}_\mu$, with 
\begin{equation}
\begin{array}{cc}
{\bm z}_\mu=\mathcal{E}_\mu({\boldsymbol \mu}_{a}^h)
     & 
\mathcal{E}_\mu: \mathbb{R}^{N_h} \rightarrow \mathbb{R}^{n_\mu} 
 \qquad (\text{encoder}) \\[3mm]
\tilde{{\boldsymbol \mu}}_{a}^h=\mathcal{D}_\mu({\bm z}_\mu),  & \mathcal{D}_\mu: \mathbb{R}^{n_\mu} \rightarrow  \mathbb{R}^{N_h} \qquad (\text{decoder}),
\end{array}
\end{equation}
where ${\bm z}_\mu$ is the encoded generating
signal and $\tilde{{\boldsymbol \mu}}_{a}^h$ its decoded
version. 
In addition, we introduce a third network,
called bridge-NN, which has the role to connect (as a bridge)
the two previous latent spaces 
\begin{equation}
\Sigma: \mathbb{R}^{n_y} \rightarrow \mathbb{R}^{n_\mu}. 
\end{equation}
The proposed DL-based DOT reconstruction is thus defined by the following 
structure (called Mod-DOT):
\begin{equation}
\mathcal{N}:= \mathcal{D}_\mu  \circ \Sigma 
\circ  \mathcal{E}_y : \mathbb{R}^{M} \rightarrow \mathbb{R}^{V}.
\end{equation}

\begin{remark}
The proposed architectures is articulated into three modules which serve to different purposes: 
the two AEs are deputed to learn the intrinsic characteristics
of the data and signal manifolds,
respectively. Their specificity 
is tailored to handle the different richness of the respective inputs. 
The bridge needs instead to learn
 the effect-to-cause relation between observation and generating signal in the reduced spaces.   
 A schematic representation of the proposed architecture is given in Fig.~\ref{fig:architect}. 
\begin{figure}[htbp]
\centering
\includegraphics[width=.7\textwidth]{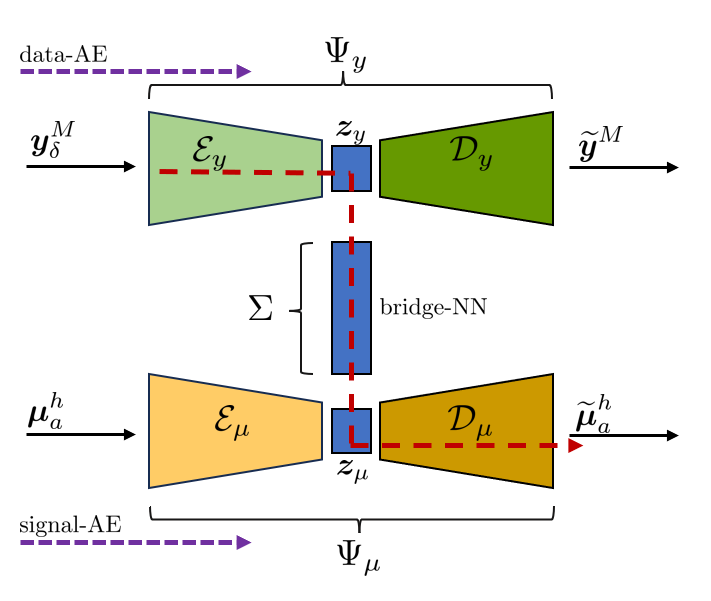}
\caption{Workflow in the Mod-DOT approach. 
Observed data and originating signal are encoded individually 
in the corresponding latent spaces via the data-AE and
the signal-AE, respectively (violet arrow paths). Latent spaces are then
connected via the bridge-NN.
At test time,  
the reconstruction process consists in the cascade application of the three networks (red arrow path)
$\mathcal{N}:= \mathcal{D}_\mu  \circ \Sigma 
\circ \mathcal{E}_y$. 
}
\label{fig:architect}
\end{figure}
    \end{remark}

\noindent{\em The Mod-DOT neural architecture.}
We detail here the AE and bridge architectures used    
for compressing the data and originating signal
manifolds and to approximate the inversion operator. 
\begin{itemize}
\item {\em data-AE}:
predominant scatter in the tissue 
severely hinders the collection of data
that show neat differences even 
for well distinct internal distributions
of the absorption coefficient. 
This results into very similar data patterns
generated by different samples. 
For this reason, this autoencoder is deliberately
chosen to be very simple. In this architecture, 
the smoothing effect typical of autoencoders provides an initial denoising of the measures.   
Both encoder and decoder part are chosen
to be made of fully connected layers coupled with 
{\tt tanh} and {sigmoid} activation functions, respectively;  
\item {\em signal-AE}:
generating signal samples are provided 
under the form of images which encode
the absorption coefficient field. In the images, 
contrast areas of altered (increased) absorption are embedded
into a background value. 
The richness of this input is 
definitely higher than the 
the observations. 
The resulting architectures are schematized in Fig.~\ref{fig:AEs},
where the inputs and outputs are also represented. For the data-AE, the input
is represented by the DOT ``sinogram",
where observations collected 
at the detectors locations for each
light source are represented as a vector. For the signal-AE, the input
is represented by the absorption coefficient field in the form
of an image or a vector. 
In the figure, we represent the case where the net is presented
with a one--channel image in which the dark blue region
represents the background value, while
the two circular regions are contrast
areas with increased absorption coefficient (see also Sect.~\ref{sec:synth} for more 
detailed description).  

\begin{figure}[htbp]
\centering
\includegraphics[width=.95\textwidth]{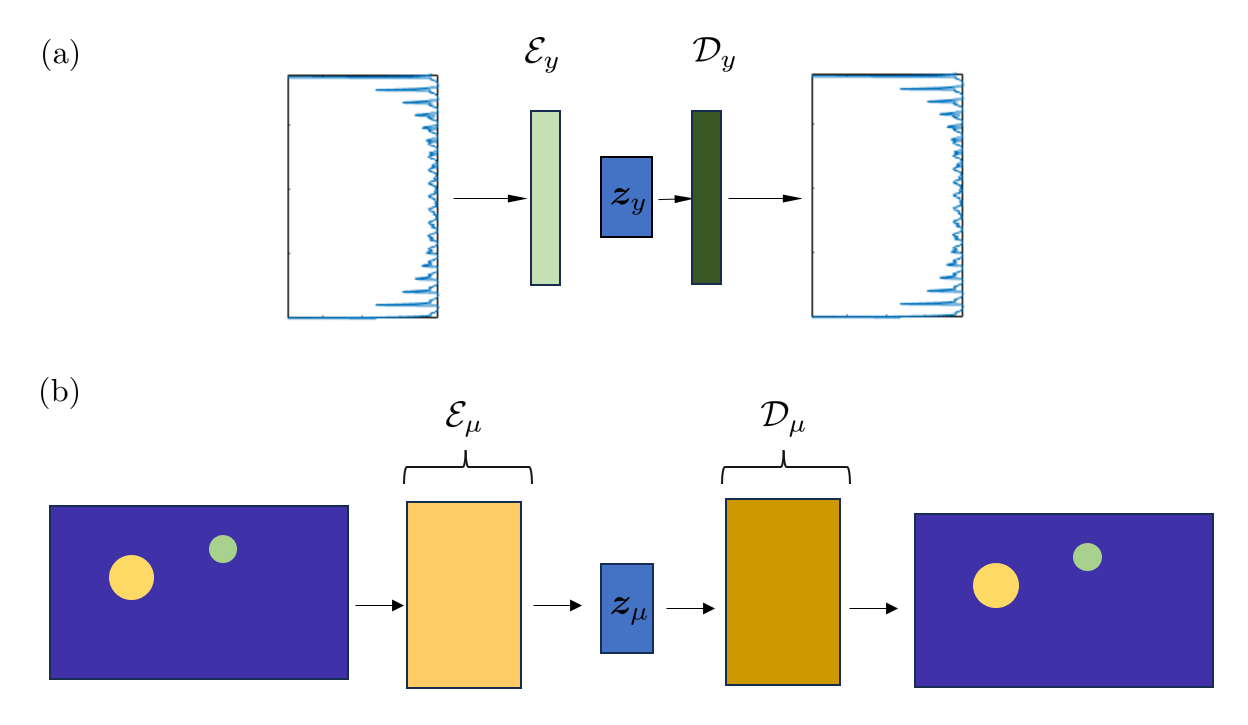}
\caption{Architecture of the AEs: (a) the data-AE receives in input the DOT ``sinogram`` 
composed by the collection of all the light observations at the detectors locations for each source and it is trained to reproduce the identity operator on it; (b) the signal-AE
receives in input the absorption coefficient field, presented here as a one--channel image where the circular areas represent contrast
regions with increased coefficient. 
It is trained to reproduce the identity
operator on the image.}
\label{fig:AEs}
\end{figure}

\item {\em bridge-NN}:
the bridge networks connects the 
data and signal latent spaces,
performing the core of the inversion. 
Fully connected dense layers are adopted in order to effectively 
perform this operation.
\end{itemize}

\null

\noindent {\em Related data-driven approaches.}
The full bridged structure finds a remarkable
interpretation as a learned SVD strategy  
(see~\cite{Boink19} for its original introduction). 
To understand this concept, it is useful to draw a parallel
with the physical models introduced in  Sect.~\ref{sec:models}: 
one may decompose 
the linearized discrete 
parameter-to--state operator as ${\vJ}=\vU{\vS}\vV^\top$, where ${\vS}$ is the diagonal matrix containing the singular values and $\vU$ and $\vV$ are the left and right singular vectors,
respectively, and express the solution via the inversion formula ${ \vJ}^{-1}=\vV {\vS}^{-1} \vU^\top$. 
This fundamental idea  of linear algebra is extended via neural
networks to more general nonlinear parameter-to-state operators,
which are not formally derived from a model but directly learnt from data.   
Namely, one  learns a generalized SVD 
decomposition via the AE architectures.  
In this framework, 
the encoder $\mathcal{E}_y$ can be seen as the corresponding of $\vU^\top$ and the decoder 
$\mathcal{D}_\mu$
can be considered as the counterpart of $\vV$. 
The bridging operator $\Sigma$, which relates the
latent spaces in the data and signal domains, plays the mapping and scaling role of the singular values in~$\vS$ in the linear algebraic SVD approach. 
 This parallelism is visually represented by
 the following scheme:
\begin{equation}
\label{eq:lSVDstruct}
\begin{tikzcd}
\vy^M_\delta \arrow{r}{\mathcal{E}_y}[swap]{\vU^\top} &  \vz_y \arrow{d}{\Sigma}[swap]{\vS^{-1}}\arrow{r}
{\mathcal{D}_y}[swap]{\vU}   &\tilde{\vy}^M\\
\vmu_a^h \arrow{r}{\mathcal{E}_\mu}[swap]{\vV^\top} &  \vz_\mu\arrow{r}{\mathcal{D}_\mu}[swap]{\vV} &\tilde{\vmu}_a^h\\
\end{tikzcd}
\end{equation}

Noticeably, the bridging operator  also acts as a regularizer of the inversion. To understand this point, it is useful to recall 
the classic inversion by Tikhonov regularization, which is 
obtained from~(\ref{eq:opt_prob2}) when $\mathcal{C}$ is the least
square functional and $\mathcal{R}$ is the $\ell_2$ norm of its argument (bound 
on the variance of the solution).
Upon linearization and discretization, one obtains a problem of the type 
\begin{equation}
(\delta {\bm q}^{h})^* \in 
\argmin_{\delta {\bm q}^h \in Q^h} ||{\bm b} - \vJ \delta
{\bm q}^h||_{\ell_2}^2 + 
\alpha ||\delta{\bf q}^h||_{\ell_2}^2={\vV}
\underbrace{({\vS}^2 + \alpha {\rm Id})^{-1} {\vS}}_{\vS_\alpha^{-1}} {\vU}^T {\bm b},
\end{equation}
where ${\bm b}$ stands for the right hand side and $\alpha$ is a regularization 
parameter. The diagonal elements of $\vS_\alpha^{-1}$, defined as 
$\frac{\sigma}{\sigma^2+\alpha}$, filter out the smallest singular values,
according to the threshold $\alpha$. This can be interpreted as a smoothed Truncated SVD, where the singular values go towards zero in a smooth fashion, 
which mitigates noise in the data, enhanced by the smallest singular values at the denominator. Truncated SVD, instead, puts to zero the singular values  after $k$ terms.
Drawing the parallel, the nonlinear SVD obtained
by the neural network learns from the data the operator $\Sigma$ and automatically  doses 
the amount of required regularization. In doing so, it learns at the same time the regularization functional $\mathcal{R}$ and the parameter $\alpha$, which are both encompassed in~$\Sigma$.

\null

\begin{figure}[htbp]
\centering
\includegraphics[width=.7\textwidth]{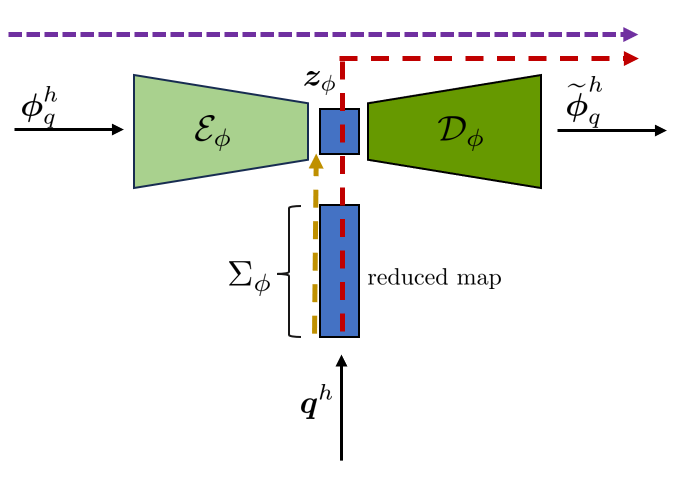}
\caption{Workflow in the DL-ROM approach (here shown for comparison with the present approach). First, the 
AE $\mathcal{D}_\phi \circ \mathcal{E}_\phi$ is trained to learn
an approximation of the identity operator over the solution manifold (violet arrow) using
FOM snapshots obtained by a computationally
intensive computation. 
This AE provides the latent low--dimensional representation ${\bm z}_\phi$ of the solution manifold. Then, the
 network $\Sigma_\phi$ is trained to learn
the map ${\bm q}^h \rightarrow {\bm z}_\phi$ (dark yellow  arrow).
Finally, the composition $\mathcal{D}_\phi \circ \Sigma_\phi$  
defines the DL-ROM approximation of the 
parameter-to-state map at test time 
(red arrow).
}
\label{fig:DL-ROM}
\end{figure}

A second remarkable tie of the present
architecture can be found with  Reduced Order Modeling (ROM) techniques applied
to the solution of parametric partial differential equations of 
type~(\ref{eq:PDE}).  
Conventional ROM techniques
rely on the assumption that the reduced order
approximation can be expressed by a linear combination of basis functions, built from
a set of full order solutions. Namely, first one computes for 
different values of the parameters the so--called Full Order Method (FOM) snapshots.  
Then, for example considering the Principal Orthogonal Decomposition (POD) approach, 
one writes the FOM snapshots 
as columns of a matrix ${\mathbf A}:=[\vphi^h_{q_1},
\vphi^h_{q_2}, \dots, \vphi^h_{q_N}]$ 
and computes its SVD decomposition 
${\bm A}=\vU{\vS}\vV^\top$. The
low--dimensional representation of the solution
is taken as 
the projection 
${\vphi}^{LR}_q:=\vV^\top {\vphi}^h_{q}$
and then one needs an algorithm that approximates
the correspondence ${\bm q}^h \rightarrow 
{\vphi}^{LR}_q$, usually via 
a reduced map $\theta$ such that 
${\vphi}^{LR}_q \approx \theta({\bm q}^h)$.
This map has been chosen in literature in different ways, including
Gaussian process regression and polynomial chaos expansion (see~\cite{franco2023deep} and references therein).
Eventually, one reconstructs  
$\vphi_q^h \approx V \theta({\bm q}^h)$.
The resulting ROM manifold is thus linear and its dimension can become extremely high and/or produce inaccurate results for problems involving nonlinearities. 
Recently, approaches have emerged 
where the solution map is approximated by
a deep neural network, yielding the so--called DL-ROM
methods~\cite{franco2023deep}. 
Adapting the notation of this latter paper to the one used in the present work for
ease of readability, a typical DL-ROM is inspired by the steps of the classical ROM and it
encompasses the three networks 
$\mathcal{E}_\phi, \mathcal{D}_\phi$ and
$\Sigma_\phi$
such that 
\begin{equation}
\mathcal{E}_\phi(\vphi_q^h))={\bm z}_\phi, \qquad 
\mathcal{D}_\phi({\bm z}_\phi)=
\tilde{\vphi}_q^h \approx \vphi_q^h,
    \qquad \Sigma_\phi({\bm q}^h) \approx 
    {\bm z}_\phi,
\label{eq:DL-ROM}
\end{equation}
where $\mathcal{E}_\phi$ is an encoder which
produces the latent space ${\bm z}_\phi$ - conceptually corresponding to
the low resolution representation ${\phi}^{LR}_q$ -, $\mathcal{D}_\phi$ decodes
back the latent space into the solution manifold
and
$\Sigma_\phi$ is a dense network which produces the parameter-to-solution 
operator. Eventually, one has $\mathcal{N}_\phi:=
\mathcal{D}_\phi \circ \Sigma_\phi: 
Q_h \rightarrow \Phi_h$ such that
 $\tilde{\vphi}_q^h =\mathcal{D}_\phi(\Sigma_\phi(q_h))
 \approx \vphi_q^h $.
We schematically depict in Fig.~\ref{fig:DL-ROM} the corresponding
architecture.  It is thus apparent that both the Mod-DOT and DL-ROM
share similarities in the use of specialized modules
connecting low dimensional spaces. 
Observe, however, that the DL-ROM
aims to solve  - in the form 
presented in~\cite{franco2023deep} - the {\em forward} problem
for values of the parameters in correspondence
of which a FOM
solution is not available due to
the high computational cost needed to explore a wide range of parameters. 
As such,
the direction of the workflow of the DL-ROM at test time
is the opposite with respect to the one for the inverse problem we consider in our application (compare for reference the red arrows in Fig.~\ref{fig:architect}
and Fig.~\ref{fig:DL-ROM}). We are not aware
of the use of DL-ROM techniques in our context. 

\null

\noindent A further connection of the present architecture can be found with 
the strategy proposed in~\cite{seo2019learning}
for the reconstruction of the conductivity
field from voltage measurements in
electric impedance tomography (EIT),
another notoriously ill--posed 
inverse problem. As in our approach,
instead of using a net which directly reconstructs
in an end-to-end fashion the 
conductivity field from the voltages, first 
a AE is
trained to learn the identity operator
on the originating signal 
(conductivity presented
as images) and then 
a nonlinear regression map is learned to
connect the measured 
obtained from the AE. This latter
net is expected to deal with 
a significantly better conditioned problem. 
Adapting also in this case the notation of~\cite{seo2019learning} to the present
one, we have 
\begin{equation}
\mathcal{E}_q({\bm q}^h)={\bm z}_q, \qquad 
\mathcal{D}_q({\bm z}_q)=
\tilde{\bm q}^h \approx {\bm q}^h,
    \qquad \Sigma_q(\vphi_q^h) \approx 
    {\bm z}_q,
\label{eq:EIT}
\end{equation}
where $\mathcal{E}_q$ is an encoder which produces the latent space ${\bm z}_q$ from the signal manifold, $\mathcal{D}_q$ decodes
back ${\bm z}_q$ into the signal manifold 
and
$\Sigma_q$ is a regression map which approximates the measure-to-parameter 
operator.
In this case, one has $\mathcal{N}_q:=
\mathcal{D}_q \circ \Sigma_q$.  
Theoretical support of the usefulness
to work with low--dimensional non--linear
manifolds for the signal images is provided
in~\cite{seo2019learning} 
based on the result that the inverse problem in EIT
with sufficiently many measures is uniquely 
solvable and Lipschitz stable on finite dimensional linear subsets of piecewise--analytic
functions. This motivates the idea that
realistic originating signal fields lie on
a non--linear manifold that is lower dimensional
than the space of all possible images. 

 \section{Numerical experiments}
\label{sec:num}

In the following we assess the performance
of the proposed Mod-DOT approach. For comparison, 
we also consider classic variational approaches and 
an End-to-End NN, called E2E-DOT, which consists in the simple juxtaposition of the data encoder and the signal
decoder, without bridge, 
that is
$\mathcal{N}_{E2E}= \mathcal{D}_\mu  \circ 
\mathcal{E}_y$, and mimics end-to-end
approaches already present in literature to address DOT reconstruction. 

\subsection{Generation of synthetic data}
\label{sec:synth}
\begin{figure}[htbp]
\centering
\includegraphics[width=0.95\textwidth]{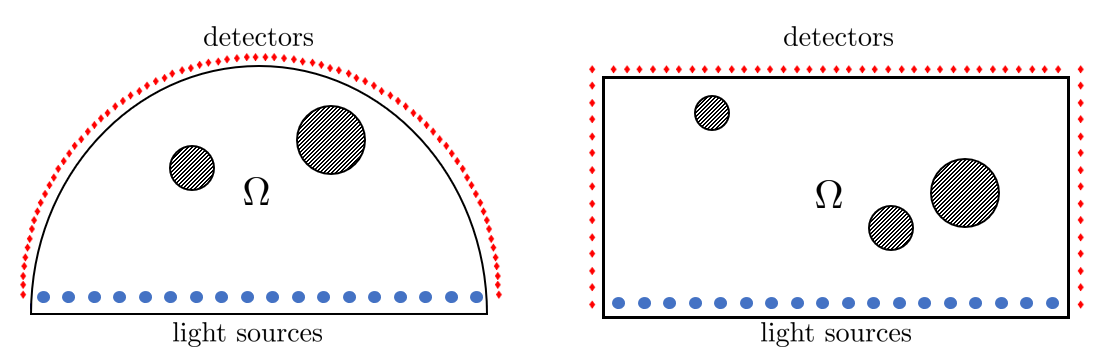}
\caption{Geometry of the semi-disk and rectangular domains considered in the numerical experiments. The detectors are uniformly disposed
	along the curved boundary for the semi-disk and top/left/right
 boundary for the rectangle domain; the sources are aligned along the bottom boundary and
	positioned 1mm inside the domain.
	The absorption 
	coefficient is $\mu_{a,0}$ except for the contrast regions (circular areas) which have
	increased absorption coefficient.}
\label{fig:domains}
\end{figure}

\begin{figure}[htbp]
\centering
\includegraphics[width=\textwidth]{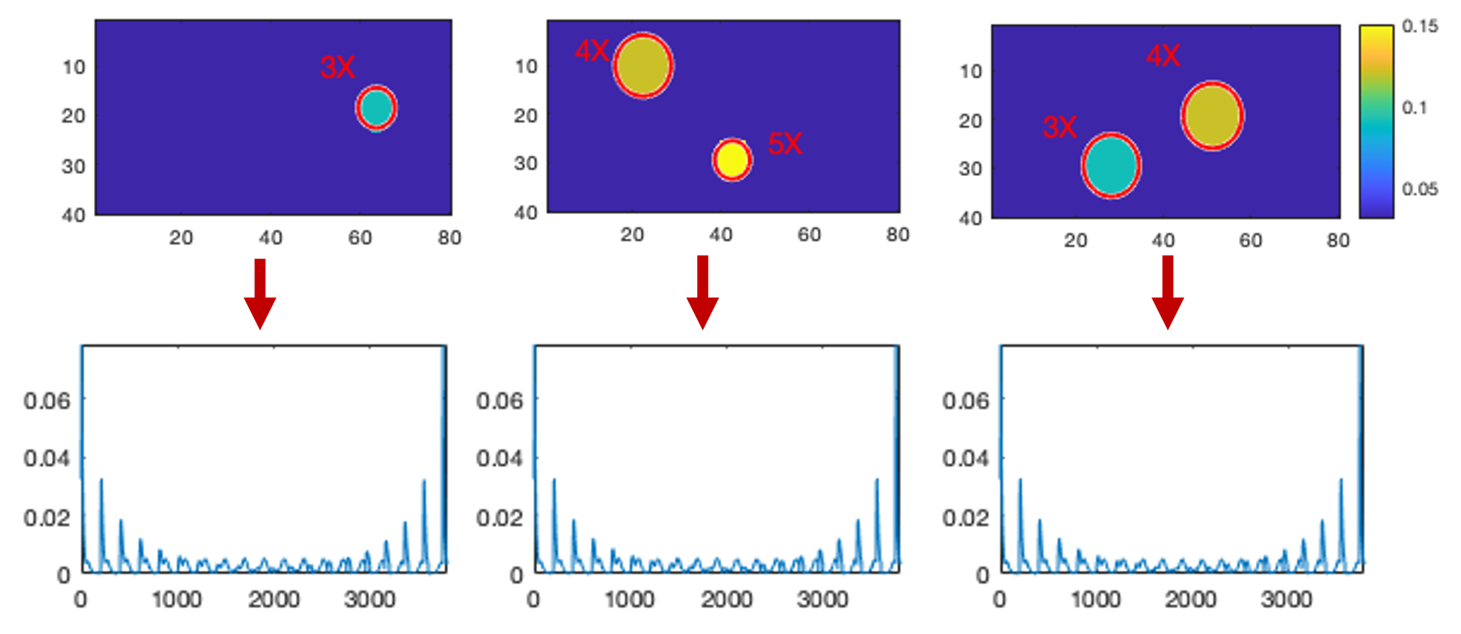}
\caption{Examples of synthetic fields of absorption
coefficient (top row) along with the corresponding computed DOT sinogram (bottom row) in the rectangular domain. Contrast regions display an increased absorption coefficient, indicated by the red label, with respect to the background value. The values
are presented as $\mu_a/D$, 
the red circles denote the ground truth
contrast regions.}
\label{fig:samples}
\end{figure}

We consider synthetic data
on 2D domains (a semi--disk or a rectangle)  in which circular contrast
regions of random size, absorption intensity and position are embedded
(see Fig.~\ref{fig:domains}). 
Light sources are positioned at constant intervals on the bottom side, 
1 mm inside the domain, so to that they are considered as volume sources.
This is clearly an approximation, but it is sufficiently accurate
for the present scopes. On this same part of the domain, 
null boundary conditions are enforced for the light field to represent
the physical plate in which lights
are embedded. The detectors are positioned on the 
left, top and right parts of the 
domain for the rectangular shape and 
on the curved portion of the domain for the
semi--disk. 
To train the network, we generated a set of 1500 samples, each including one or two circular contrast regions with 
random radius and position inside the domain. The 
contrast regions were
chosen to have absorption coefficient equal to $3,4$ or $5$-fold the background coefficient. When there were two contrast regions,  each region was allowed to have different
absorption coefficient.
The test dataset used to evaluate the performance of the proposed architecture and to compare it with 
other approaches contained 150 samples built according to the same rules as above.  
Observations are, where specified,    
polluted with Gaussian noise ($1\%$ to $5\%$
with respect to the local value).
Measures of light fluence 
at the detectors locations are generated by solving the DA model~(\ref{eq:eq_diff}) via a quadratic Finite Element code implemented via an in-house software 
in MATLAB 
using a small discretization parameter
for high resolution. We consider 9000 train 
samples.
The parameters used in the 
numerical tests are listed in Tab.~\ref{tab:param}. 
\begin{table}[h]
\caption{Parameter values used in the 
generation of the synthetic data}
\label{tab:param}
\begin{tabular}{@{}llll@{}}
\toprule
  Name & Symbol & Value & Units \\ \midrule 
  domain size (rectangle) &  & $10 \times 5$ & cm  \\
  domain size (semi--disk) &  & radius $5$ & cm  \\
 max radius of contrast  region & & 1.0 & cm \\
min radius of contrast region & & 0.5 & cm \\
background absorption coeff & $\mu_{a,0}$ & 
0.01 & cm$^{-1}$ \\
contrast region absorption coeff & & (3 to 5)$\times \mu_{a,0}$ & \\
scattering coeff & $\mu_{s}$ & 
0.1  {\text to}  1 & cm$^{-1}$ \\ 
no. of sources & $n_s$ & 19 & \\
no. of detectors & $n_d$ & 200 & \\
 \botrule
\end{tabular}
\end{table}
Fig.~\ref{fig:samples} shows some examples
of the absorption coefficient fields along 
with the corresponding DOT sinograms computed by the
FEM code for the rectangular domain. 
Notice again how the light fluence measures 
obtained for different absorption fields are very similar:
they have indeed a maximal point-to-point difference
of $10^{-4}$, which makes the reconstruction
problem very hard. 

We consider 150 test samples for evaluation.  
To assess the performance of the 
algorithms by first labeling in the reconstructed
signal image the voxels belonging to the contrast region.
To do this, the image is binarized and the connected
components are labeled out. 
Then, we evaluate the following indices:
\begin{itemize}
\item[$I_1:$] True Positive Ratio ($TPR$).
For each test sample, it is given by
$$
TPR=\displaystyle \frac{
\text{True Positives}}{\text{True Positives}
+\text{False Negatives}},
$$
where True Positives are the number of voxels
in which the algorithm correctly recognized
the contrast region, while False Negatives
are the number of voxels
in which the algorithm did not recognize the
presence of a contrast region existing in
the ground truth
\item[$I_2:$] Absolute Bias Error ($ABE$).
For each test sample, it is given by
$$
ABE=
\displaystyle \frac{1}{V} \sum_{i=1}^{V}
|\vmu_{a,GT}^{h,i}-\vmu_a^{h,i}|
$$
where $\vmu_{a,GT}^{h,i}$ is the ground truth
coefficient evaluated at the voxel 
\item[$I_3:$] Structural Similarity Index ($SSIM)$.
This is a perceptual metric that 
is usually adopted in video industry to quantify image quality degradation caused by processing such as data compression or by losses in data transmission.
Luminance, contrast and structure are the three parameters
that are included in its definition 
(see~\cite{renieblas2017structural} for its definition and computation in the
radiology field).
\end{itemize}

\subsection{Implementation of the neural networks}

In this section we outline the structure of the networks implemented for the Mod-DOT and E2D-DOT structures.

\begin{itemize}
    \item {\em data-AE}: it is a compact architecture which receives in input the sinogram $y_\delta^M$ of dimension $3800 (=200 \times 19)$.
Its encoder part $\mathcal{E}_y$ consists of a single fully connected (FC) layer with a hyperbolic tangent activation function. The generated latent space has size 800.
The decoder part $\mathcal{D}_y$ also features a single fully connected (FC) layer with a Sigmoid activation function. It takes an 800-dimensional input vector and generates a 3800-dimensional output vector.
In total, the data-AE encompasses $6,084,600$ trainable parameters.
\item {\em signal-AE}: it is a richer architecture which takes in input images of dimension $40 \times 80$
(or equivalently vectors of size 3200) describing the absorption coefficient field. For its construction, 
we conducted numerous experiments with different architectures and dimensions. Eventually, we considered two versions, one with fully connected layers and one with convolutional layers.
which are described in detail in Tab.\ref{tab:AE2conv}.
This encoder generates a bottleneck with dimension 800 ($4\times 10\times 20$), equal to the dimension of the latent space of the data-AE.
Again, this latter choice is the result of various experiments. 
\begin{table}[h]
\begin{tabular}{@{}lllllll@{}}
\toprule
 Name & Layer type &  Output shape &   kernel size & stride & padding & dilation\\
\midrule
$\mathcal{E}_\mu^1$ & Conv2d  & $[-1, 16, 40, 80]$ & $(3, 3)$ & $(1,1)$ & $(1,1) $ & -\\
\midrule
$\mathcal{E}_\mu^2$ & ReLU & $ [-1, 16, 40, 80]$ & - & - & - & - \\
\midrule
$\mathcal{E}_\mu^3$ & MaxPool2d & $[-1, 16, 20, 40]$  & $2$ & $2$ & - & $1$ \\
\midrule
$\mathcal{E}_\mu^4$ & Conv2d  & $[-1, 8, 20, 40]$& $(3, 3)$ & $(1,1) $ & $(1,1) $ & - \\
\midrule
$\mathcal{E}_\mu^5$ & ReLU & $[-1, 8, 20, 40]$ & - & - & - & -  \\
\midrule
$\mathcal{E}_\mu^6$ & MaxPool2d  & $[-1, 8, 20, 40]$  & $2$ & $2$  & - & $1$\\
\midrule
$\mathcal{E}_\mu^7$& Conv2d  & $[-1, 4, 10, 20]$  & $(3, 3)$ & $(1,1) $  & $(1,1) $ & -\\ 
\midrule
$\mathcal{E}_\mu^8$ & ReLU & $[-1, 4, 10, 20]$& - & - & - & -  \\
\midrule\midrule
$\mathcal{D}_\mu^1$ & ConvTranspose2d  & $[-1, 8, 20, 40]$  & $(2, 2)$ & $(2, 2)$ & - &  -  \\
\midrule
$\mathcal{D}_\mu^2$  & ReLU & $[-1, 8, 20, 40]$ & - & - & - & -   \\
\midrule
$\mathcal{D}_\mu^3$  & ConvTranspose2d  & $[-1, 16, 40, 80]$ & $(2, 2)$ &  $(2, 2)$ & - & - \\
\midrule
$\mathcal{D}_\mu^4$  & ReLU & $[-1, 16, 40, 80]$ & - & - & - & -   \\
\midrule
$\mathcal{D}_\mu^5$  & ConvTranspose2d  & $[-1, 1, 40, 80]$  & $(3, 3)$ & $(1,1)$ & $(1, 1)$ & - \\
\midrule
$\mathcal{D}_\mu^6$  & Sigmoid & $[-1, 1, 40, 80]$& - & - & - & -   \\
\botrule 
\end{tabular}
    \caption{Structure of the signal-AE with convolutional layers. In total, there are $2421$ trainable parameters.
    The first 8 layers compose the encoder, the remaining the decoder.}
    \label{tab:AE2conv}
\end{table}

\item {\em bridge-NN}: it is a network composed of seven fully connected layers with hyperbolic tangent activation functions. Each layer has the same input and output dimensions of 800. This network has 4,485,600 trainable parameters
\item {\em denoiser NN}: both the Mod-DOT and E2D-DOT nets were complemented with a downstream network based on convolutional layers which helps in removing the remaining artifacts upon inversion.  
The network architecture, called DEN, is described in Tab.~\ref{tab:Den}. In total, it comprises 185,217 trainable parameters. 

\begin{table}[h]
\begin{tabular}{@{}llllll@{}}
\toprule
   Name & Layer type & Output shape &  kernel size & stride & padding \\
    \midrule
    $DEN_1$ & Conv2d &$ [-1, 32, 40, 80] $ & $(3, 3)$ & $(1,1)$ & $same$\\
    \midrule
    $DEN_2$ & LeakyReLU & $[-1, 32, 40, 80]$ & - & -  & -\\
    \midrule
    $DEN_3$ & Conv2d  & $[-1, 64, 40, 80]$ & $(3, 3)$ & $(1,1)$ & $same$\\
    \midrule
    $DEN_4$ & LeakyReLU & $[-1, 64, 40, 80]$  & - & -  & -  \\
    \midrule
    $DEN_5$ & Conv2d  & $[-1, 128, 40, 80]$  & $(3, 3)$ & $(1,1)$ & $same$\\
    \midrule
    $DEN_6$ & LeakyReLU & $[-1, 128, 40, 80]$  & - & - & -  \\
    \midrule
    $DEN_7$ &  ConvTranspose2d  & $[-1, 64, 40, 80]$ & $(3, 3)$ & $(1,1)$ & $(1,1)$\\
    \midrule
    $DEN_8$ & LeakyReLU & $[-1, 64, 40, 80]$  & - & - & - \\
    \midrule
    $DEN_9$ &  ConvTranspose2d  & $[-1, 32, 40, 80]$  & $(3, 3)$ & $(1,1)$ & $(1,1)$\\
    \midrule
    $DEN_{10}$ & LeakyReLU & $[-1, 32, 40, 80]$  & - & - & -\\
    \midrule
    $DEN_{11}$ &  ConvTranspose2d  & $[-1, 1, 40, 80]$ & 
    $(3, 3)$ & $(1,1)$ & $(1,1)$\\ 
    \midrule
    $DEN_{12}$ & LeakyReLU & $[-1, 1, 40, 80] $  & - & - & -\\
    \midrule
\end{tabular}
    \caption{Structure of the denoising network based on convolutional layers.}
    \label{tab:Den}
\end{table}
\end{itemize} 

\null 

\subsubsection{Training strategy} \label{subsec:TrainingAlg}
According to the modular philosophy
adopted for the Mod-DOT architecture, a decoupled or partially decoupled training
strategy is suggested. This approach allows
to take advantage of intrinsic regularities
in the data and signal manifolds, as discussed above. 
First, the two AEs are trained separately, in a sort of pre--train
phase,
then the trained weights 
of the encoder 
$\mathcal{E}_y$ and of the decoder
$\mathcal{D}_\mu$ are 
used as initial guess to train the 
complete bridged network
$\mathcal{N}=
\mathcal{D}_\mu \circ 
 \Sigma \circ \mathcal{E}_y
$.
The utility of the pre--train phase has been assessed via an extensive experimentation, and it is especially evident for the 
signal-AE component. On the contrary, the bridge is trained in a coupled fashion with the pre-trained portions of the two AEs. 
Observe that an analogous decoupled strategy is used also in the DL-ROM approach for the AE part. On the contrary with respect
to our approach, in the DL-ROM method also the reduced map is trained separately. We have however observed that in the case of the solution of the inverse problem, a fully decoupled approach like this latter does not yield satisfactory results, since the bridge-NN (which conceptually corresponds to the reduced map) has to approximate a very complex map to 
bridge the latent spaces. 
For pre-training, we used the Adam optimizer, batch size 64, and learning rate $5 \times 10^{-5}$.
The number of epochs was set to $500$ for the data-AE and $10000$ for the signal-AE. For both the AEs, we used the Means Square Error cost function. 

\null

\noindent In the training phase of both the Mod-DOT and E2D-DOT networks, $10000$ epochs were used to allow for prolonged learning periods, considering the complexity of the 
target to achieve. A batch size of $64$ was utilized for each iteration, and the dataset was shuffled at the start of every epoch.
The Adam optimization algorithm was employed with learning rate $5 \times 10^{-5}$

Unlike the pretraining networks, various choices for the loss function were investigated, and, namely, we considered the following options ($\#N$ is the number of samples):
\begin{itemize}
\item[(a)] mean square error:
$
\displaystyle loss_{MSE}=\frac{1}{\#N}\sum_{i=1}
^{\#N}[\mathcal{N}(\vy_{\delta,i}^M)-
\vmu_{a,i}^h]^2$
\item[(b)] mean square error + 
$\ell_1$ norm for feature selection
(with $\alpha$ regularization parameter):\\
$
 \displaystyle loss_{MSE+\ell_1}=\frac{1}{\#N}\sum_{i=1}
^{\#N}[\mathcal{N}(\vy_{\delta,i}^M)-
\vmu_{a,i}^h]^2 + \alpha
\sum_{i=1}
^{\#N}[\mathcal{N}(\vy_{\delta,i}^M)-
\vmu_{a,i}^h]
$
\item[(c)] mean square error + 
control on AE error:\\
$
 \displaystyle loss_{MSE+AE}=\frac{1}{\#N}\sum_{i=1}
^{\#N}[\mathcal{N}(\vy_{\delta,i}^M)-
\vmu_{a,i}^h]^2 + 
[\Psi_y(\vy_{\delta,i}^M)-\vy_{\delta,i}^M]^2+
[\Psi_\mu(\vmu_{a,i}^h)-\vmu_{a,i}^h]^2
$
\end{itemize}

\noindent All the considered neural networks have been implemented 
in the 
{\tt pytorch} framework. Computations 
have been run on the HPC cluster INDACO owned by
University of Milano, equipped with 
Nvidia A100 GPUs (2.3GHz, 32 cores, 1TB RAM).

\subsection{Implementation of the classic variational approaches}
\label{ssec:varappr}
In order to compare the performance of the proposed NN approach 
with classic methods, we consider two strategies based 
on the Rytov linearization using 
different regularization procedures: 
\begin{itemize}
    \item Elastic Net regularization~\cite{Causin19}: we
solve the problem 
\begin{equation}
(\mua^h)^* \in \argmin_{\mua^h} \frac12\|\vJ \mua-{\bm b}\| + \alpha\lp\vartheta\|\mua^h\|_1+(1-\vartheta)\|\mua^h\|_2^2\rp,   \end{equation}
where the regularizator is a convex combination of the $\ell_2$ and $\ell_1$ norms, with $\vartheta\in[0,1]$. This approach aims
at combining the sparse--promoting and peak-enhancing property
of the $\ell_1$ norm and the robustness of the Tikhonov regularization. The solution is computed 
using the MATLAB implementation of {\tt glmnet}~\cite{glmnet}, 
a highly efficient package that fits generalized linear and similar models via penalized maximum likelihood. 
Cross validation is used to obtain the optimal regularization parameter
\item  Bregman iterations~\cite{Benfenati14,Benfenati2020}:  this is a strategy where at each step where the regularization operator is substituted by its the Bregman Divergence computed at the previous iterate. If $\mathcal{R}$ satisfies suitable hypotheses (namely, being a convex, proper, lower semicontinuous function) and denoting by $\partial \mathcal{R}(\boldsymbol{\xi})$ the subdifferential of $\mathcal{R}$ computed at ${\boldsymbol{\xi}}$ (see \cite{rockafellar2015convex} for technical details), then the Bregman divergence of $\mathcal{R}$ is given by
\[
D_\mathcal{R}^{\vp}({\boldsymbol{\xi}}, {\boldsymbol{\eta}}) = \mathcal{R}({\boldsymbol{\xi}}) - \mathcal{R}({\boldsymbol{\eta}}) - \langle \vp,
{\boldsymbol{\xi}}-{\boldsymbol{\eta}}\rangle, \quad \vp \in\partial \mathcal{R}({\boldsymbol{\xi}}),
\]
where $\vp^k\in\partial\mathcal{R}(\mua^k)$.
An iterative procedure 
is used, such that the functional 
to be minimized at the $(k+1)$--th step is 
given by 
\begin{equation}
\left\{
\begin{array}{l}
\displaystyle \mua^{k+1}  \in   \argmin_{\mua^h} \displaystyle\frac12\|\vJ\,\mua^h
-{\bm b}\|_2^2 + \alpha D_\mathcal{R}^{\vp^k}(\mua^h,\mua^k)- \alpha \langle \vp^k, 
\mua^h\rangle \\[3mm]
 \displaystyle\vp^{k+1}  =   \vp^k - \frac 1\alpha \vJ^\top\lp\vJ\,\mua^h -{\bm b}\rp.
\end{array}
\right. 
\end{equation}
The above Bregman procedure is coupled with the choice of $\ell_1$ regularization. 
\end{itemize}

\subsection{Numerical results: qualitative and quantitative evaluation}

\subsubsection{Comparison of Mod-DOT vs classic approaches}
\begin{figure*}[!t]
\newcommand{\factor}{0.18}
\begin{center}
	\subfloat[\label{fig:gt}]{\includegraphics[width=\factor\textwidth]{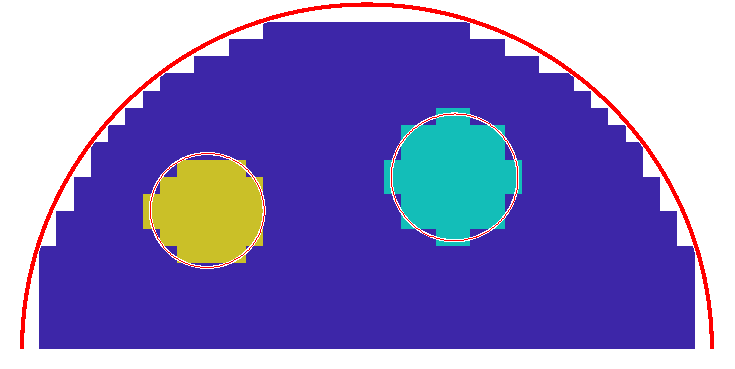}\hspace{0.025\textwidth}\includegraphics[width=\factor\textwidth]{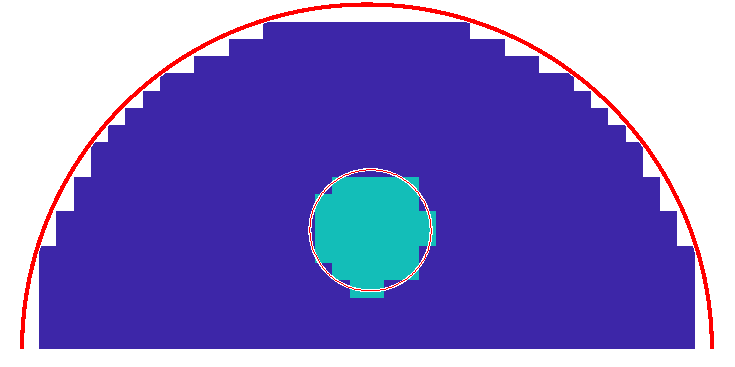}\hspace{0.025\textwidth}\includegraphics[width=\factor\textwidth]{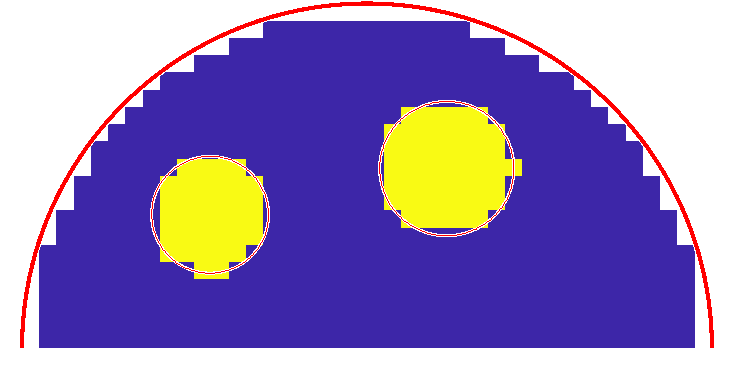}\hspace{0.025\textwidth}\includegraphics[width=\factor\textwidth]{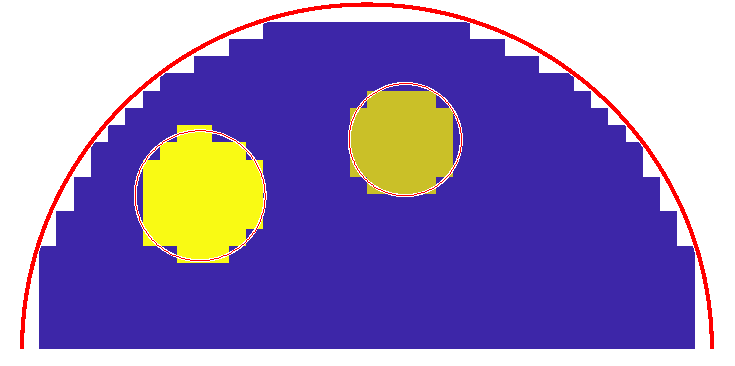}\hspace{0.025\textwidth}\includegraphics[width=\factor\textwidth]{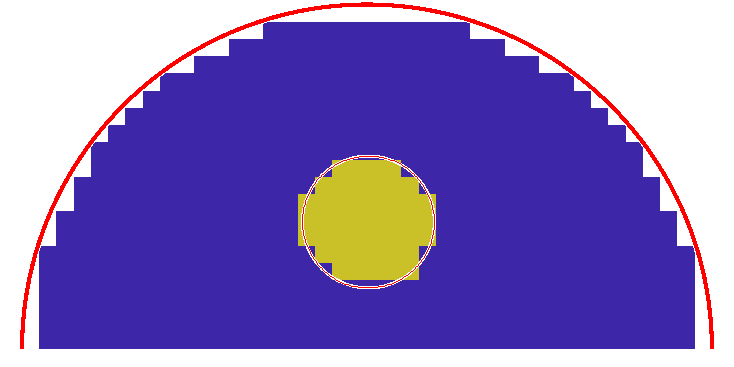}}
	
\subfloat[\label{fig:nfree}]{\includegraphics[width=\factor\textwidth]{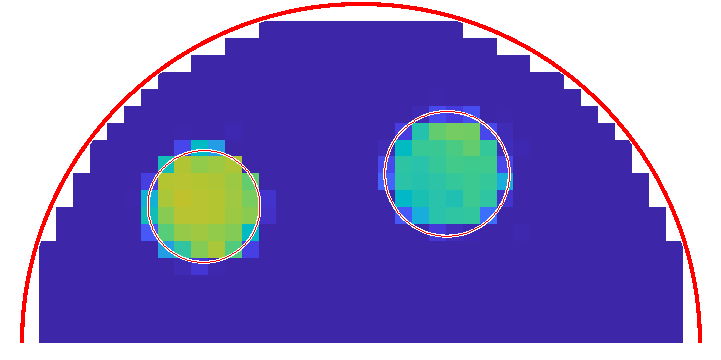}\hspace{0.025\textwidth}\includegraphics[width=\factor\textwidth]{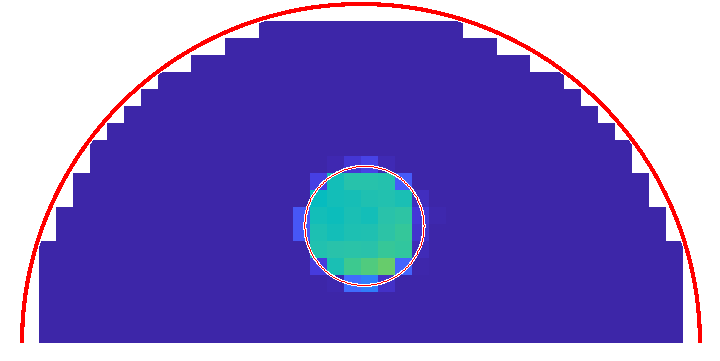}\hspace{0.025\textwidth}\includegraphics[width=\factor\textwidth]{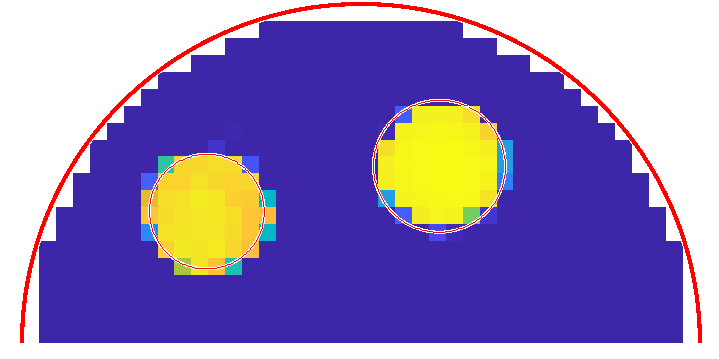}\hspace{0.025\textwidth}\includegraphics[width=\factor\textwidth]{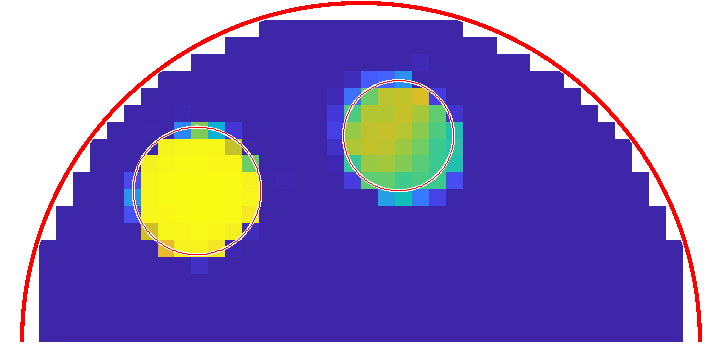}\hspace{0.025\textwidth}\includegraphics[width=\factor\textwidth]{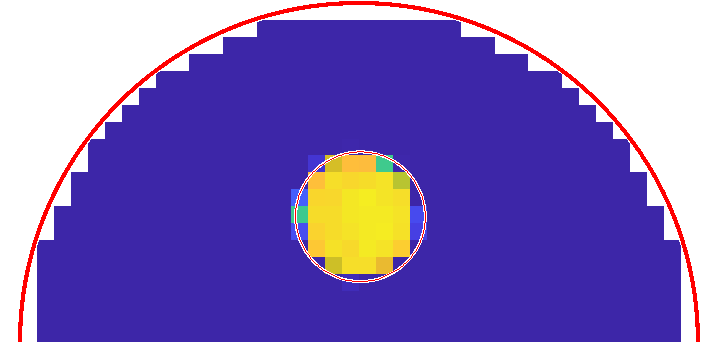}}

\subfloat[\label{fig:G1}]{\includegraphics[width=\factor\textwidth]{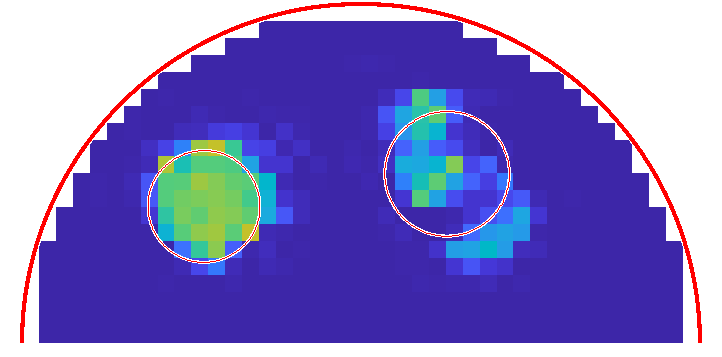}\hspace{0.025\textwidth}\includegraphics[width=\factor\textwidth]{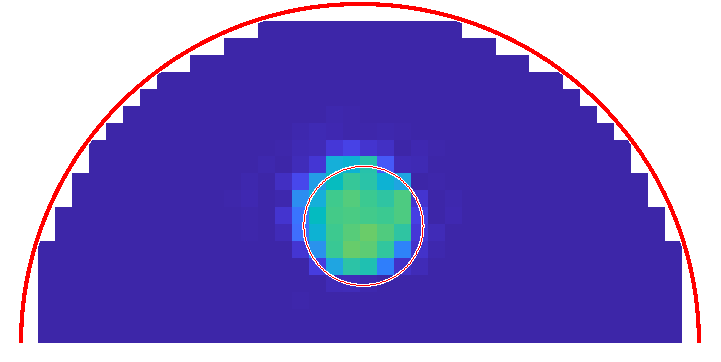}\hspace{0.025\textwidth}\includegraphics[width=\factor\textwidth]{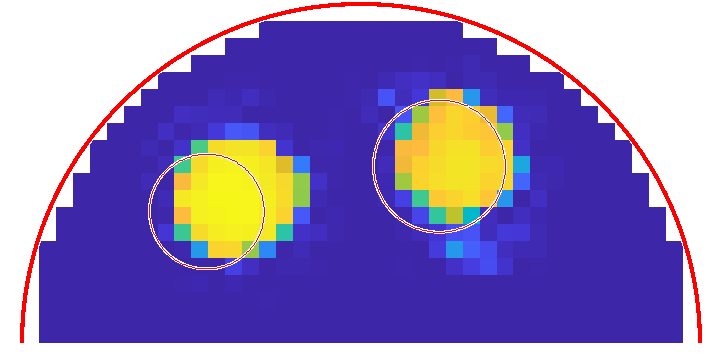}\hspace{0.025\textwidth}\includegraphics[width=\factor\textwidth]{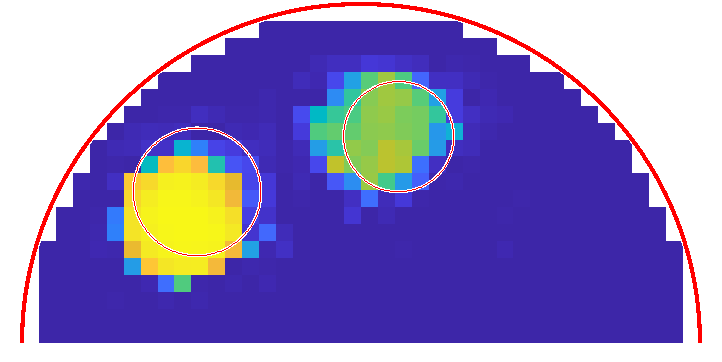}\hspace{0.025\textwidth}\includegraphics[width=\factor\textwidth]{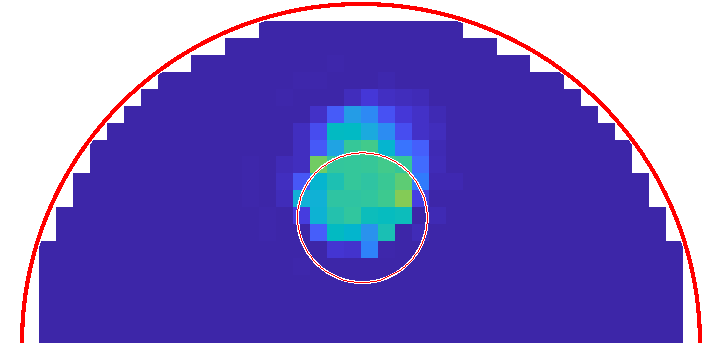}}
	
	\subfloat[\label{fig:G3}]{\includegraphics[width=\factor\textwidth]{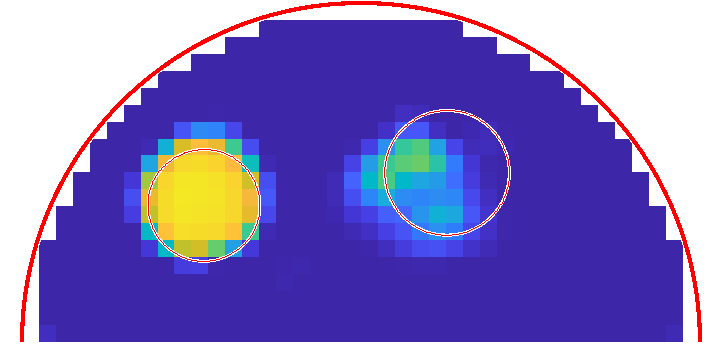}\hspace{0.025\textwidth}\includegraphics[width=\factor\textwidth]{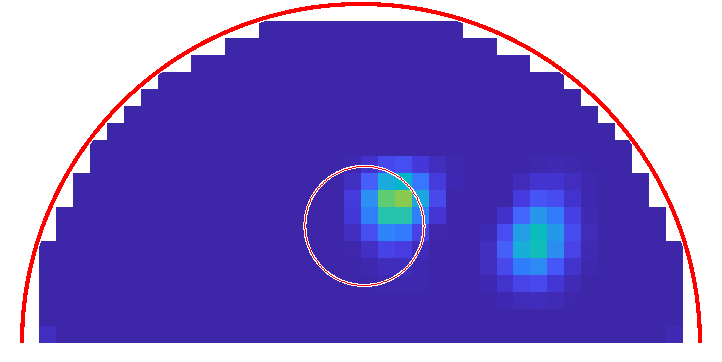}\hspace{0.025\textwidth}\includegraphics[width=\factor\textwidth]{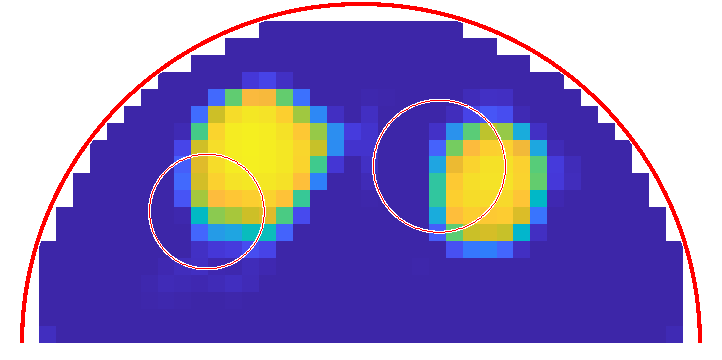}\hspace{0.025\textwidth}\includegraphics[width=\factor\textwidth]{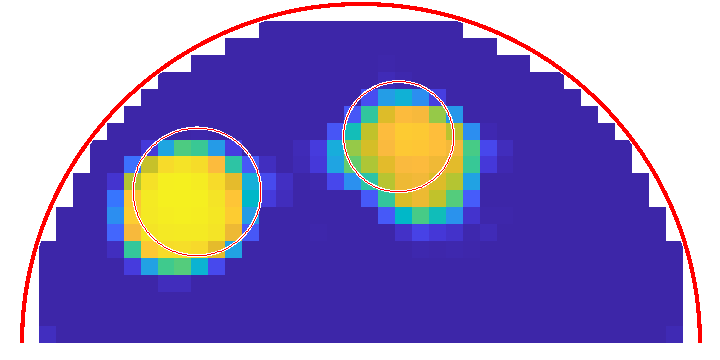}\hspace{0.025\textwidth}\includegraphics[width=\factor\textwidth]{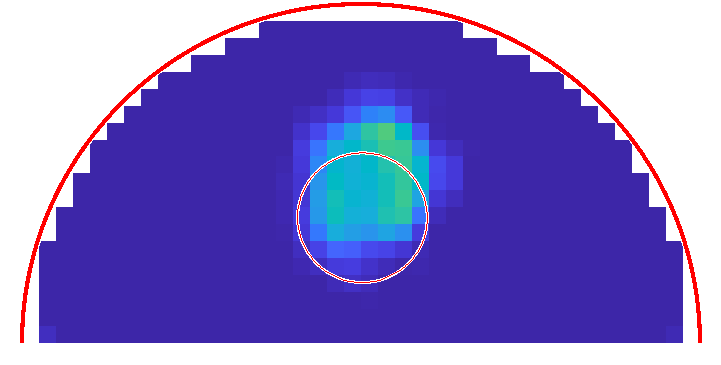}}	
	
\subfloat[\label{fig:G5}]{\includegraphics[width=\factor\textwidth]{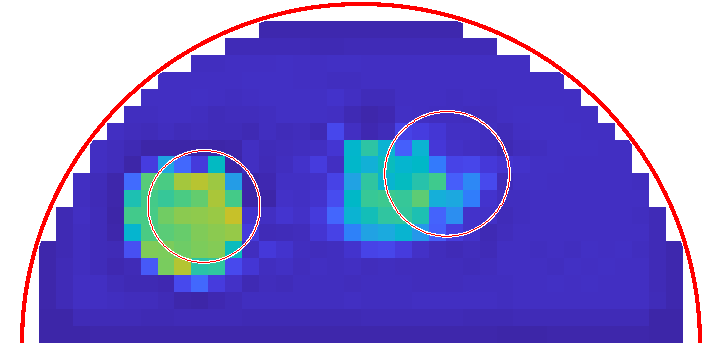}\hspace{0.025\textwidth}\includegraphics[width=\factor\textwidth]{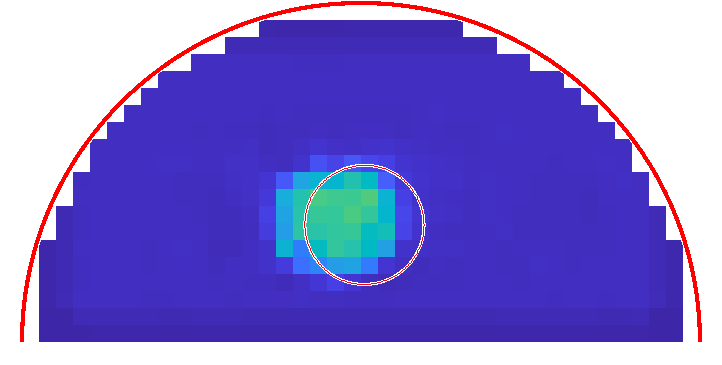}\hspace{0.025\textwidth}\includegraphics[width=\factor\textwidth]{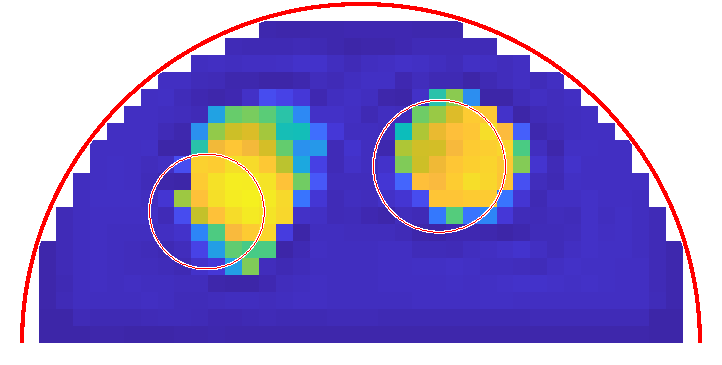}\hspace{0.025\textwidth}\includegraphics[width=\factor\textwidth]{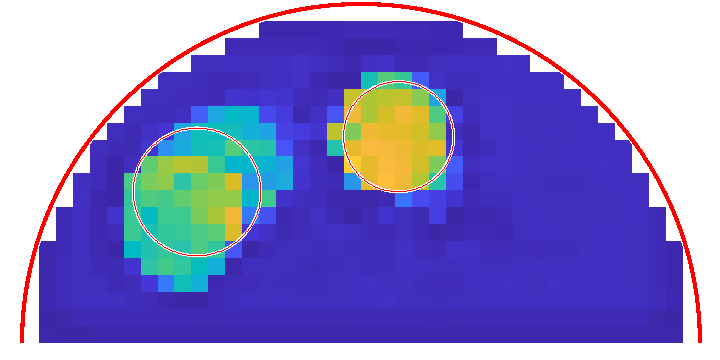}\hspace{0.025\textwidth}\includegraphics[width=\factor\textwidth]{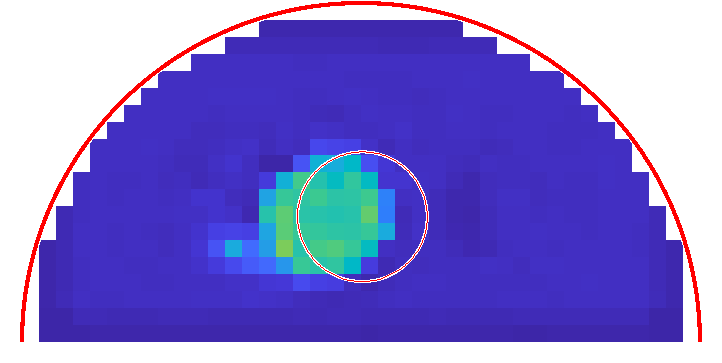}}

\end{center}
\subfloat{\qquad \quad \includegraphics[width=.85\textwidth]{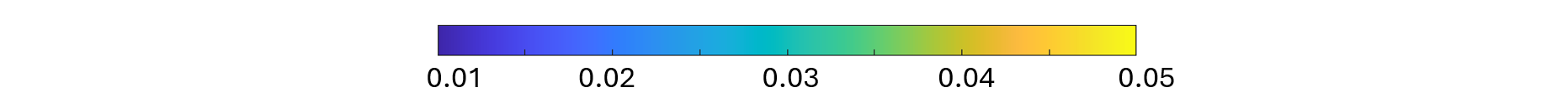}}
\caption{Reconstruction of the absorption coefficient distribution obtained with the Mod-DOT approach. \eqref{fig:gt} Ground truth. The 
red circles denote the true contrast regions along with the corresponding absorption coefficient, expressed 
as a multiplier of the background value. (a): noiseless data; (b): 1\% noise; (c): 3\% noise; (d): 5\% noise. The red circles denote the ground truth contrast regions.}
\label{fig:fResults}
\end{figure*}

\begin{figure}[!t]
	\newcommand{\factor}{0.2}
	\begin{center}
		\subfloat[\label{fig:ela0}]{\includegraphics[width=\factor\textwidth]{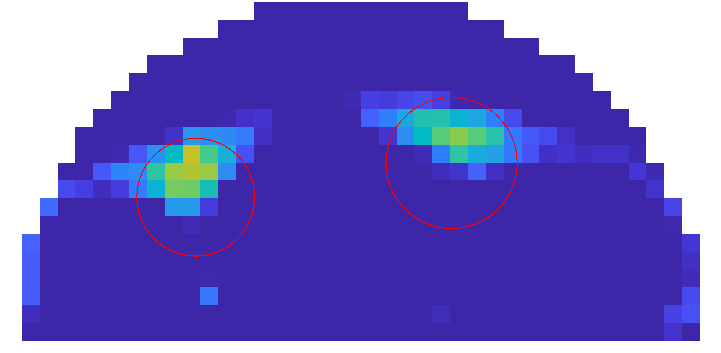}\hspace{0.025\textwidth}\includegraphics[width=\factor\textwidth]{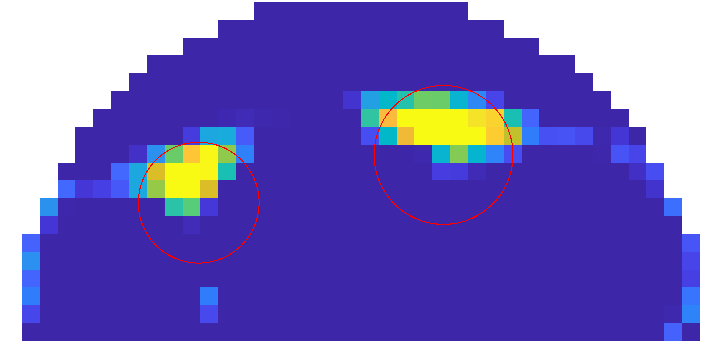}\hspace{0.025\textwidth}\includegraphics[width=\factor\textwidth]{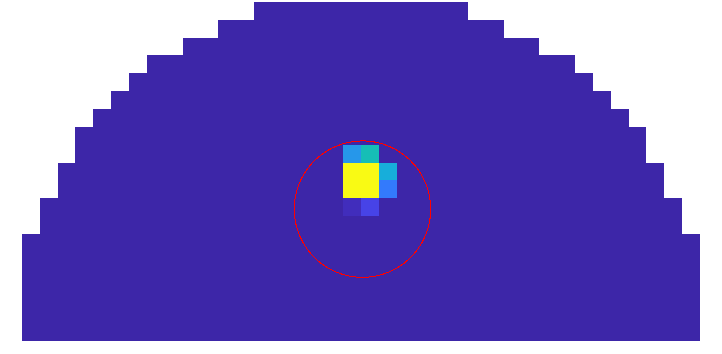}}
		
	\subfloat[\label{fig:ela1}]{\includegraphics[width=\factor\textwidth]{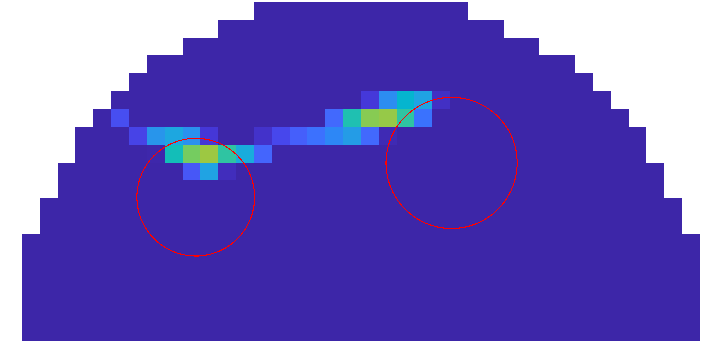}\hspace{0.025\textwidth}\includegraphics[width=\factor\textwidth]{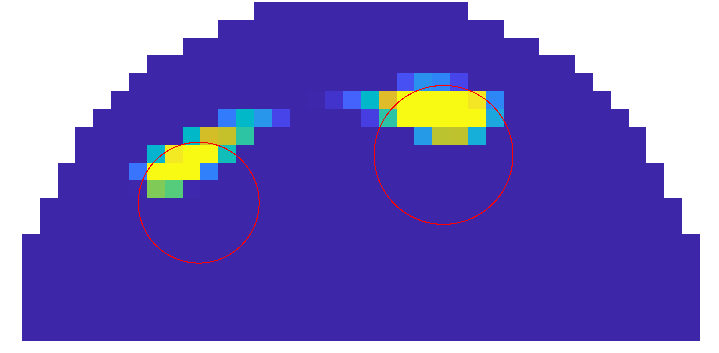}\hspace{0.025\textwidth}\includegraphics[width=\factor\textwidth]{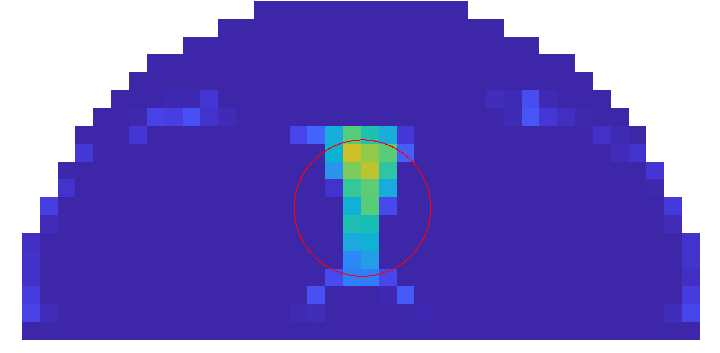}}

	\subfloat[\label{fig:ela3}]{\includegraphics[width=\factor\textwidth]{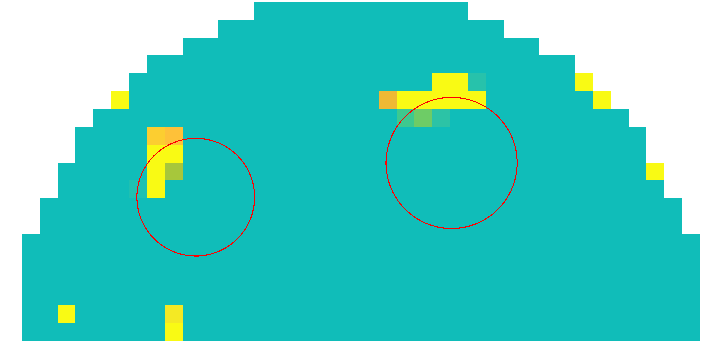}\hspace{0.025\textwidth}\includegraphics[width=\factor\textwidth]{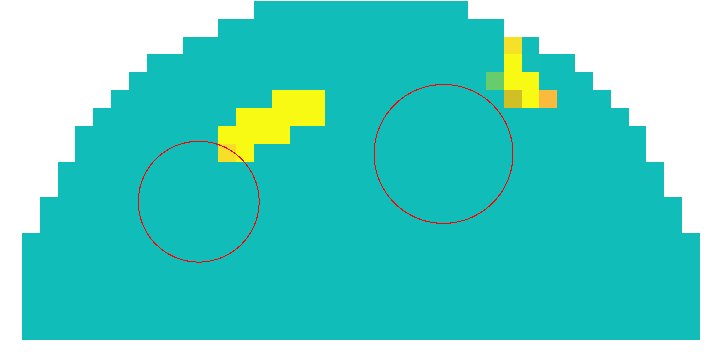}\hspace{0.025\textwidth}\includegraphics[width=\factor\textwidth]{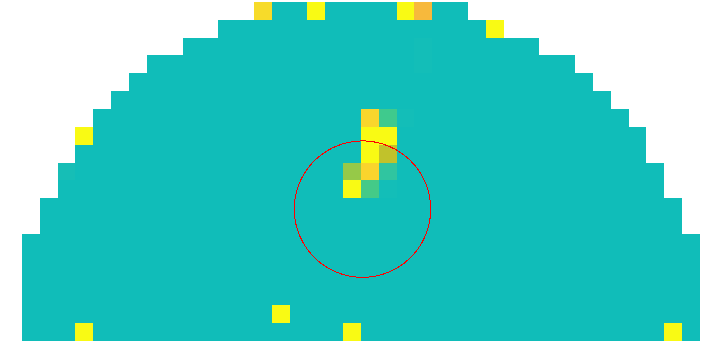}}

 \subfloat{\includegraphics[width=.85\textwidth]{Fig8C}}

	\end{center}
	\caption{Reconstruction via the Elastic Net approach of the test samples of columns 1, 3 and 5 in \cref{fig:gt}.
  (a): noiseless data; (b): 1\% noise; (c): 3\% noise. The red circles denote the ground truth contrast regions.}
	\label{fig:ela}
\end{figure}

\begin{figure}[!t]
	\newcommand{\factor}{0.2}
	\begin{center}
		\subfloat[\label{fig:breg0}]{\includegraphics[width=\factor\textwidth]{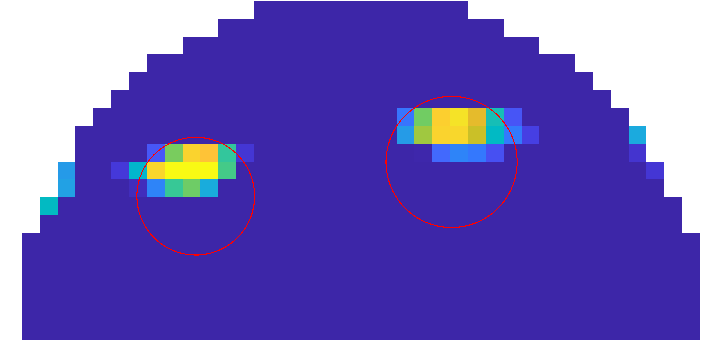}\hspace{0.025\textwidth}\includegraphics[width=\factor\textwidth]{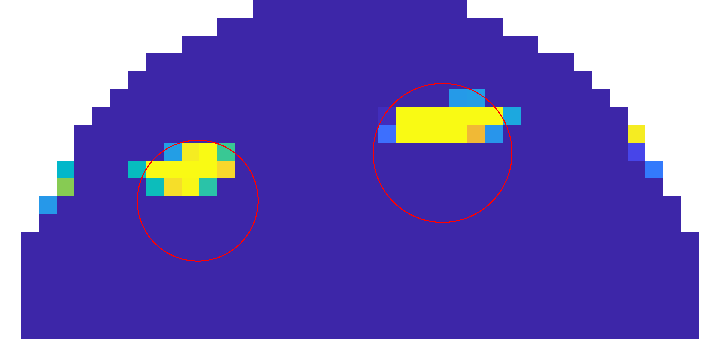}\hspace{0.025\textwidth}\includegraphics[width=\factor\textwidth]{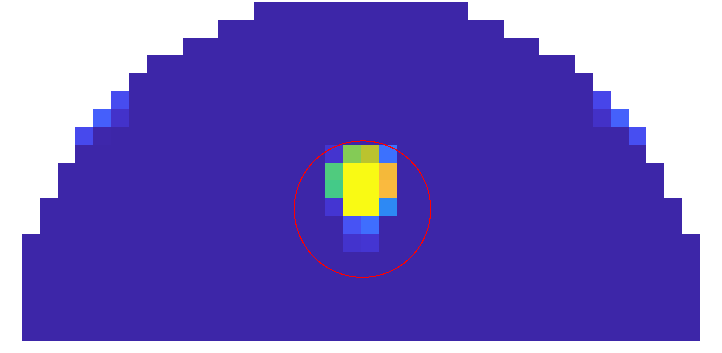}}
		
\subfloat[\label{fig:breg1}]{\includegraphics[width=\factor\textwidth]{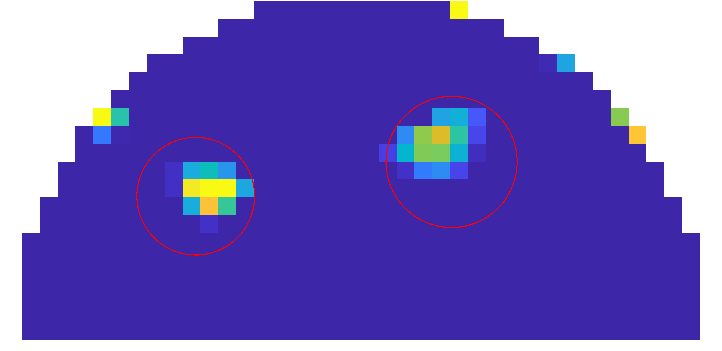}\hspace{0.025\textwidth}\includegraphics[width=\factor\textwidth]{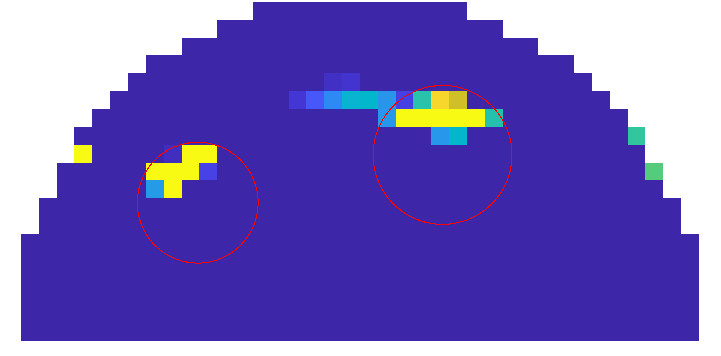}\hspace{0.025\textwidth}\includegraphics[width=\factor\textwidth]{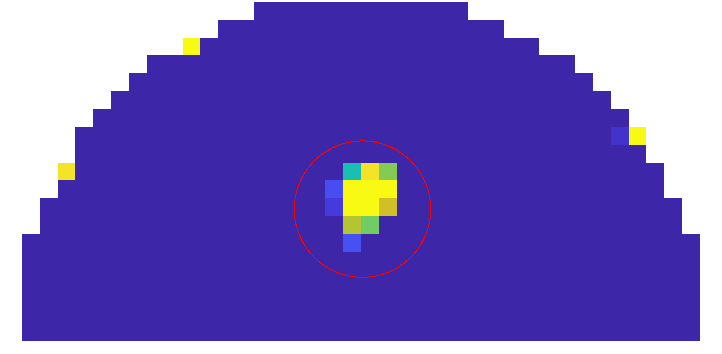}}
		
		\subfloat[\label{fig:breg3}]{\includegraphics[width=\factor\textwidth]{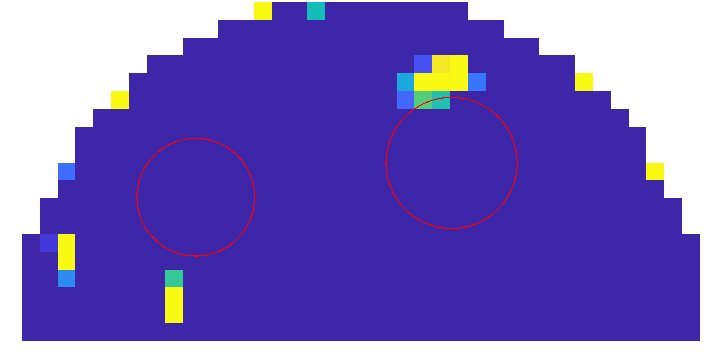}\hspace{0.025\textwidth}\includegraphics[width=\factor\textwidth]{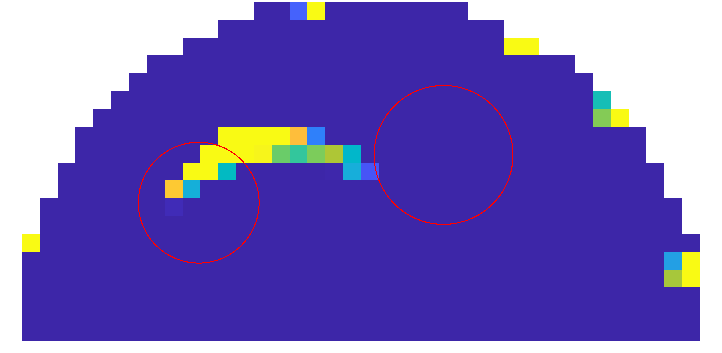}\hspace{0.025\textwidth}\includegraphics[width=\factor\textwidth]{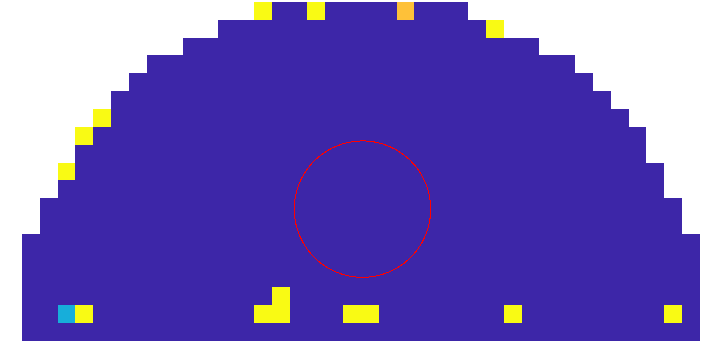}}

       \subfloat{\includegraphics[width=.85\textwidth]{Fig8C}}
  
	\end{center}
	\caption{Reconstruction via the Bregman approach. (a): noiseless data; 
	(b): 1\% noise; (c): 3\% noise. The red circles denote the ground truth contrast regions.}
	\label{fig:breg}
\end{figure}

We consider for this study the semi-disk geometry and we compare the performance of the
Mod-DOT approach with fully connected layers in the signal-AE  vs the reconstruction obtained by the Elastic Net and Bregman iterations approaches.    
Fig.~\ref{fig:fResults} presents the absorption coefficient reconstructed by the Mod-DOT for five different test (unseen)  samples in the test set with varying levels of noise. In each figure, the red circles depict the ground truth contrast regions. With noiseless data (second row of Fig.~\ref{fig:fResults}), the contrast region is almost perfectly recovered, both in location, shape and intensity. The successive rows refer to the reconstruction obtained when Gaussian noise pollutes the data. As it is to be expected, increasing noise levels worsen the overall quality of the reconstruction, even if the gross location and shape of the perturbed region are still correctly captured. 
Fig.~\ref{fig:ela} (Elastic-Net) and Fig.~\ref{fig:breg} (Bregman)  
show the results corresponding to the test samples of columns 1, 3 and 5 in \cref{fig:gt}. Observe as already in the noiseless case (first row) both these classic methods do not achieve a good reconstruction. As the noise level increases the reconstruction quality 
further deteriorates, so that noise level 5\% is not reported here for the classic methods due to complete lack of significance. Simulations performed reducing the voxel size, from the original 0.25~cm, in the
reconstruction show partially improved,
although this does not solve the difficulty in handling noisy data.

We summarize in Tab.~\ref{tab:results} the ACR and TPR metrics for this test case. The table presents the averages 
on the 150 samples of the test set for the Mod-DOT approach and, for comparison,
for the two considered classic approaches. Notice that the ACR results have been binned according to
the intensity of the contrast region. 
\begin{table*}[htbp]
	\begin{center}
	\begin{tabular}{lllll}
	\toprule
		noise & ${3\mu_{a,0}}$ (GT: 3e-2)&  ${4\mu_{a,0}}$ (GT: 4e-2) &  ${5\mu_{a,0}}$ (GT: 5e-02) & TPR\\
		\toprule
		&\multicolumn{4}{c}{Mod-DOT}\\
		\midrule
		0\% & 3.05e-02 $\pm$ 2.00e-03 & 4.08e-02 $\pm$ 3.72e-03 & 4.74e-02 $\pm$ 2.19e-03 & 0.94\\

		1\% & 3.07e-02 $\pm$ 2.77e-03 & 3.76e-02 $\pm$ 5.54e-03 & 4.38e-02 $\pm$ 4.01e-03 & 0.76\\

		3\% & 3.06e-02 $\pm$ 5.56e-03 & 3.61e-02 $\pm$ 6.10e-03 & 4.11e-02 $\pm$ 4.65e-03 & 0.52\\

		5\% & 3.00e-02 $\pm$ 3.82e-03 & 3.15e-02 $\pm$ 4.44e-03 & 3.39e-02 $\pm$ 6.04e-03 & 0.46\\
		
		\midrule		
		&\multicolumn{4}{c}{Elastic Net}\\
		\midrule
		0\% & 2.73e-02 $\pm$ 4.74e-03 & 3.45e-02 $\pm$ 4.96e-03 & 3.90e-02 $\pm$ 8.69e-03 & 0.45\\
		1\% & 2.91e-02 $\pm$ 9.36e-03 & 4.30e-02 $\pm$ 8.70e-03 & 5.01e-02 $\pm$ 9.80e-03 & 0.17\\
		3\% & 9.55e-02 $\pm$ 1.92e-02 & 1.25e-01 $\pm$ 2.83e-02 & 1.23e-01 $\pm$ 3.36e-02 & 0.05\\
            \midrule
&\multicolumn{4}{c}{Bregman}\\
		\midrule
		0\% & 4.02e-02 $\pm$ 8.95e-03 & 5.93e-02 $\pm$ 1.61e-02 & 8.54e-02 $\pm$ 2.98e-02 & 0.26\\
		1\% & 4.77e-02 $\pm$ 1.81e-02 & 5.85e-02 $\pm$ 1.84e-02 & 8.34e-02 $\pm$ 2.85e-02 & 0.17\\
		3\% & 1.28e-01 $\pm$ 5.50e-02 & 1.45e-01 $\pm$ 9.06e-02 & 1.36e-01 $\pm$ 6.80e-02 & 0.03\\
		\midrule		
	\end{tabular}
\end{center}
\caption{ACR (binned according to
the intensity of the contrast region) and TPR metrics for different noise levels. Values
		are averaged over 150 test samples.}\label{tab:results}
\end{table*}
For all the approaches, increasing levels of noise significantly worsen the TPR, 
as well as the ACR. However, the Mod-DOT reconstruction  
is significantly less affected by quality loss. 
The contrast regions reconstructed with this approach still show a notable contrast with respect to the background value, and 46\% of the
extension of the perturbed regions is correctly recovered on average even for the highest noise level (vs  
5\% and 3\% for the variational approaches, respectively).
One should notice that in the Elastic Net approach the reconstructed intensities constantly
approximate the nominal one by defect for low noise level: this is not fully surprising, since 
regularization comes at the cost of a certain degree of smearing/blurring of the solution. The approach based on Bregman iterations increase the contrast, as already observed in previous works, although it struggles with noisy data.

\subsubsection{Comparison of networks with fully connected vs convolutional layers}
In this study we show that
convolutory layers perform better
than fully connected layers in the 
signal-AE, especially when
using noisy data for DOT reconstruction. 
With this aim, we consider the rectangular geometry  with
circular contrast regions of different size, location and 
intensity. In this case
the absorption coefficient samples are presented to the nets as $80 \times 40$ images. Training is carried out as specified above.
We study the performance metrics 
both for the E2E and Mod-DOT approaches. In each case, we apply
denoising and we consider the loss
function which gives best results. 
Tab.~\ref{tab:indici_medi_NOISEFREE}
refers to the noiseless case, while 
Tab.~\ref{tab:indici_medi_NOISE}
refers to the noisy data. 
It is apparent that in the noiseless case there is no evident advantage
to use convolutory layers in terms
of the metrics. The number of parameters is, however, always in favour of the convolutory filters. For example, 
the Mod-DOT-FC architecture has
10,274,817 parameters while the 
Mod-DOT-CONV architecture has
7,712,426 parameters.
The difference is more relevant 
when we consider noisy data. 
As a matter of fact, not only the 
the version of the net with convolutory
filters has less parameters, but also the metrics show a better quality of the reconstruction. Based on these findings, in the following studies we always adopt a signal-AE based on convolutory layers. 

\begin{table}[htbp]
    \centering
    
    \begin{tabular}{lllll}
        \toprule
        structure & ABE & MSE & TPR & SSIM \\
        \toprule
        E2E - FC  \small{(MSE+L1)} &$1.08974\cdot 10^{-3}$ & $ 3.87490\cdot 10^{-5}$ & $0.93115$ & $0.98792$\\
        \midrule
        Mod-DOT - FC \small{(MSE+L1)} & $ 1.55653\cdot 10^{-3}$ & $8.72685\cdot 10^{-5}$ & $0.88245$ & $0.97509$ \\
         \midrule
        E2E - CONV \small{(MSE)} & $ 1.11889\cdot 10^{-3}$ & $3.60208\cdot 10^{-5}$ & $0.92200$ & $0.98850$ \\
        \midrule
        Mod-DOT - CONV \small{(MSE+L1)} & $ 1.60053\cdot 10^{-3}$ & $9.59041\cdot 10^{-5}$ & $0.86057$ & $0.97317$ \\
        \bottomrule 
    \end{tabular}
    \caption{Performance metrics for the different architectures we consider. For each method it is
    specified the loss function which performs better. FC=fully connected, CONV=convolutional layers.}
    \label{tab:indici_medi_NOISEFREE}
\end{table}

\begin{table}[htbp]
    \centering
    
    \begin{tabular}{llllll}
    \toprule
        structure  & noise & ABE & MSE & TPR & SSIM \\
          \toprule
        \multirow{3}{*}{E2E - FC \small{(MSE+L1)}} & $1\%$ & $9.0472 \cdot 10^{-3}$   & $5.8909\cdot 10^{-4}$  & $ 0.7939$ & $0.8275$ \\
         & $3\%$ & $1.3144 \cdot 10^{-2}$ & $9.3043\cdot 10^{-4}$  & $0.7627$ & $0.7586$ \\
         & $5\%$ & $1.3882\cdot 10^{-2}$ &  $9.4602\cdot 10^{-4}$ & $0.6722$  & $0.7497$ \\
        \midrule
        \multirow{3}{*}{Mod-DOT - FC \small{(MSE+L1)}} & $1\%$ & $2.3437 \cdot 10^{-3}$  & $1.4226 \cdot 10^{-4}$  & $0.7925$  & $0.9592$ \\
         & $3\%$ & $3.2026\cdot 10^{-3}$  & $2.2142 \cdot 10^{-4}$  & $0.7152$ &  $0.9424$  \\
         & $5\%$ & $4.2008 \cdot 10^{-3}$ & $2.8598 \cdot 10^{-4}$  & $0.6988$  &  $0.9217$ \\
         \midrule
        \multirow{3}{*}{E2E - CONV \small{(TOT)}} & $1\%$ & $ 4.55500\cdot 10^{-3}$ & $2.82388\cdot 10^{-4}$ & $0.7054$ & $0.9287$ \\
         & $3\%$ &$ 1.55794\cdot 10^{-2}$  & $1.08184\cdot 10^{-3}$ &  $0.5339$ & $0.7707$ \\
         & $5\%$ &$ 1.41335\cdot 10^{-2}$ & $1.08453\cdot 10^{-3}$ & $0.6160$ & $0.7624$ \\
        \midrule
        \multirow{3}{*}{Mod-DOT - CONV \small{(MSE)}} & $1\%$ &$ 2.17595\cdot 10^{-3}$ & $ 1.42636\cdot 10^{-4}$ & $0.81388$ & $0.96252$ \\
         & $3\%$ &$  2.83007\cdot 10^{-3}$ & $1.99531\cdot 10^{-4}$ & $ 0.73764$ &  $0.94974$ \\
         & $5\%$ & $ 3.25266\cdot 10^{-3}$ &$2.41429\cdot 10^{-4}$ & $0.71242$ &  $0.94322$ \\
         \bottomrule
    \end{tabular}
    \caption{Performance metrics for the different architectures we consider. For each method it is
    specified the loss function which performs better. FC=fully connected, CONV=convolutional layers.}
    \label{tab:indici_medi_NOISE}
\end{table}

\subsubsection{Comparison of Mod-DOT vs E2E architectures}
\begin{figure}[htbp]
\centering
\includegraphics[width=\textwidth]{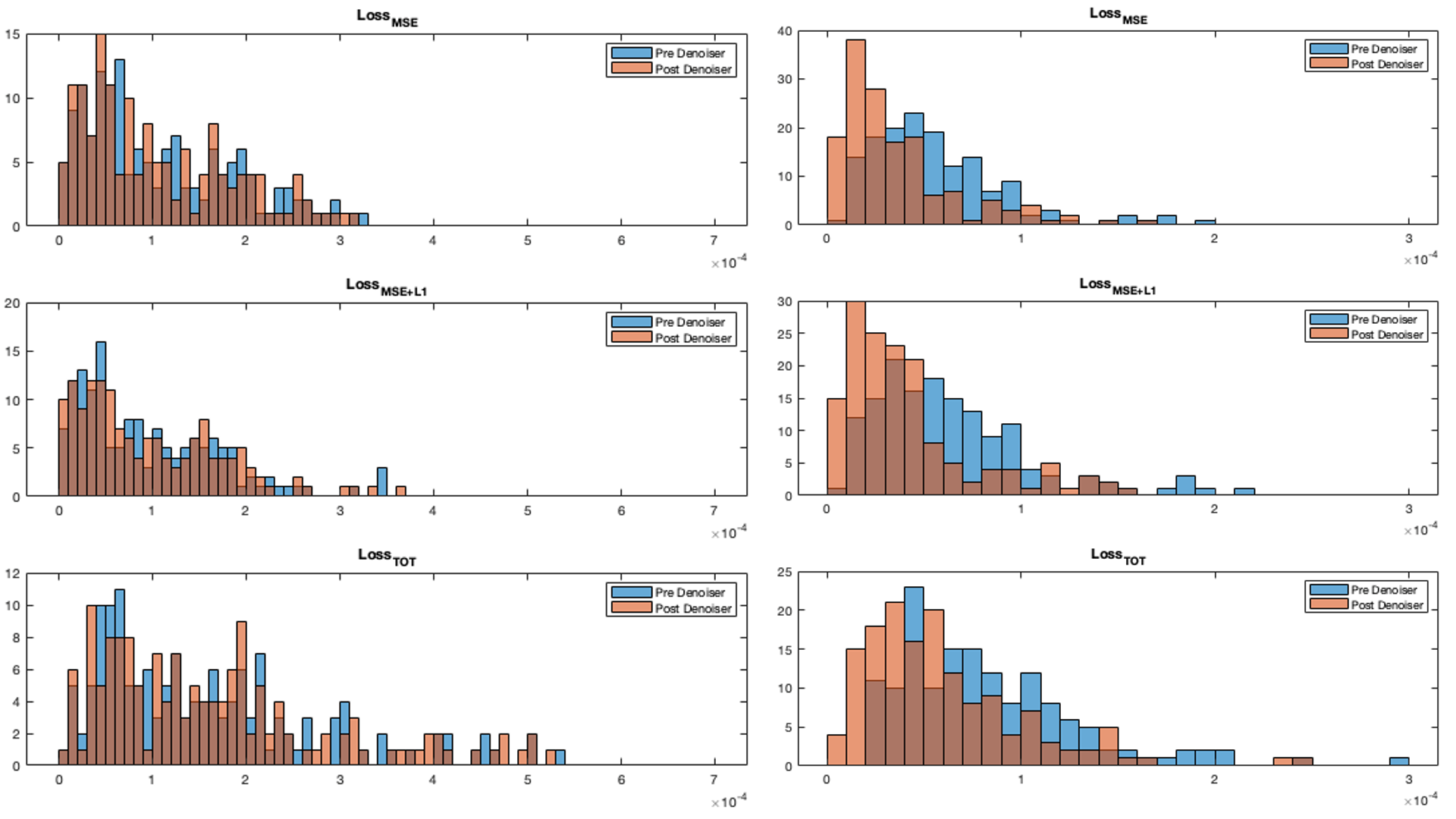}
\caption{Histogram of the MSE obtained from 150 test samples. The reconstruction is obtained via the 
E2E (left column) vs Mod-DOT
architecture (right column) with noiseless data and after passage through the denoising network.
The superposed histograms in each
figure represent the MSE   
before (light blue bars) and after (orange bars) passing through the denoising network. 
The rows correspond to the different choices of the loss function as discussed in
Sect.~\ref{subsec:TrainingAlg}.
Note that the ranges of the graphs in each column 
are different for better readability.}
\label{fig:bars_noiseless}
\end{figure}

\begin{figure}[htbp]
\centering
\includegraphics[width=\textwidth]{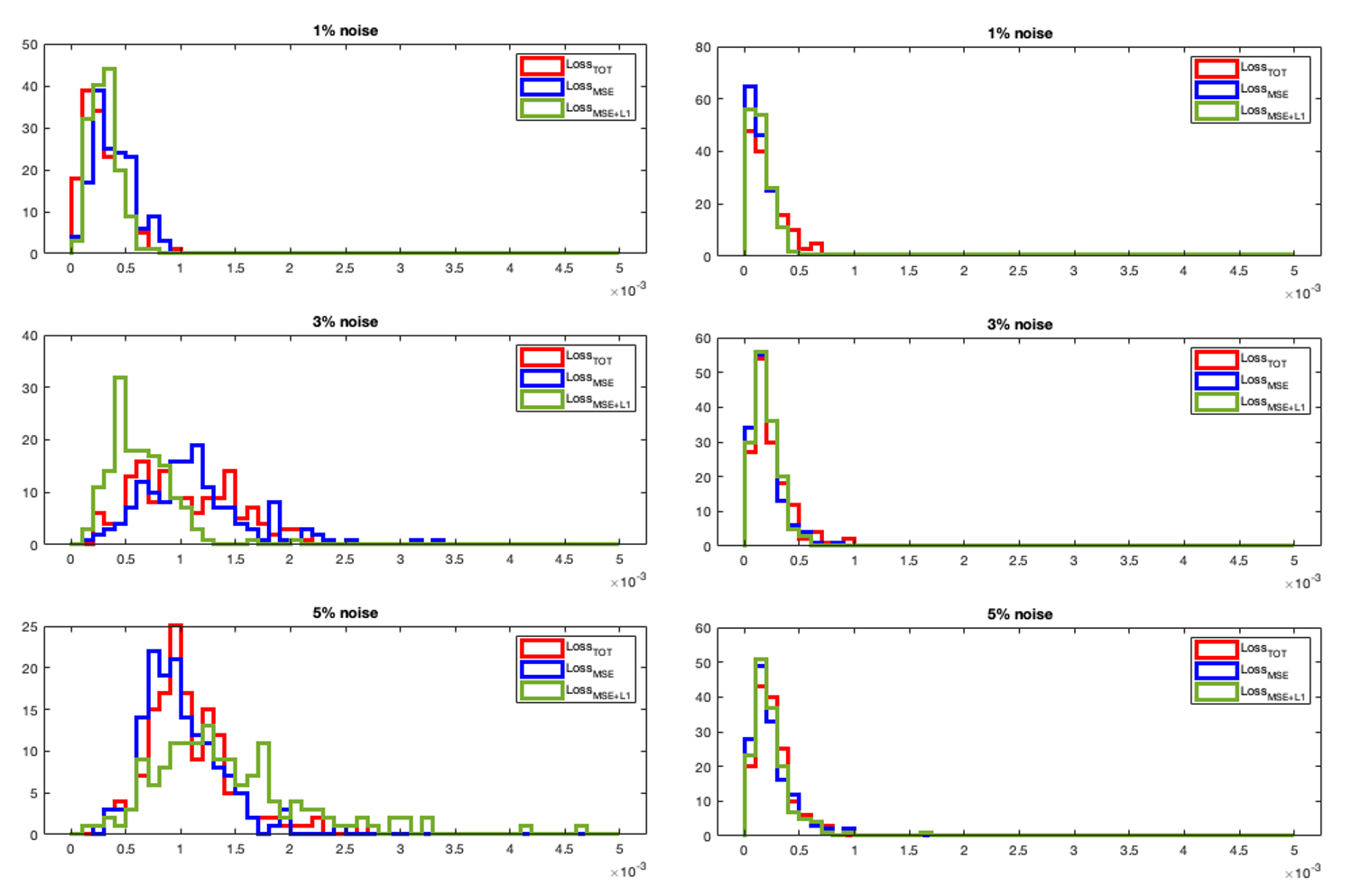}
\caption{Histogram of the MSE obtained from 150 test samples. The reconstruction is obtained via the 
E2E (left column) vs Mod-DOT
architecture (right column) with noisy data and after passage through the denoising network. 
The rows correspond to the different noise levels, while colors in each
graph correspond to different choices of the loss function.
Note that the ranges of the graphs
are different for better readability.}
\label{fig:bars_noise}
\end{figure}

We compare in this study the Mod-DOT and the corresponding E2E 
architecture. Notice that the E2E architecture represents
a typical end-to-end neural network-based workflow for problem inversion. 
In the context of DOT reconstruction, it is comparable to the architectures
used in the recent literature works~\cite{feng2019back,yoo2019deep,deng2021deep}. We consider the same domain as in the previous section and
we assess the performance of the two architectures based on 150 test (unseen) samples. 
We first consider noiseless data, then we consider noisy data as 
in the previous study. 
Noticeably, the E2E architecture
outperforms the 
Mod-DOT when we use noiseless data
(see Fig.~\ref{fig:bars_noiseless}),
while the Mod-DOT architecture is
superior when considering noisy, more realistic, data (see Fig.~\ref{fig:bars_noise}).

\begin{figure*}[htbp]
\newcommand{\factor}{0.18}
\begin{center}
	\subfloat[\label{fig:gt_E2E}]{\includegraphics[width=\factor\textwidth]{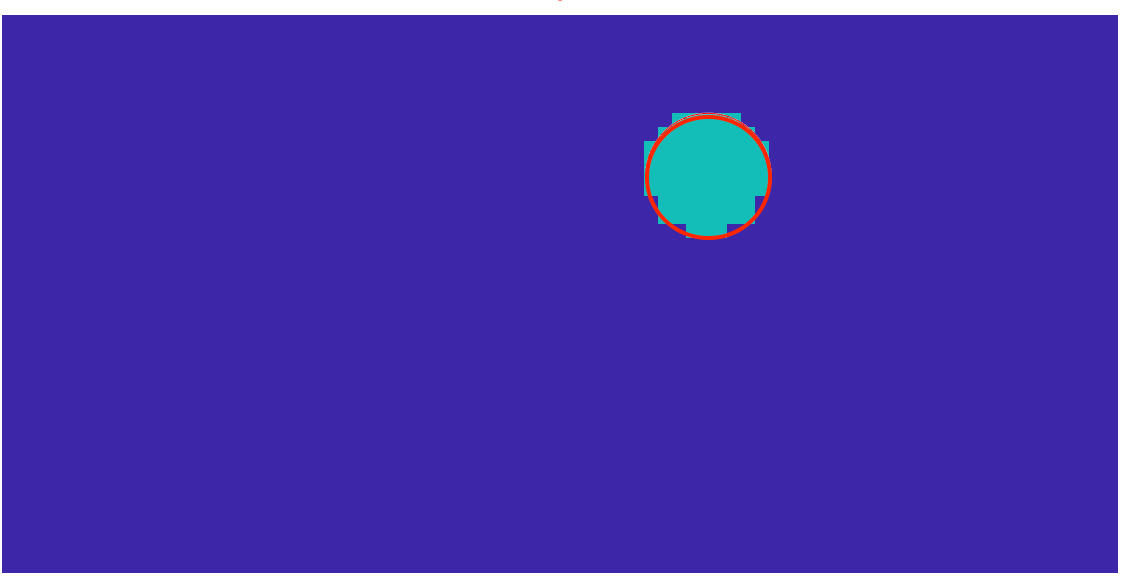}\hspace{0.025\textwidth}\includegraphics[width=\factor\textwidth]{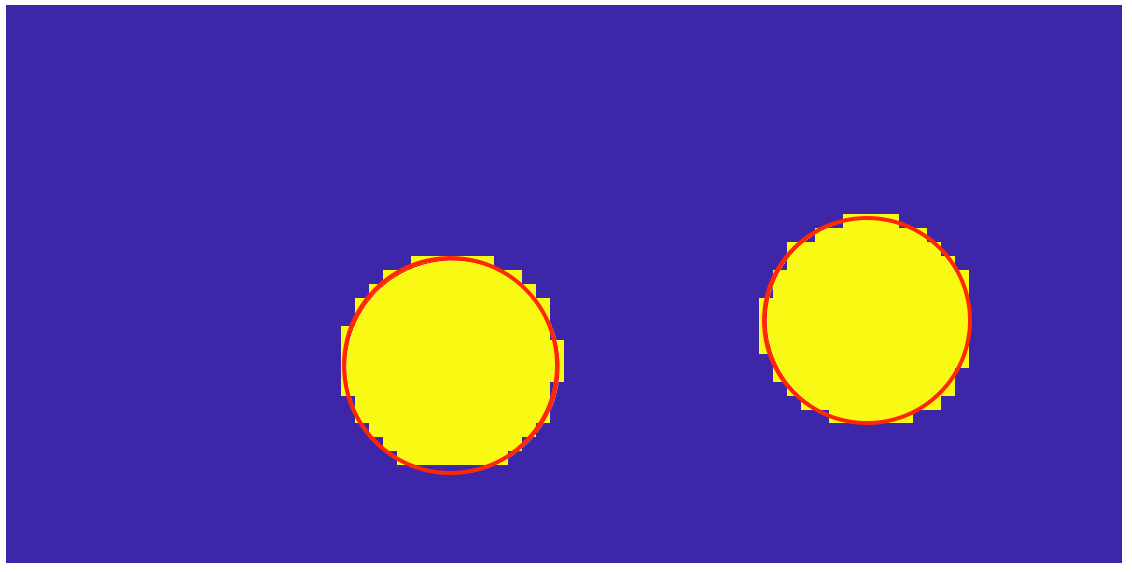}\hspace{0.025\textwidth}\includegraphics[width=\factor\textwidth]{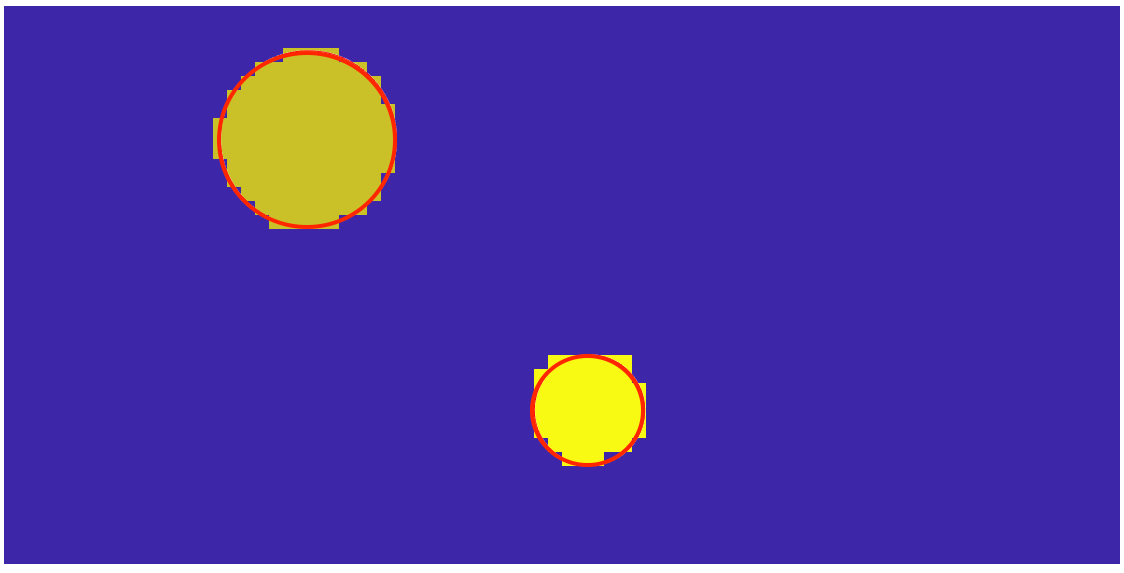}\hspace{0.025\textwidth}\includegraphics[width=\factor\textwidth]{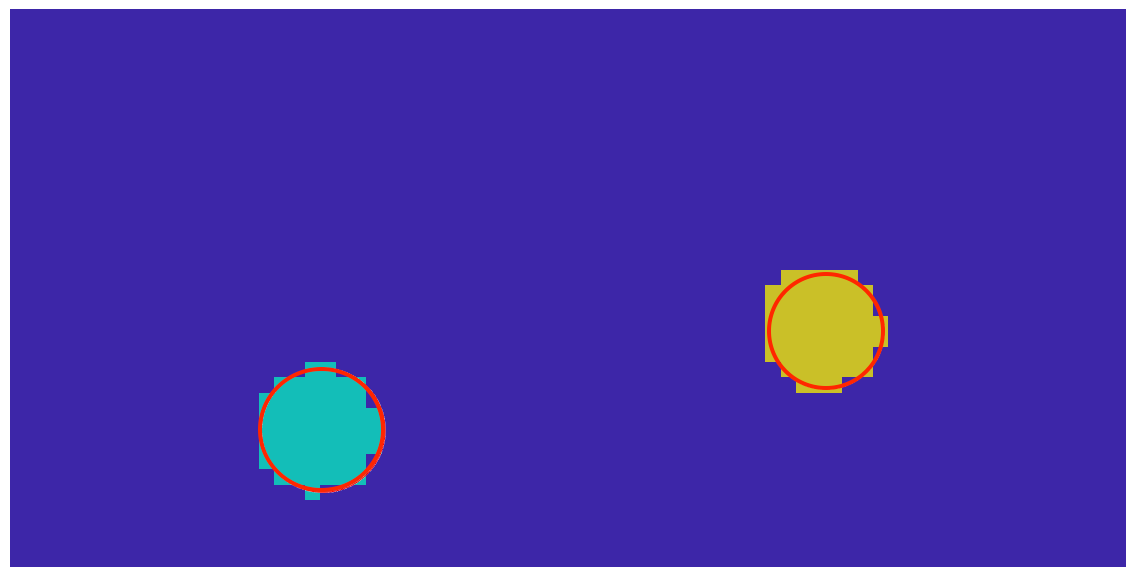}\hspace{0.025\textwidth}\includegraphics[width=\factor\textwidth]{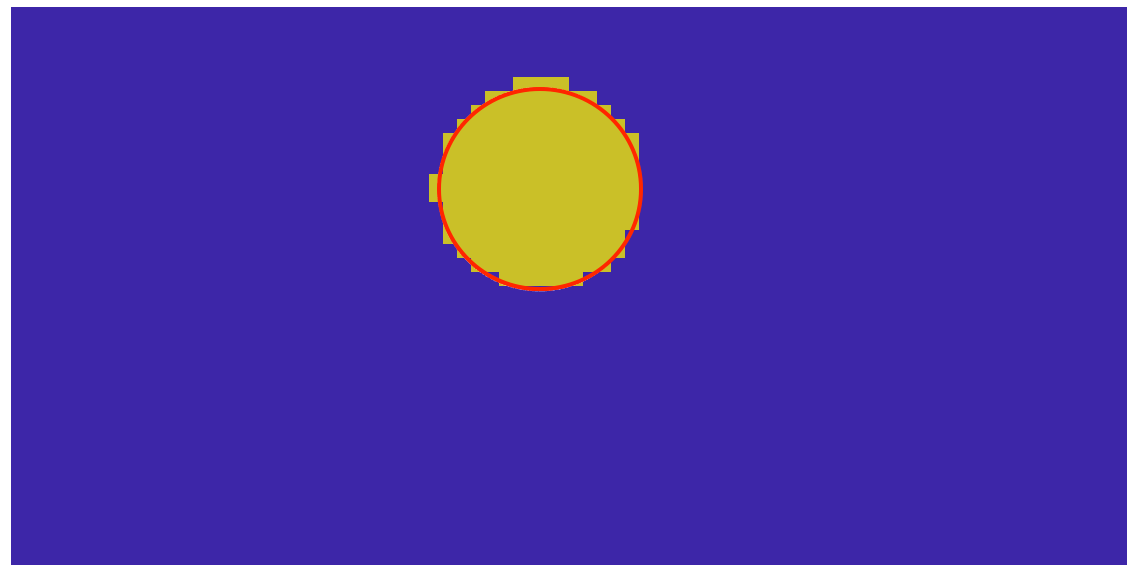}}
	
	\subfloat[\label{fig:PreDEN_E2E}]{\includegraphics[width=\factor\textwidth]{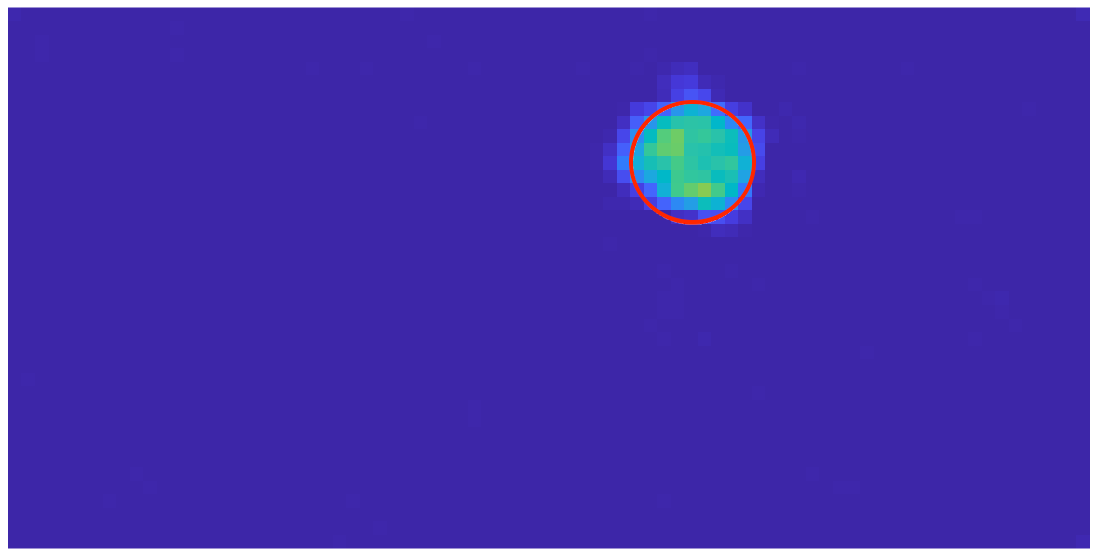}\hspace{0.025\textwidth}\includegraphics[width=\factor\textwidth]{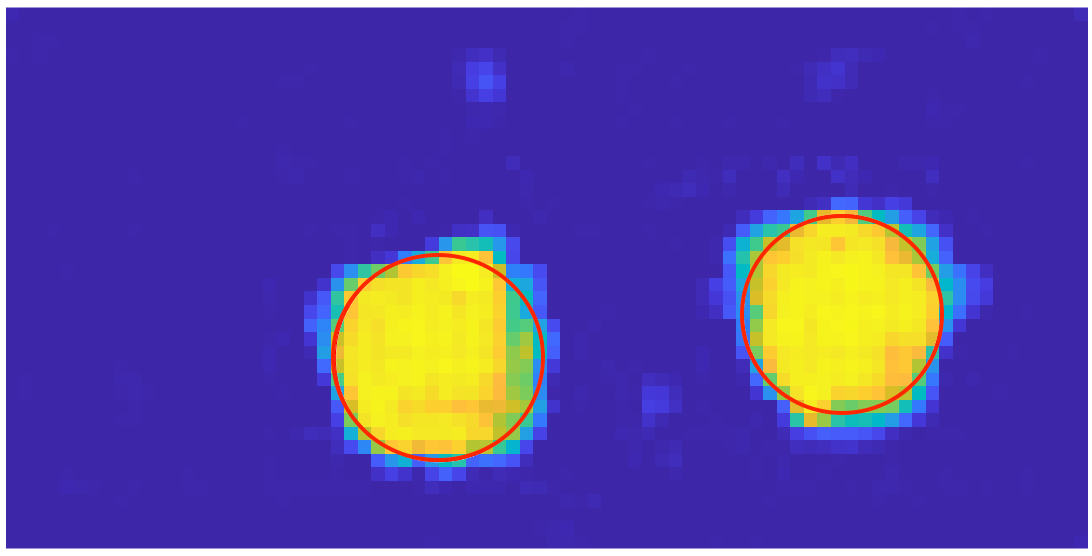}\hspace{0.025\textwidth}\includegraphics[width=\factor\textwidth]{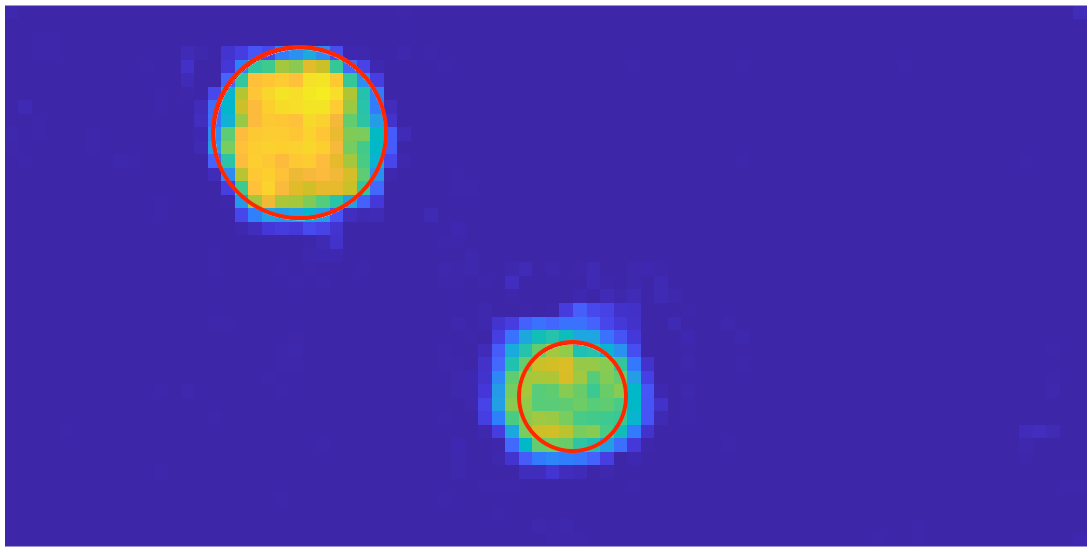}\hspace{0.025\textwidth}\includegraphics[width=\factor\textwidth]{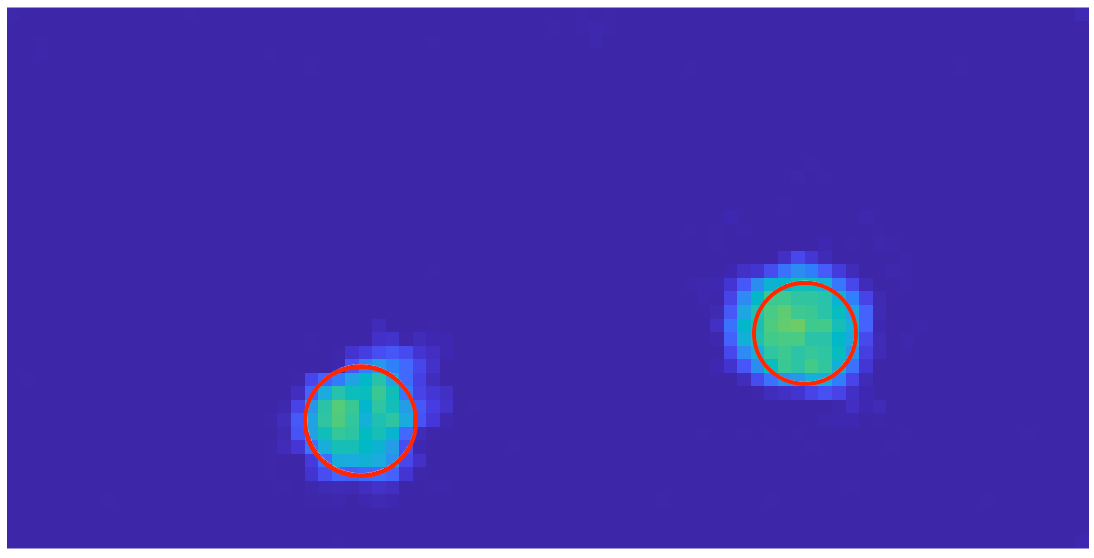}\hspace{0.025\textwidth}\includegraphics[width=\factor\textwidth]{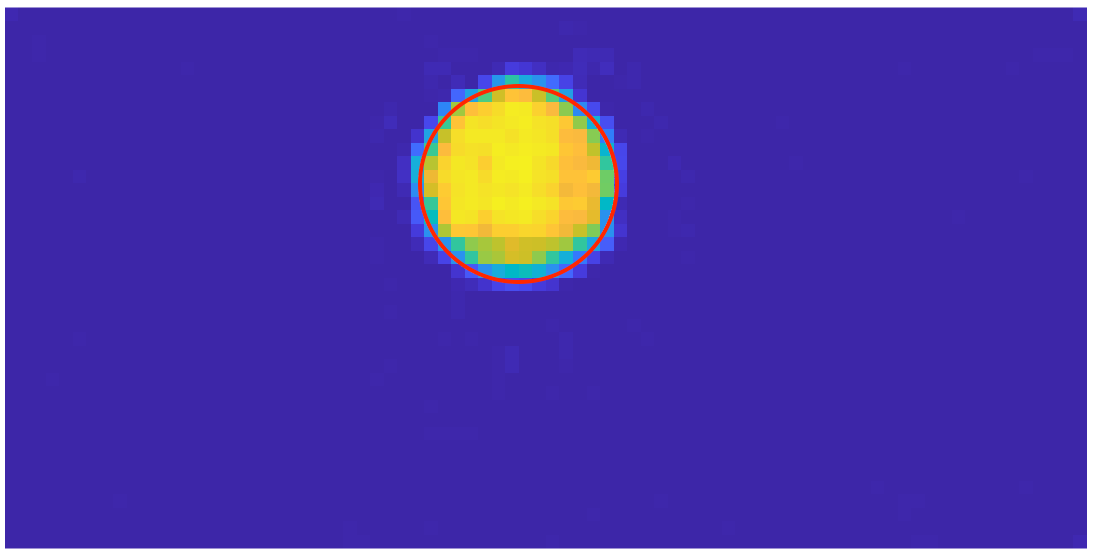}}
	
	\subfloat[\label{fig:PostDEN_E2E}]{\includegraphics[width=\factor\textwidth]{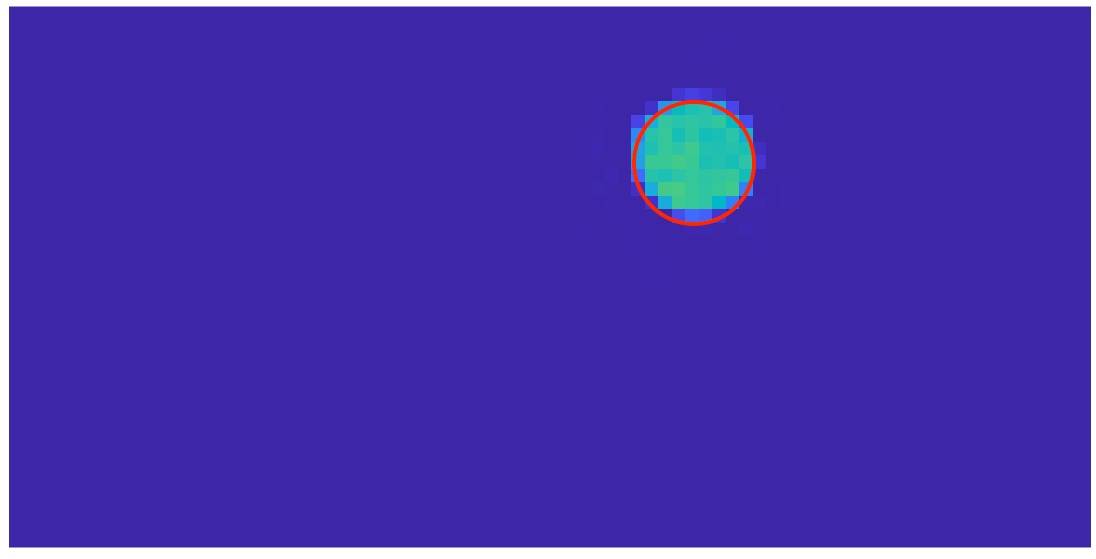}\hspace{0.025\textwidth}\includegraphics[width=\factor\textwidth]{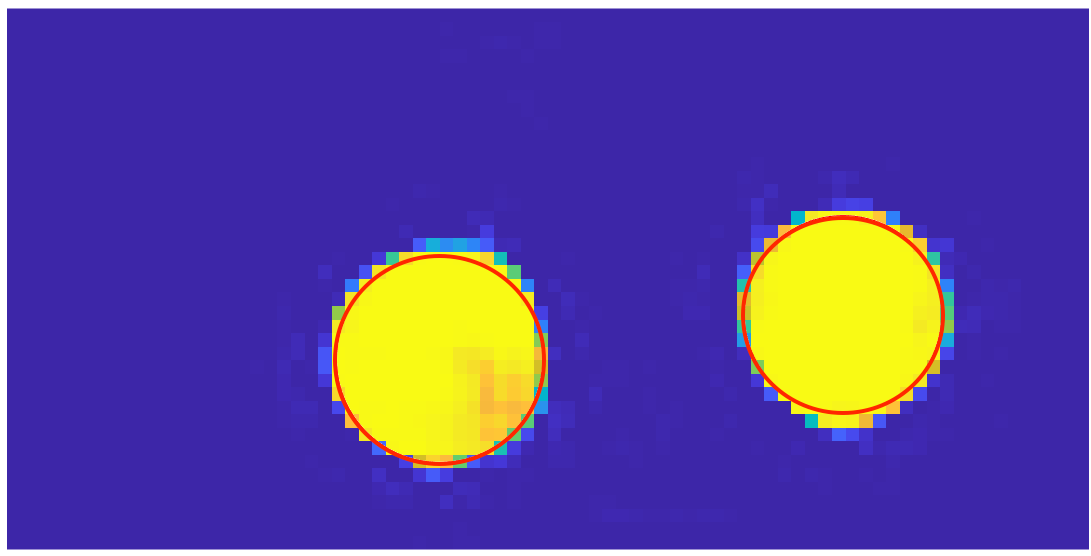}\hspace{0.025\textwidth}\includegraphics[width=\factor\textwidth]{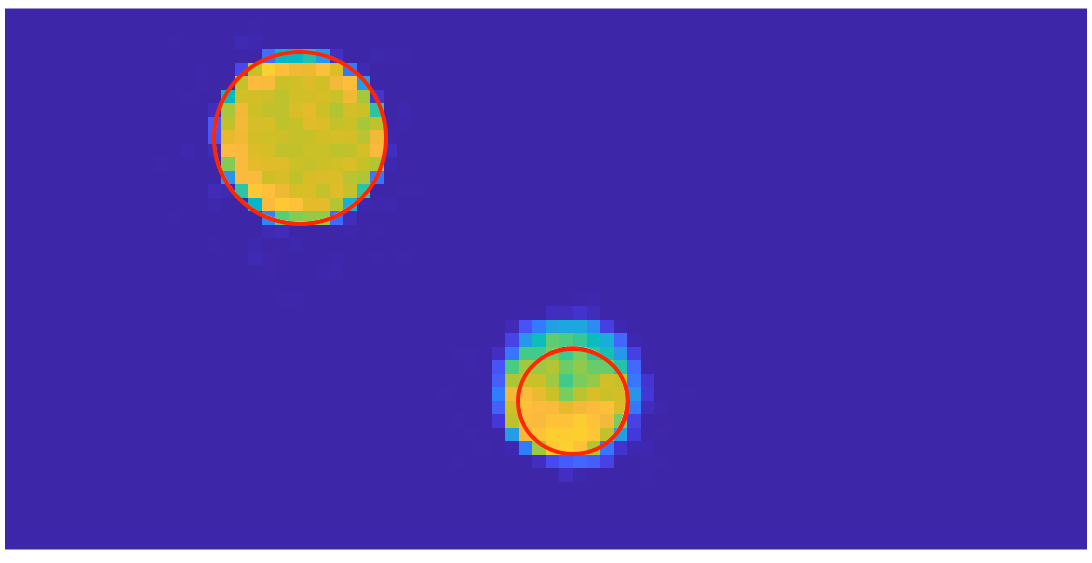}\hspace{0.025\textwidth}\includegraphics[width=\factor\textwidth]{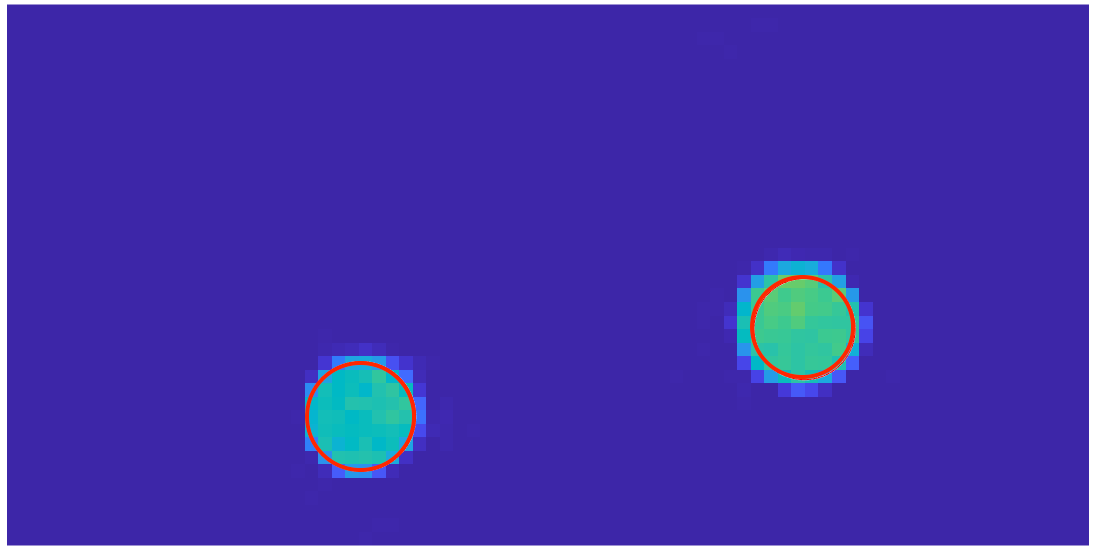}\hspace{0.025\textwidth}\includegraphics[width=\factor\textwidth]{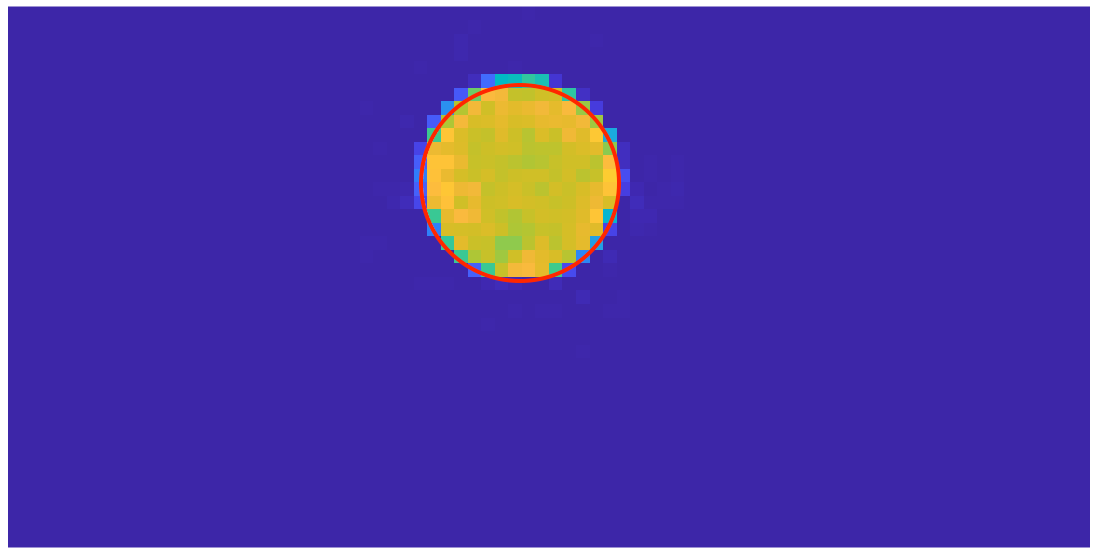}}
		
	\subfloat[\label{fig:G1_E2E}]{\includegraphics[width=\factor\textwidth]{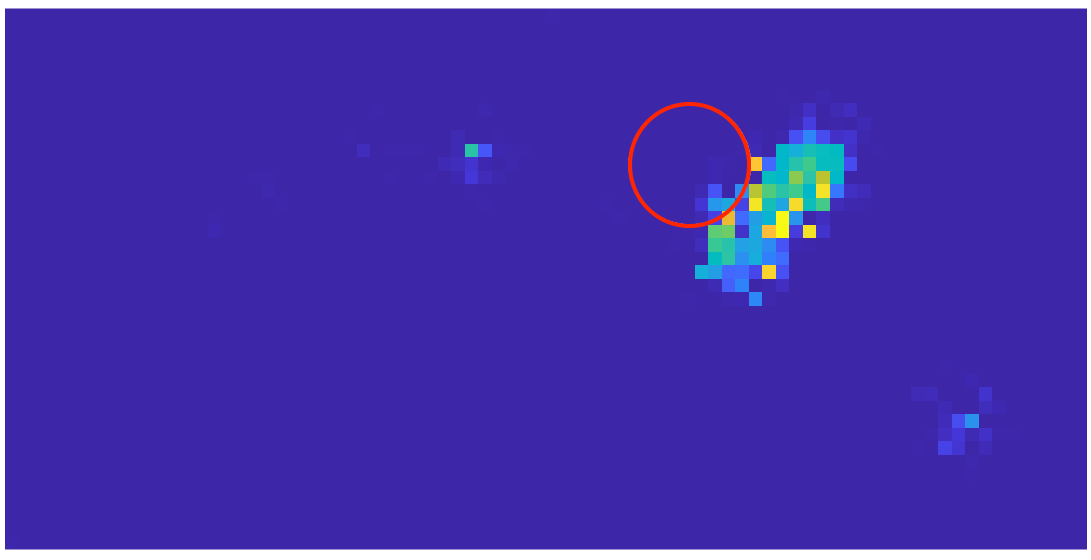}\hspace{0.025\textwidth}\includegraphics[width=\factor\textwidth]{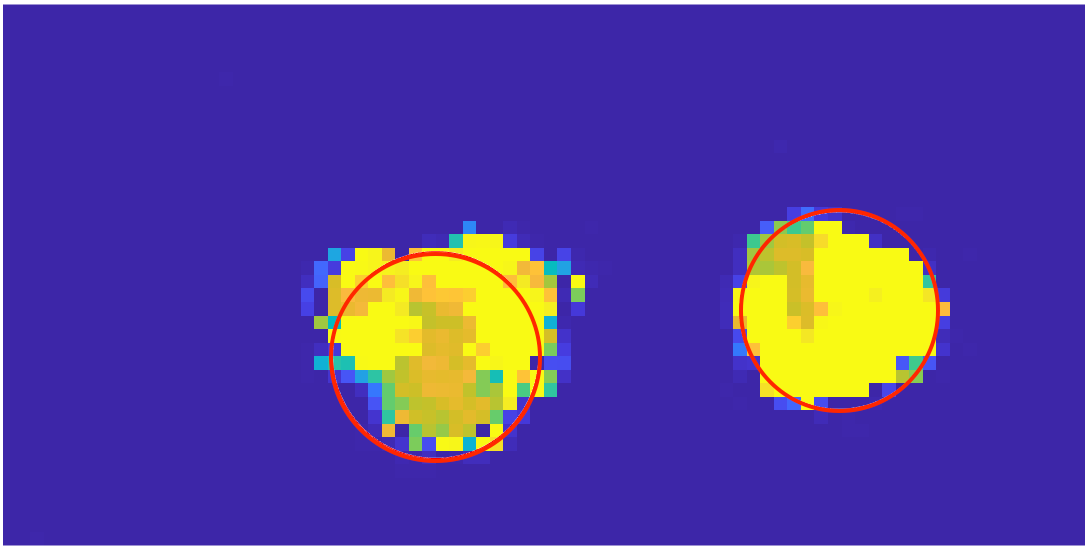}\hspace{0.025\textwidth}\includegraphics[width=\factor\textwidth]{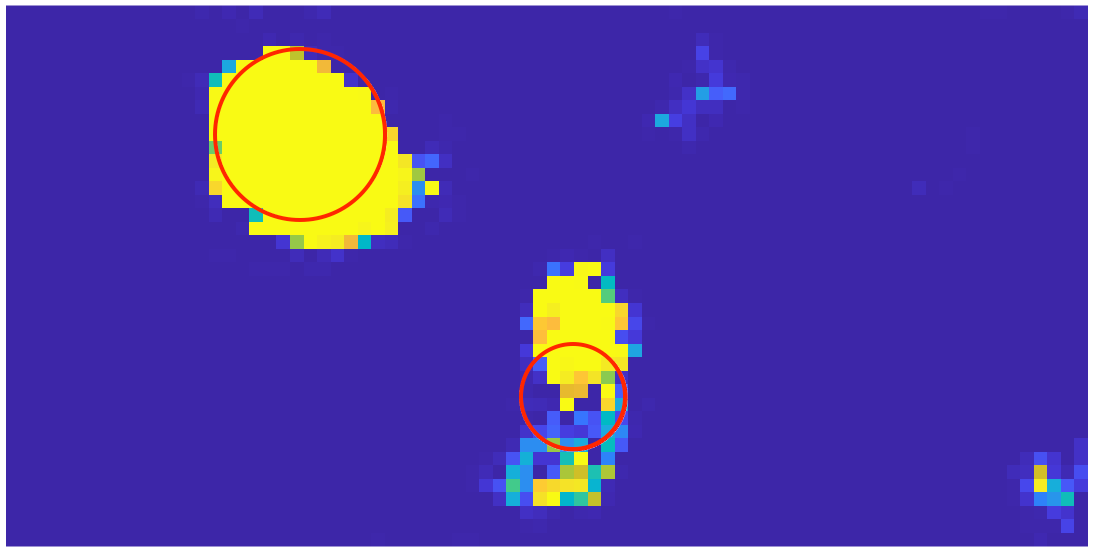}\hspace{0.025\textwidth}\includegraphics[width=\factor\textwidth]{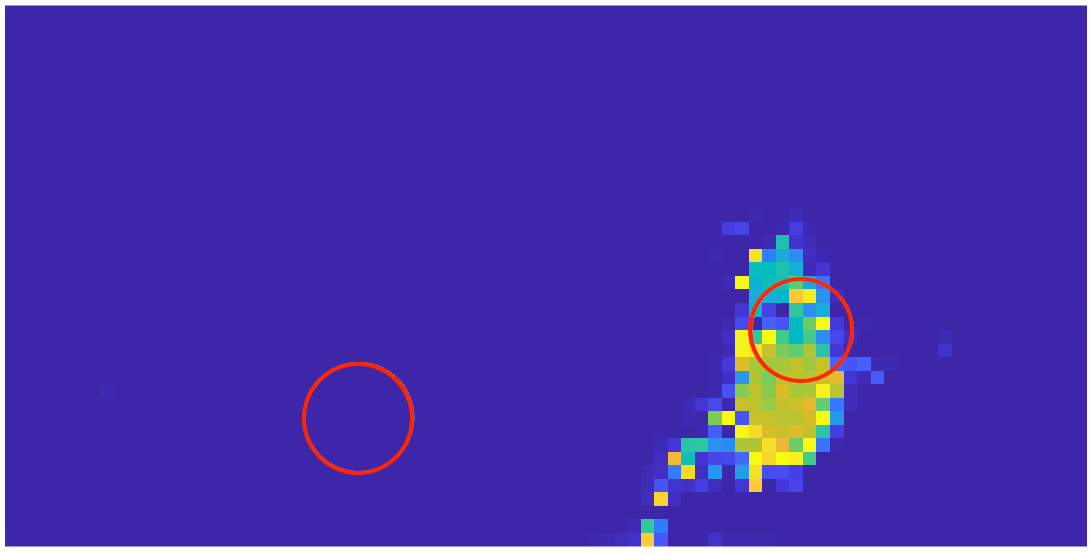}\hspace{0.025\textwidth}\includegraphics[width=\factor\textwidth]{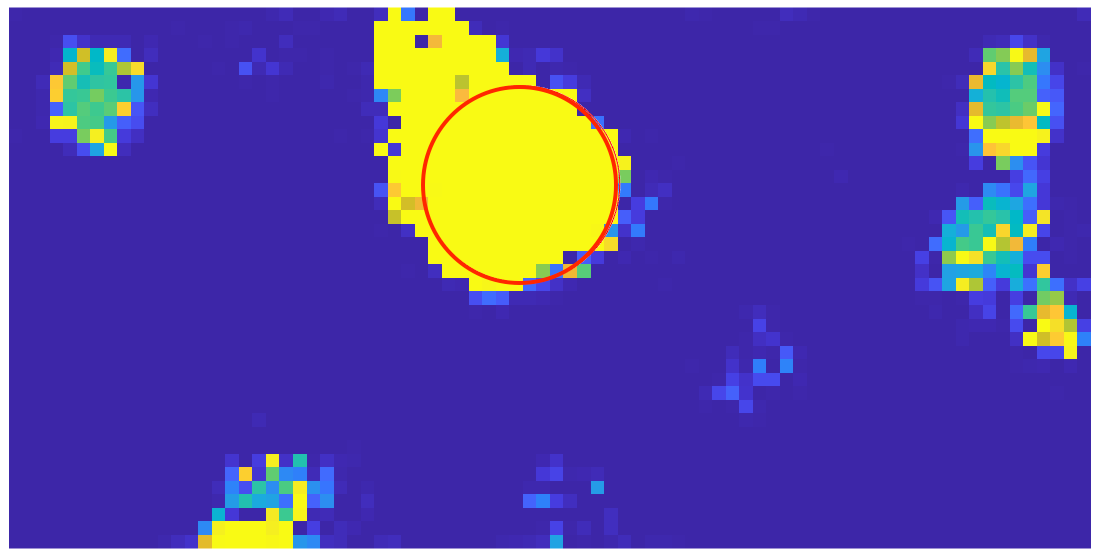}}
	
    \subfloat[\label{fig:G3_E2E}]{\includegraphics[width=\factor\textwidth]{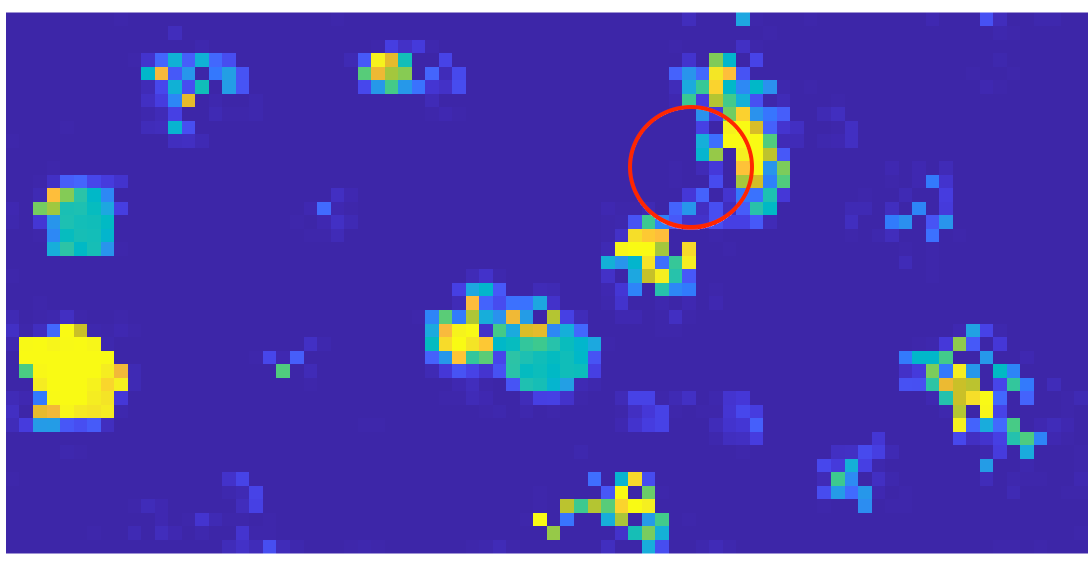}\hspace{0.025\textwidth}\includegraphics[width=\factor\textwidth]{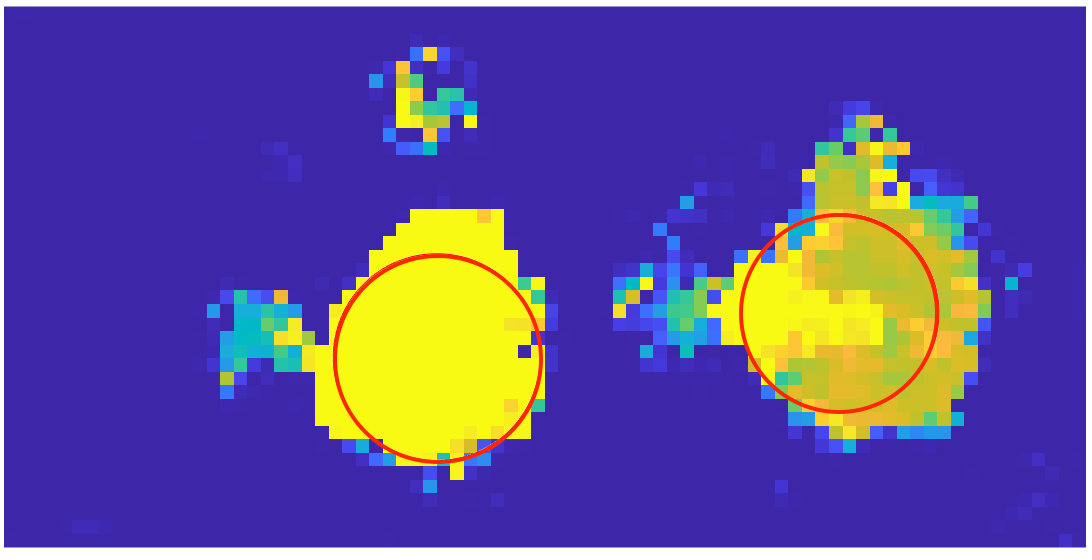}\hspace{0.025\textwidth}\includegraphics[width=\factor\textwidth]{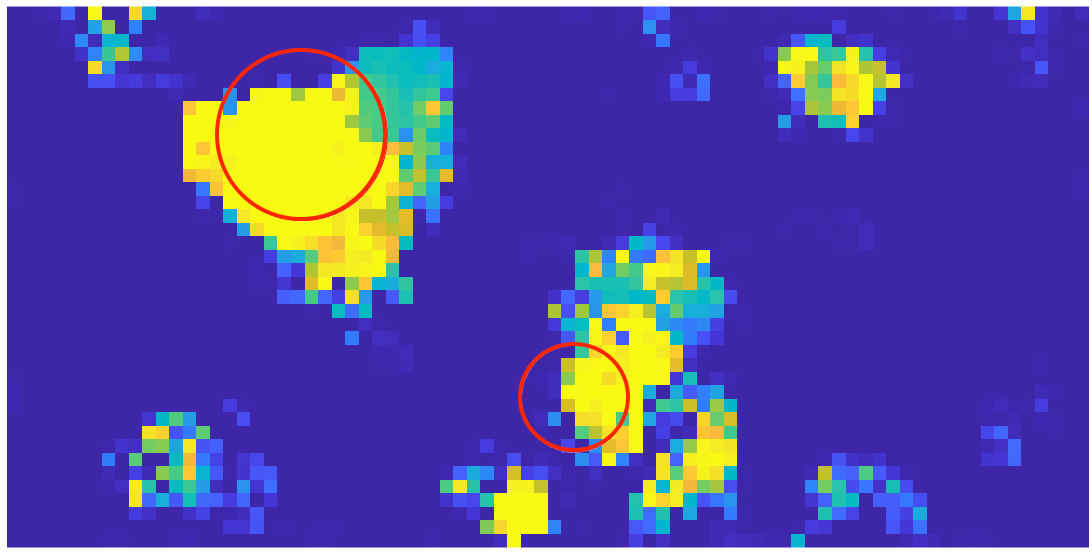}\hspace{0.025\textwidth}\includegraphics[width=\factor\textwidth]{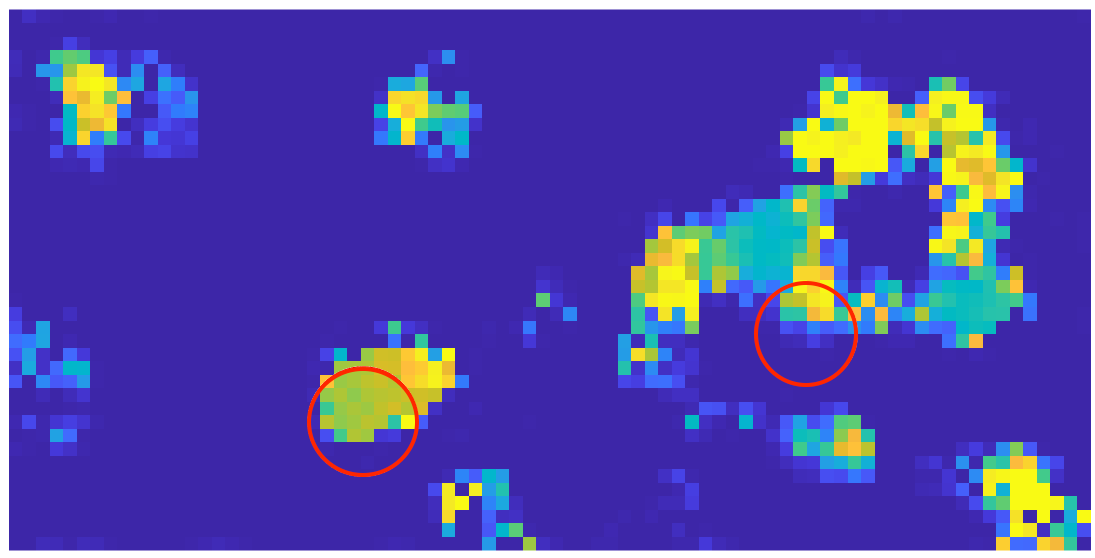}\hspace{0.025\textwidth}\includegraphics[width=\factor\textwidth]{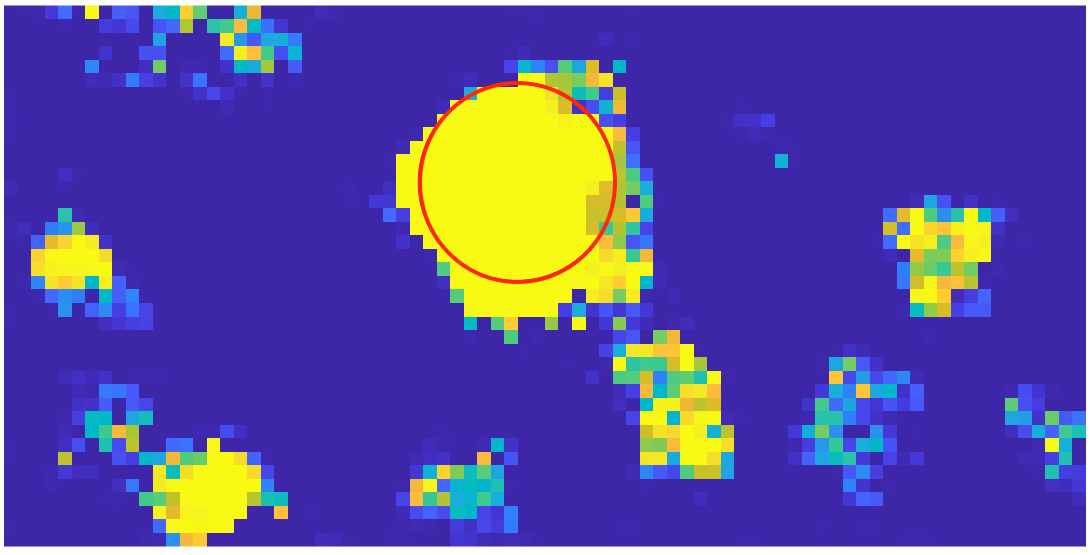}}

    \subfloat[\label{fig:G5_E2E}]{\includegraphics[width=\factor\textwidth]{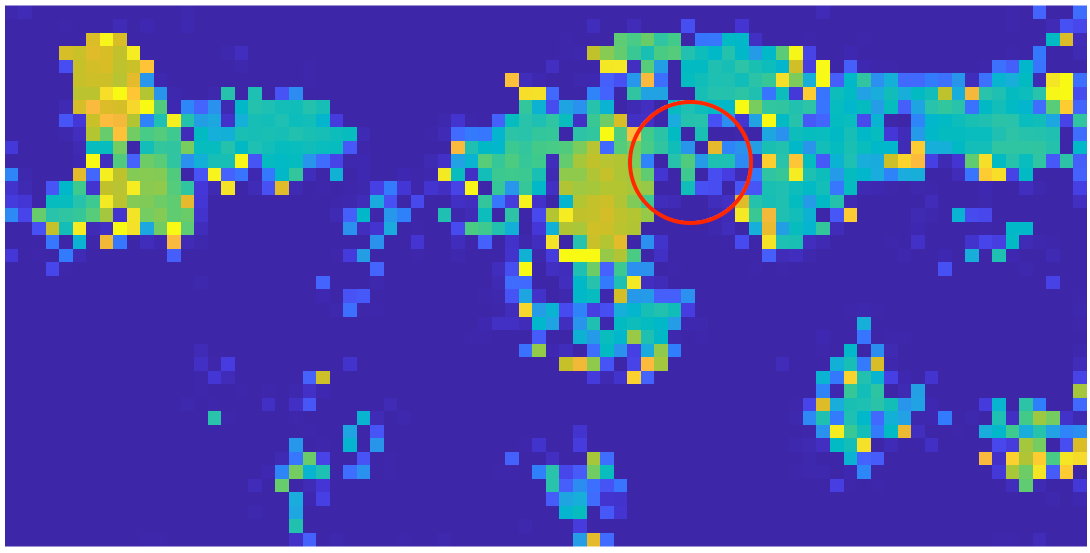}\hspace{0.025\textwidth}\includegraphics[width=\factor\textwidth]{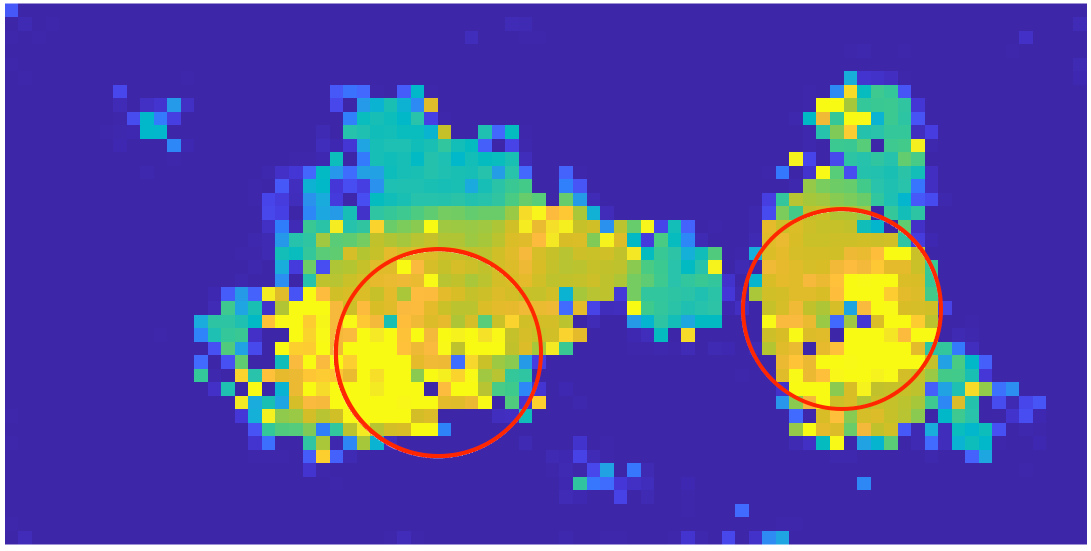}\hspace{0.025\textwidth}\includegraphics[width=\factor\textwidth]{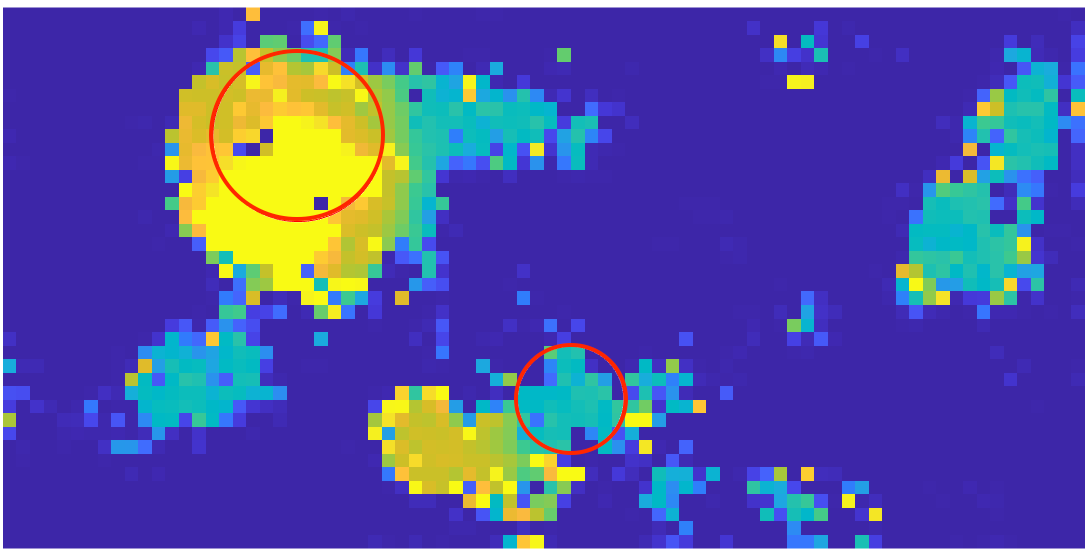}\hspace{0.025\textwidth}\includegraphics[width=\factor\textwidth]{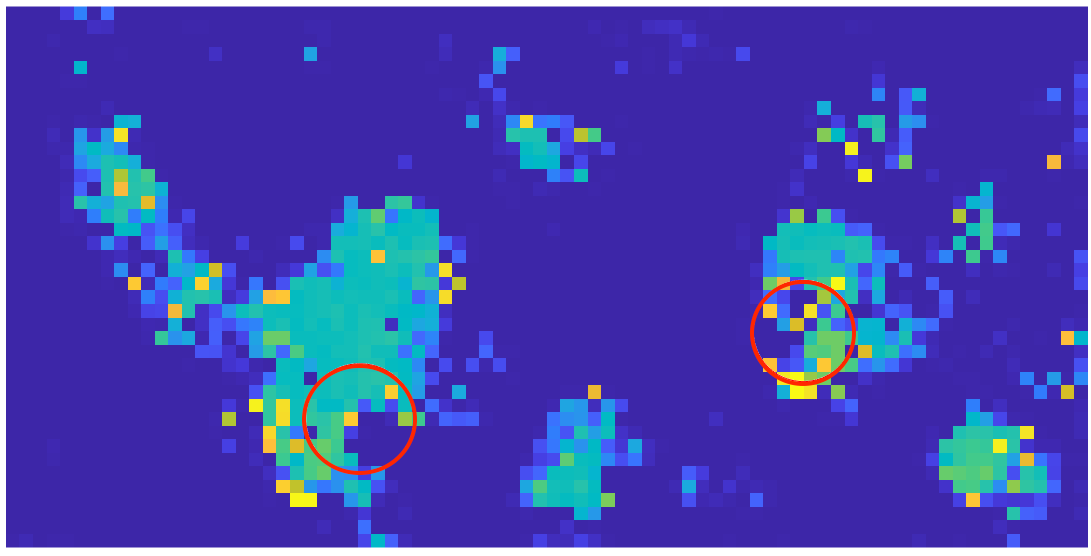}\hspace{0.025\textwidth}\includegraphics[width=\factor\textwidth]{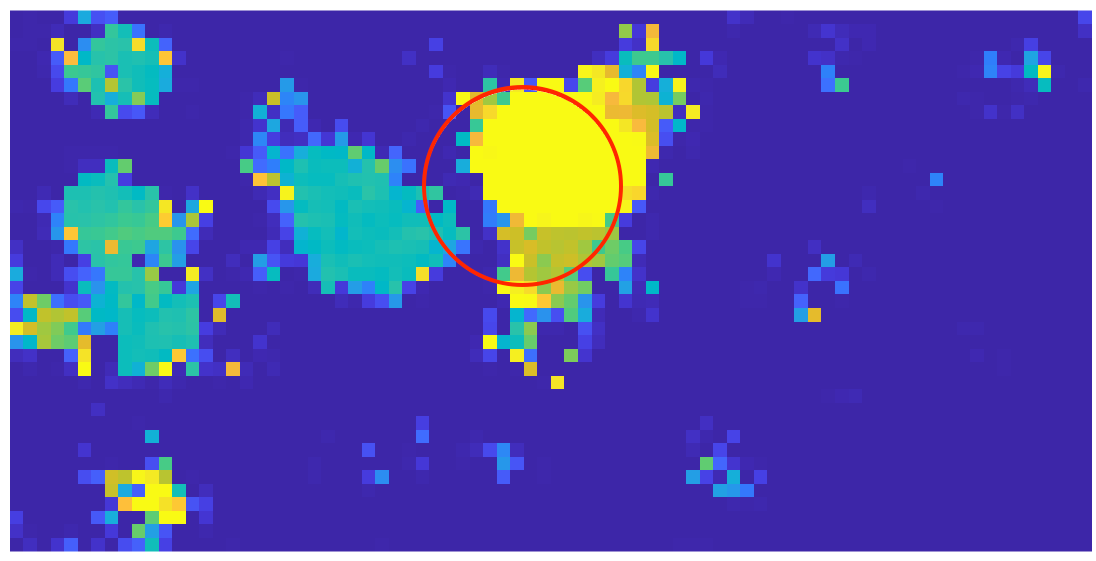}}
    
\subfloat{\includegraphics[width=1\textwidth]{Fig8C}}

\end{center}
\caption{Reconstruction of the absorption coefficient distribution obtained with the E2E approach. (a): ground truth; (b): reconstruction without pretraining the AEs, before the application of the denoiser; (c): reconstruction with pretraining of the AEs 
before denoising; (d): reconstruction with pretraining of the AEs after denoising.(e): 1\% noise; (f): 3\% noise; (g): 5\% noise.
Rows (e) to (g) use pretrain and denoising.
Values are presented as $\mu_a/D$, the 
red circles denote the true contrast regions.
}
\label{fig:fResults_E2E}
\end{figure*}

\begin{figure*}[htbp]
\newcommand{\factor}{0.18}
\begin{center}

\subfloat[\label{fig:gt_Br}]{\includegraphics[width=\factor\textwidth]{Fig13_A11}\hspace{0.025\textwidth}\includegraphics[width=\factor\textwidth]{Fig13_A12}\hspace{0.025\textwidth}\includegraphics[width=\factor\textwidth]{Fig13_A13}\hspace{0.025\textwidth}\includegraphics[width=\factor\textwidth]{Fig13_A14}\hspace{0.025\textwidth}\includegraphics[width=\factor\textwidth]{Fig13_A15}}

	\subfloat[\label{fig:noPre}]{\includegraphics[width=\factor\textwidth]{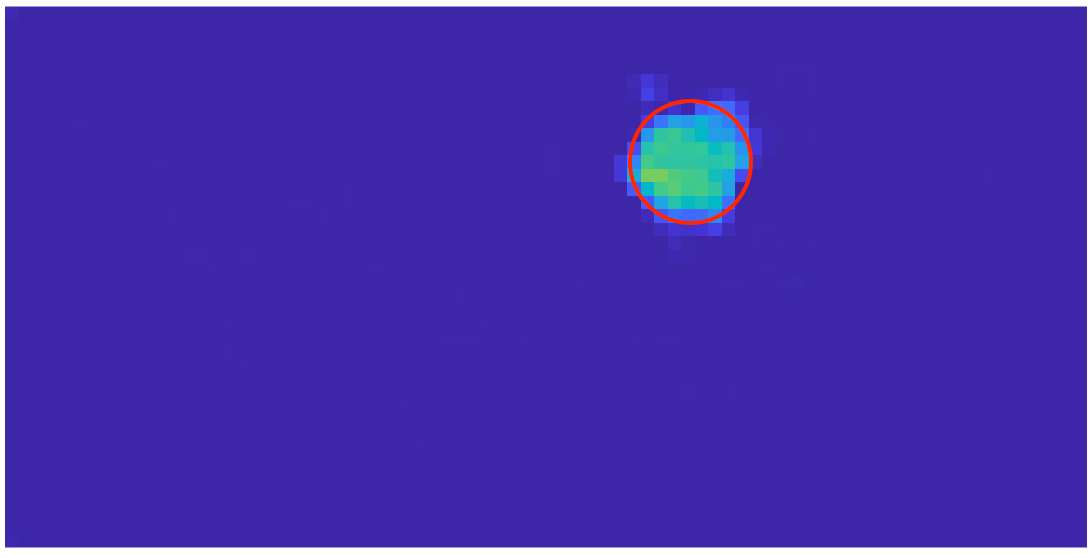}\hspace{0.025\textwidth}\includegraphics[width=\factor\textwidth]{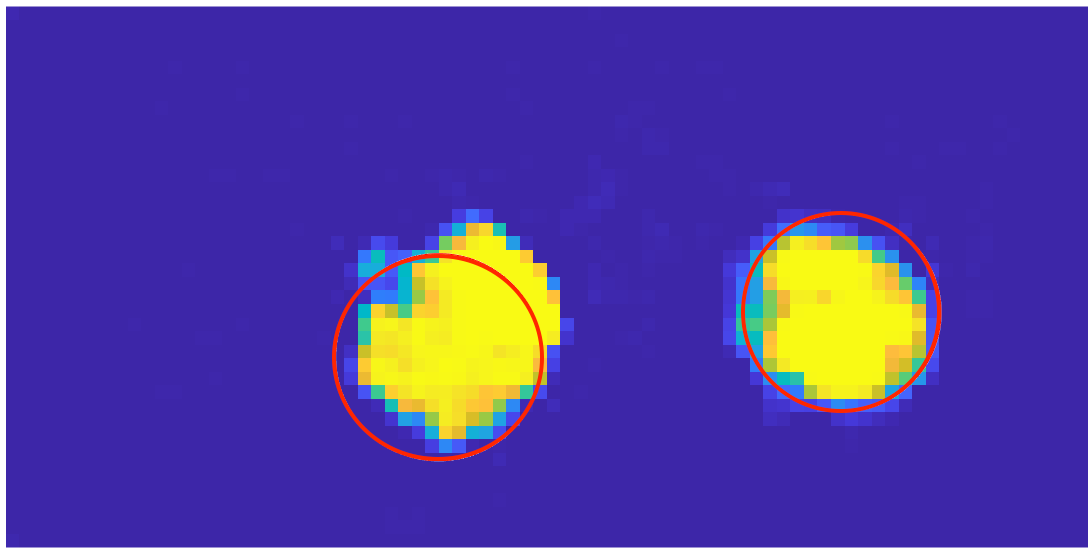}\hspace{0.025\textwidth}\includegraphics[width=\factor\textwidth]{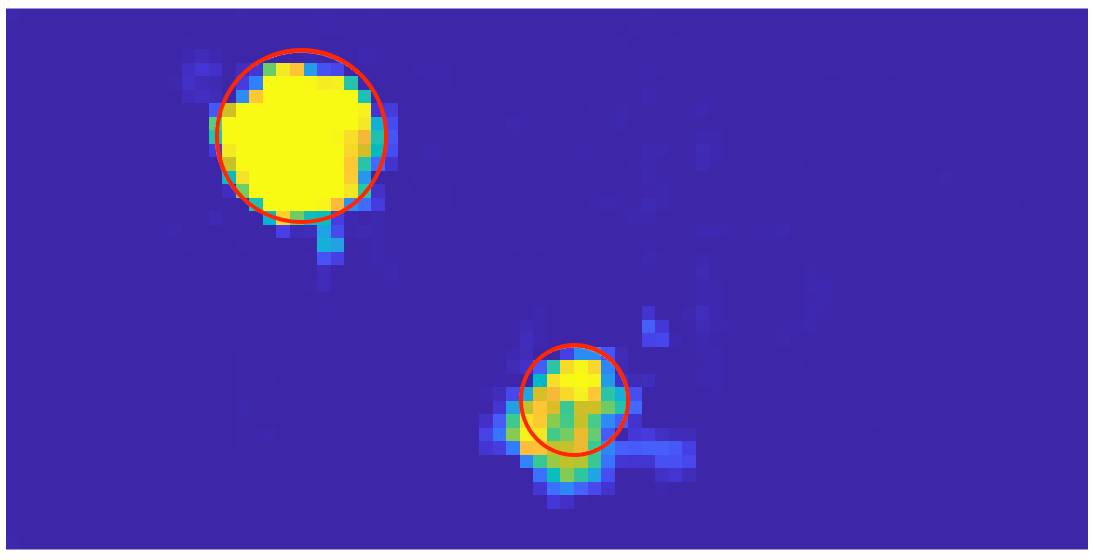}\hspace{0.025\textwidth}\includegraphics[width=\factor\textwidth]{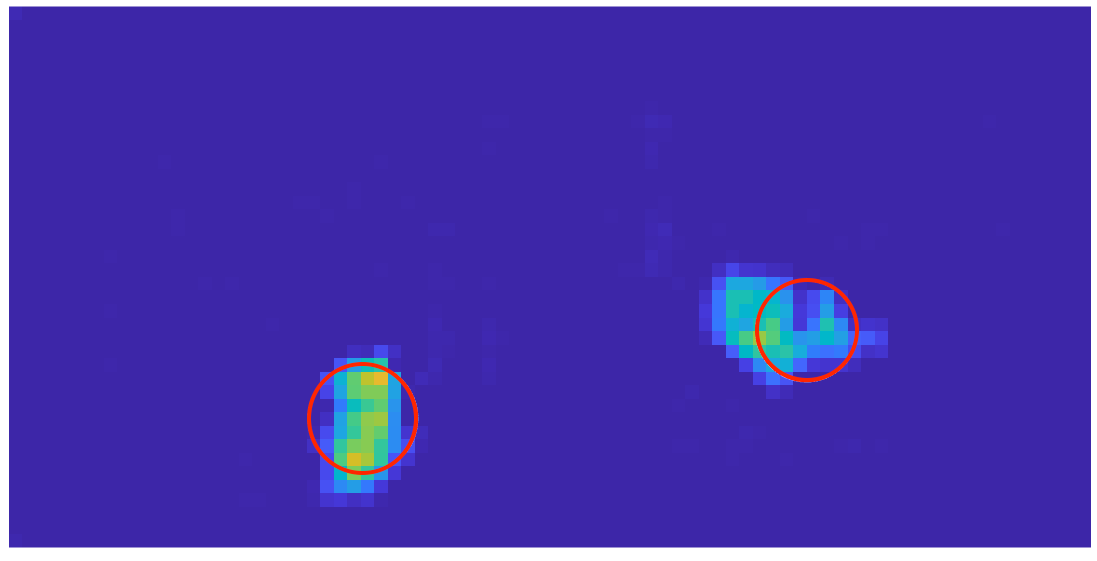}\hspace{0.025\textwidth}\includegraphics[width=\factor\textwidth]{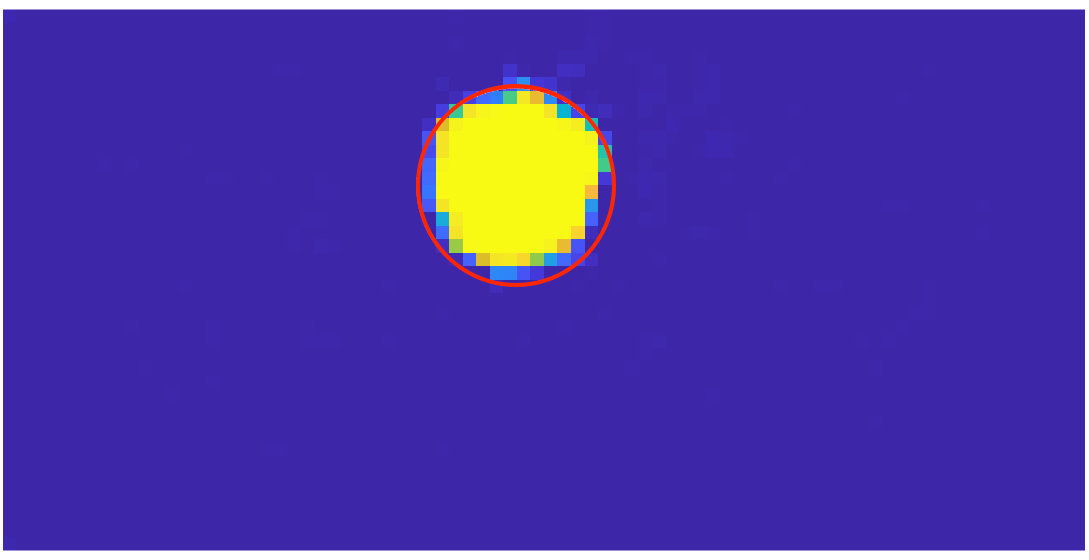}}
	
	\subfloat[\label{fig:PreDEN_DOT}]{\includegraphics[width=\factor\textwidth]{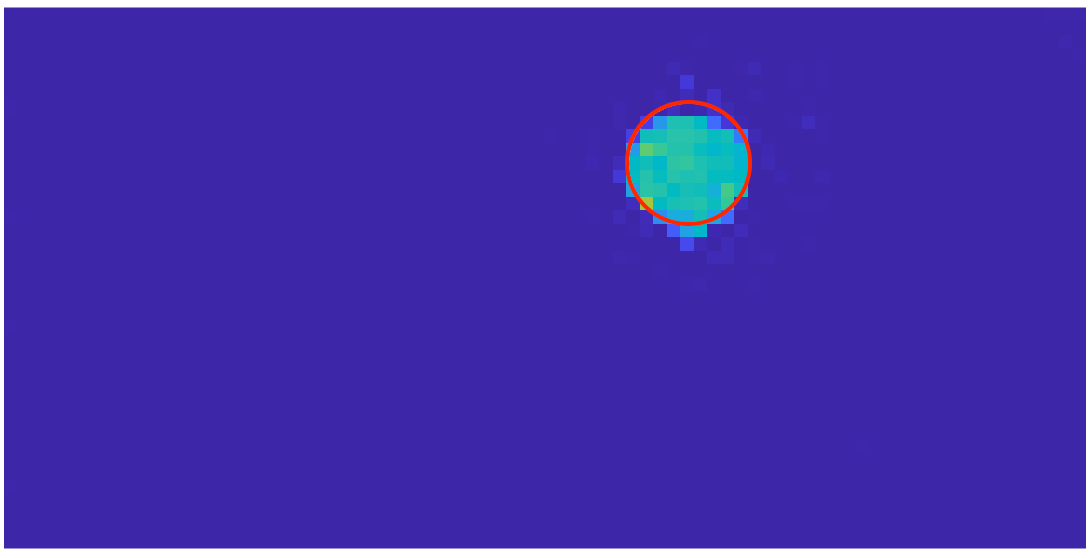}\hspace{0.025\textwidth}\includegraphics[width=\factor\textwidth]{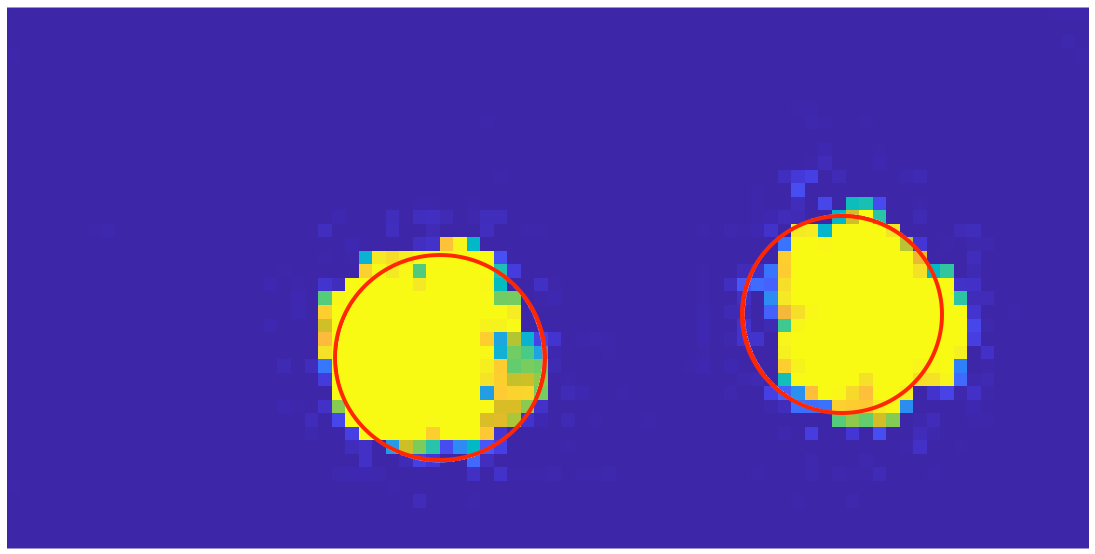}\hspace{0.025\textwidth}\includegraphics[width=\factor\textwidth]{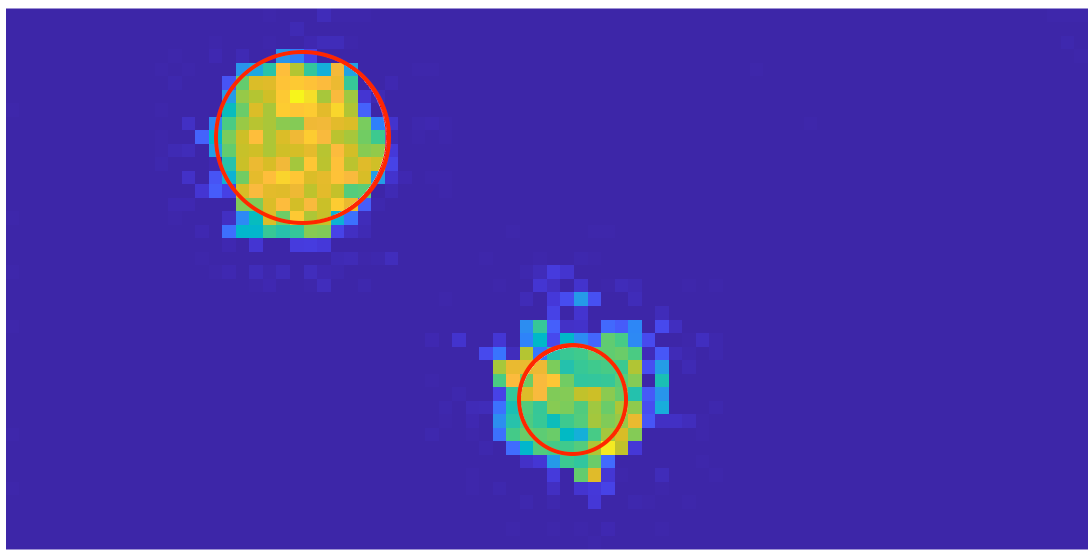}\hspace{0.025\textwidth}\includegraphics[width=\factor\textwidth]{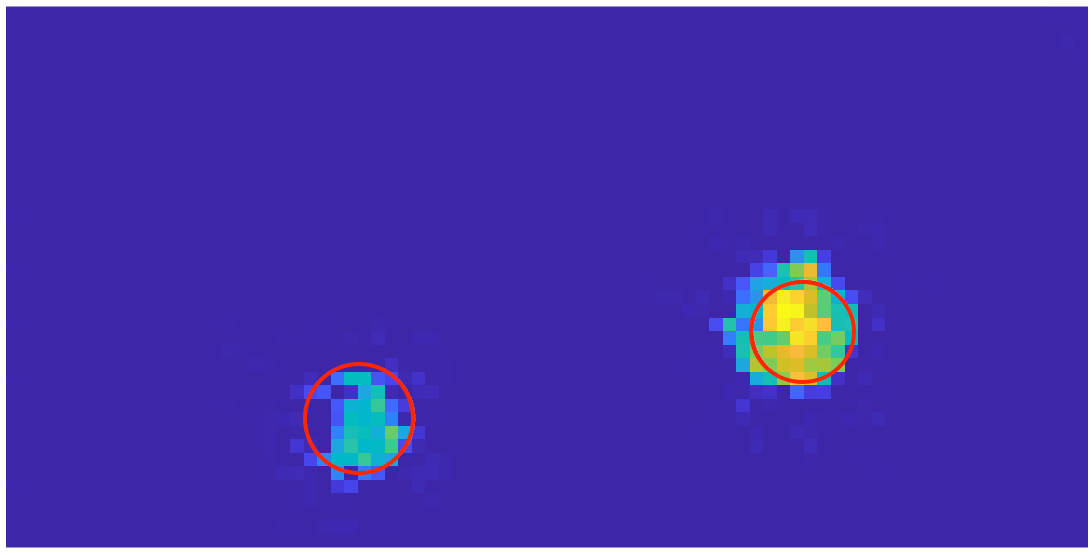}\hspace{0.025\textwidth}\includegraphics[width=\factor\textwidth]{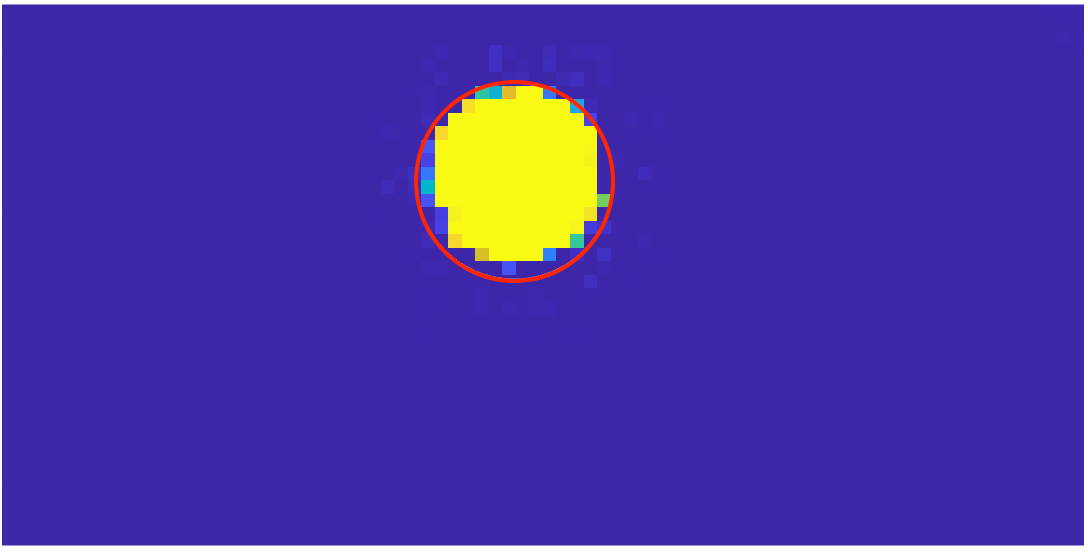}}
	
	\subfloat[\label{fig:PostDEN_DOT}]{\includegraphics[width=\factor\textwidth]{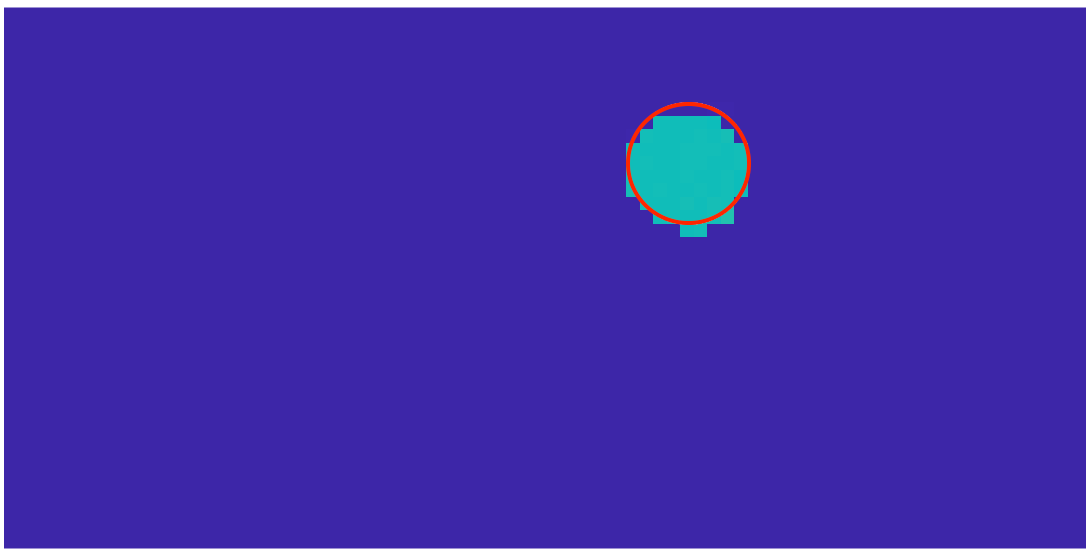}\hspace{0.025\textwidth}\includegraphics[width=\factor\textwidth]{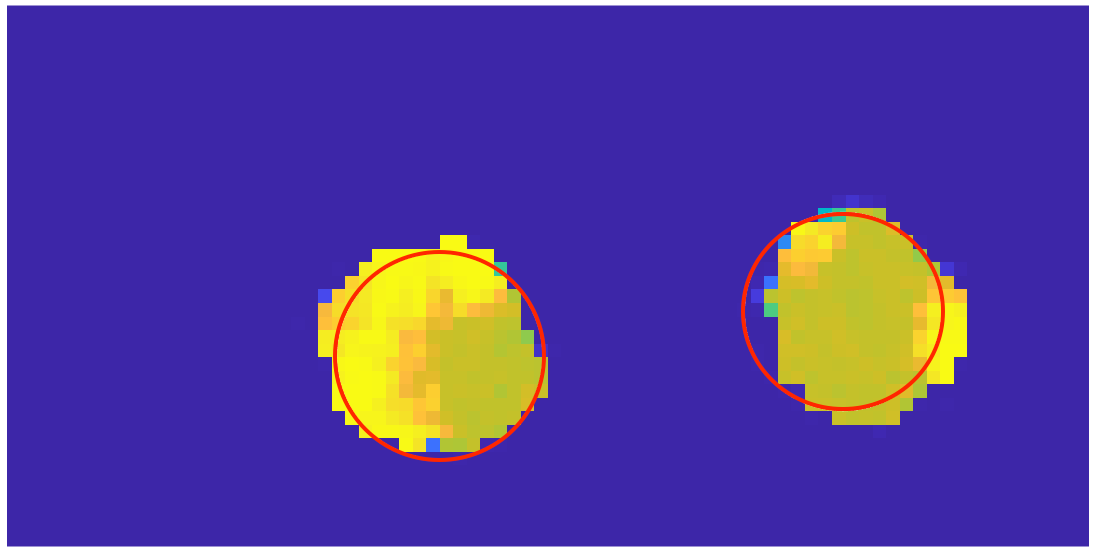}\hspace{0.025\textwidth}\includegraphics[width=\factor\textwidth]{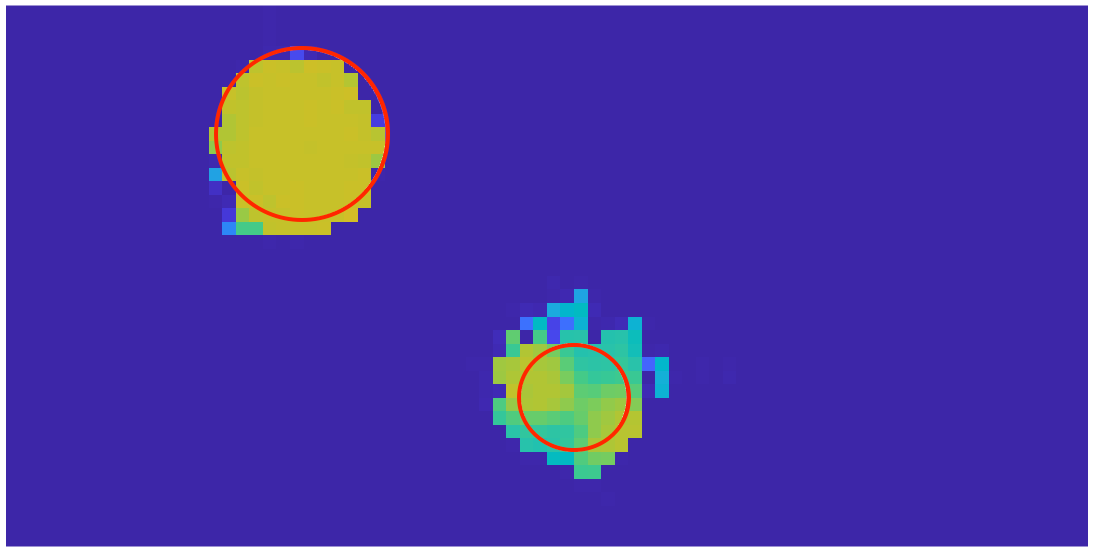}\hspace{0.025\textwidth}\includegraphics[width=\factor\textwidth]{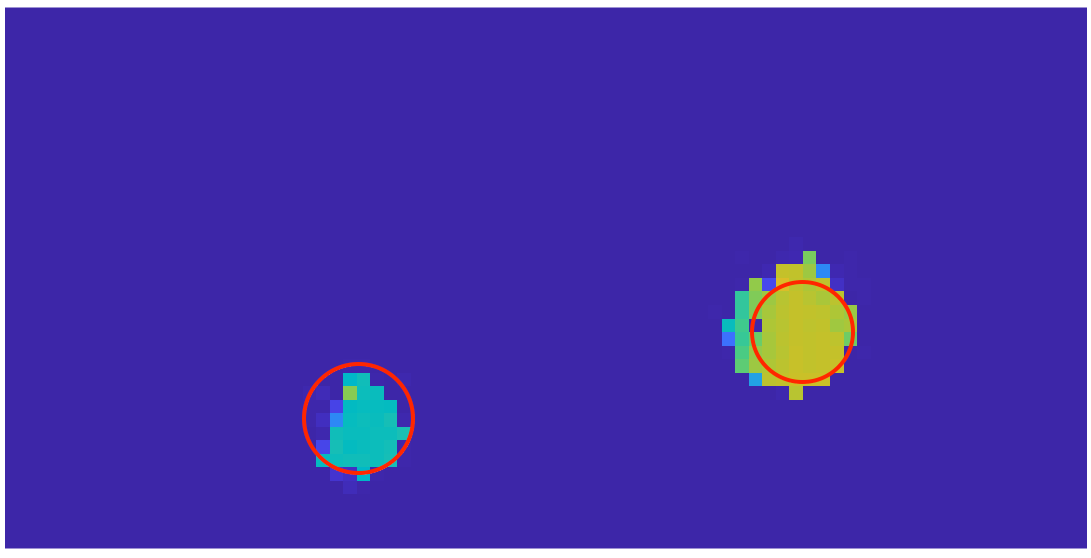}\hspace{0.025\textwidth}\includegraphics[width=\factor\textwidth]{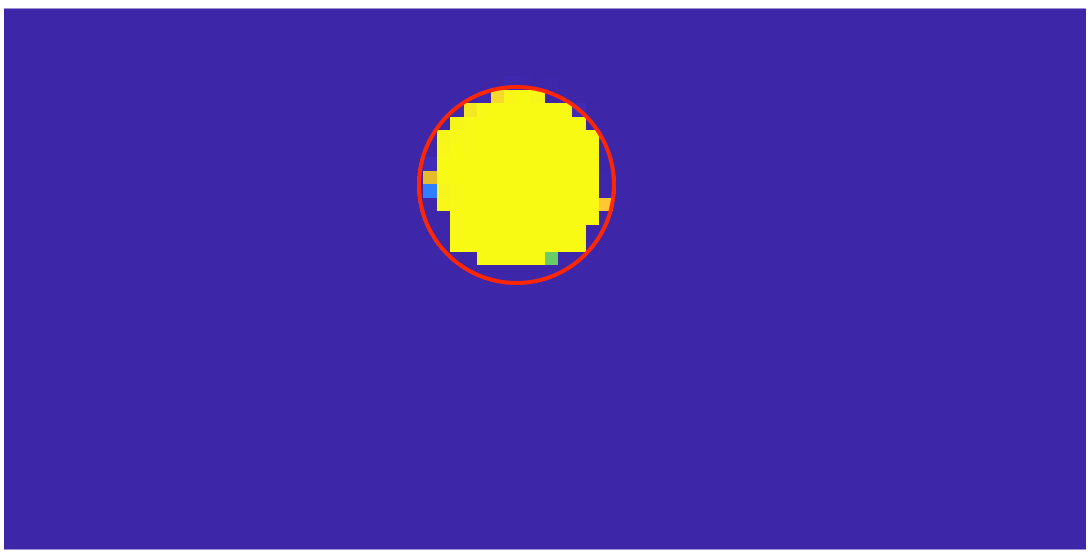}}
	
	\subfloat[\label{fig:G1_DOT}]{\includegraphics[width=\factor\textwidth]{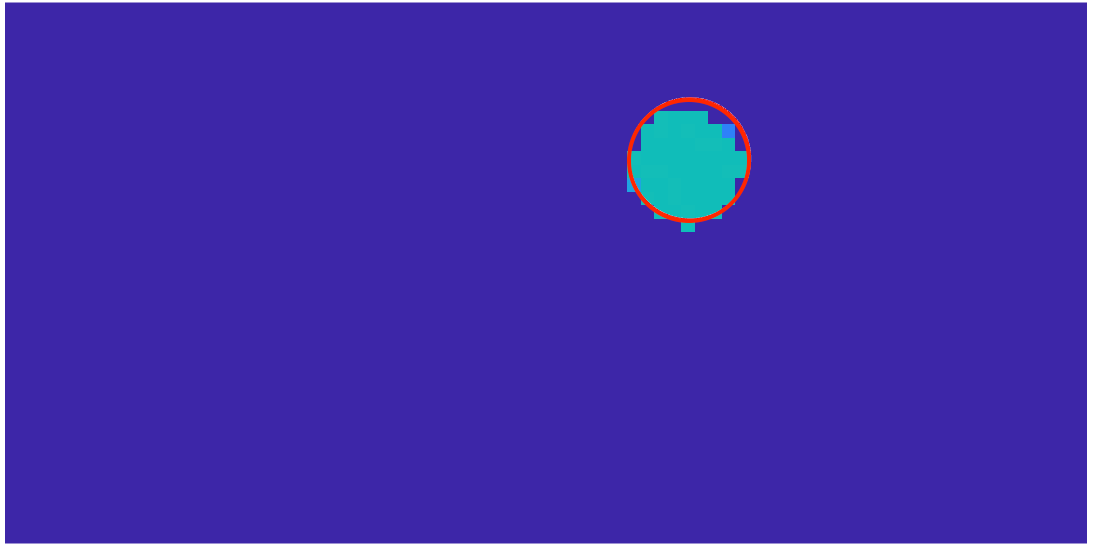}\hspace{0.025\textwidth}\includegraphics[width=\factor\textwidth]{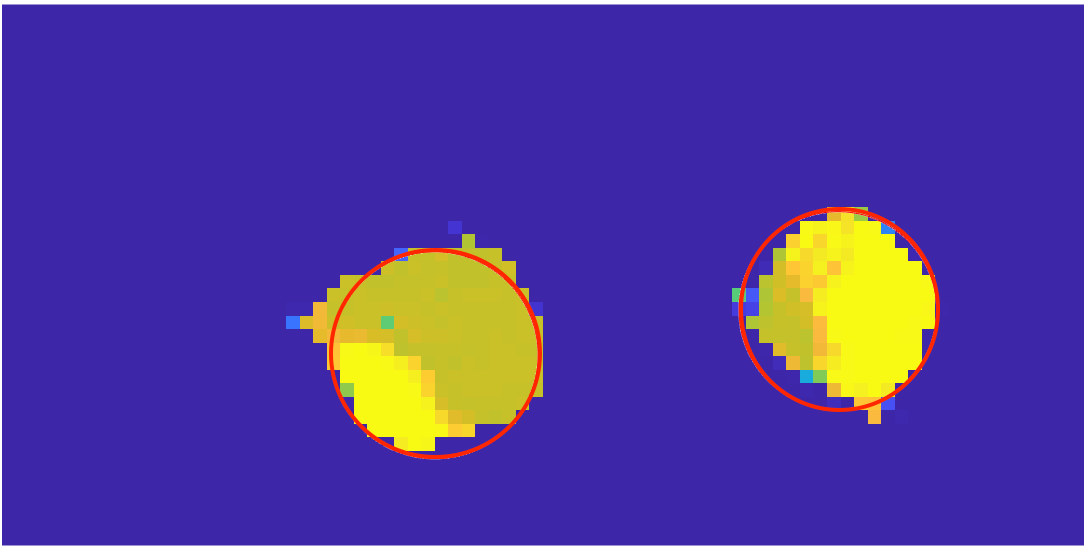}\hspace{0.025\textwidth}\includegraphics[width=\factor\textwidth]{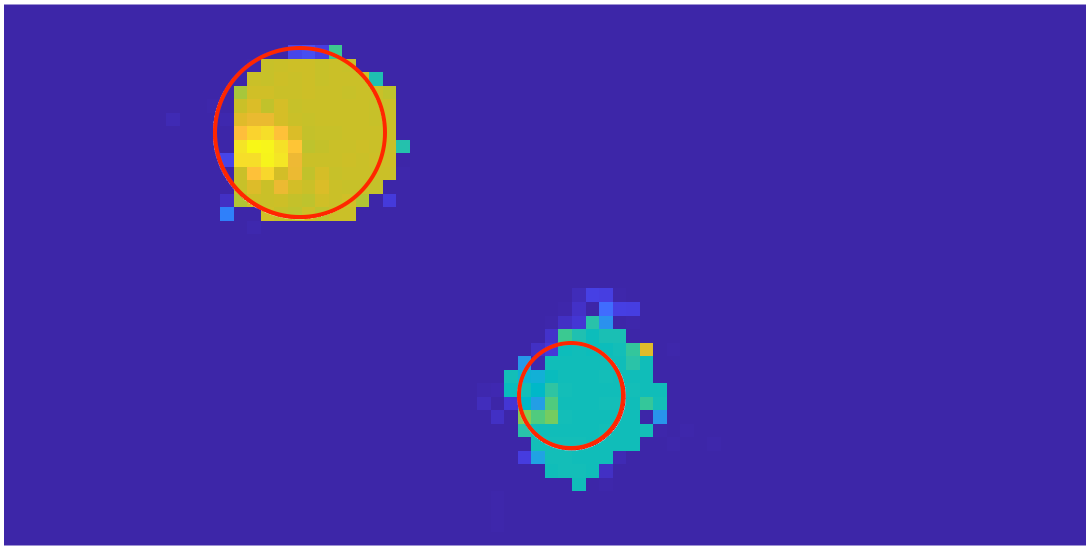}\hspace{0.025\textwidth}\includegraphics[width=\factor\textwidth]{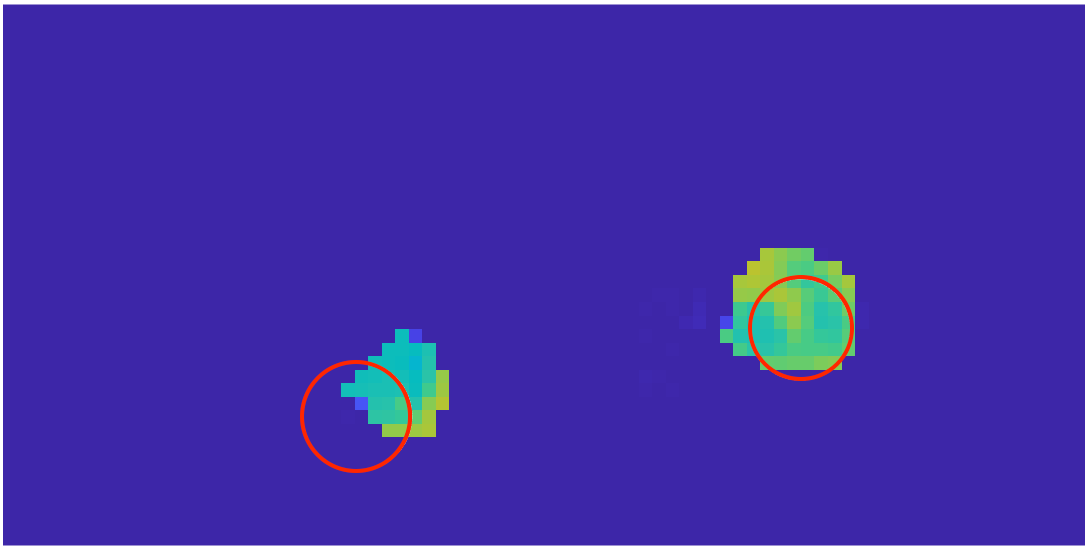}\hspace{0.025\textwidth}\includegraphics[width=\factor\textwidth]{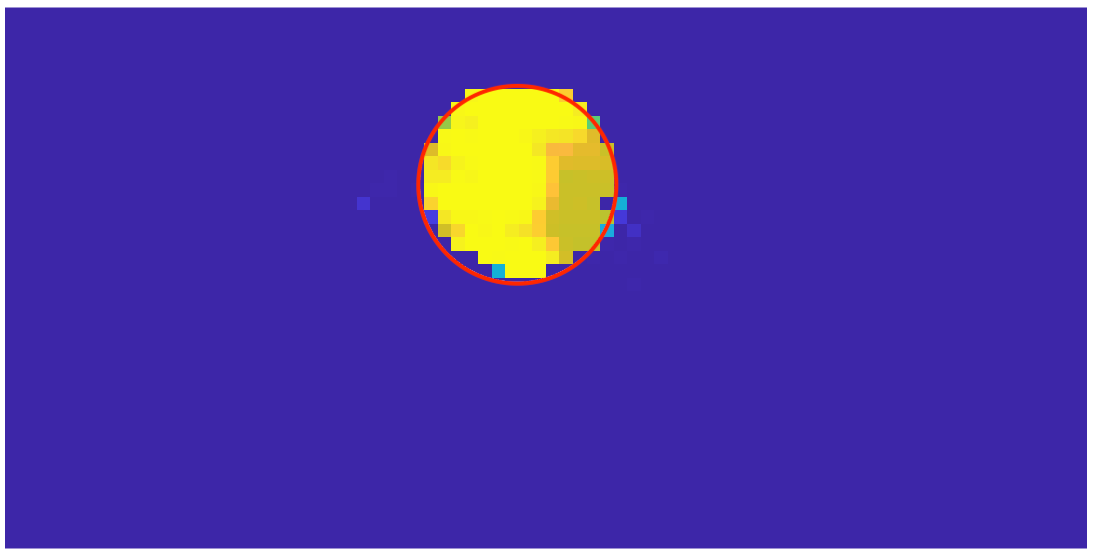}}

    \subfloat[\label{fig:G3_DOT}]{\includegraphics[width=\factor\textwidth]{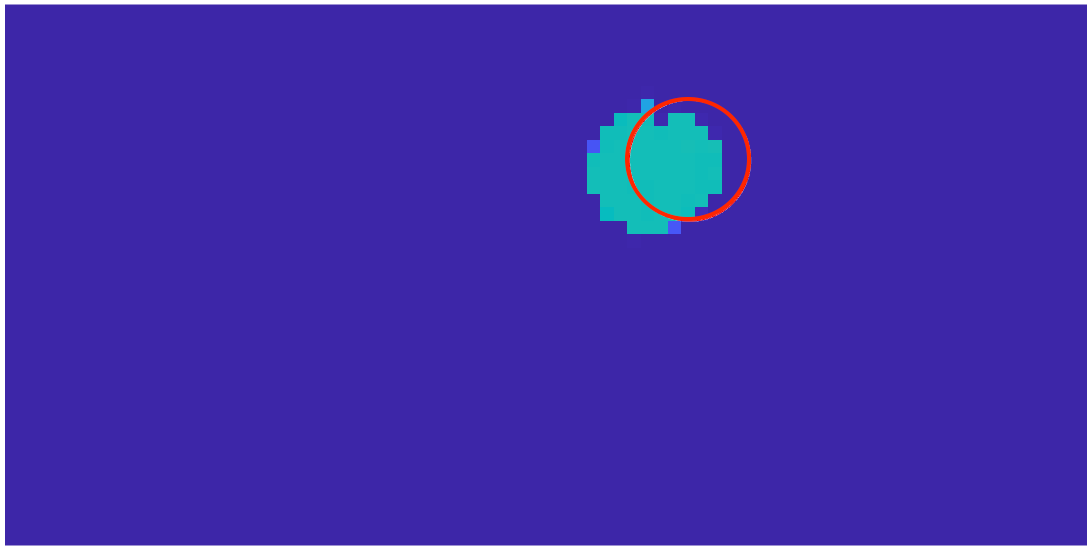}\hspace{0.025\textwidth}\includegraphics[width=\factor\textwidth]{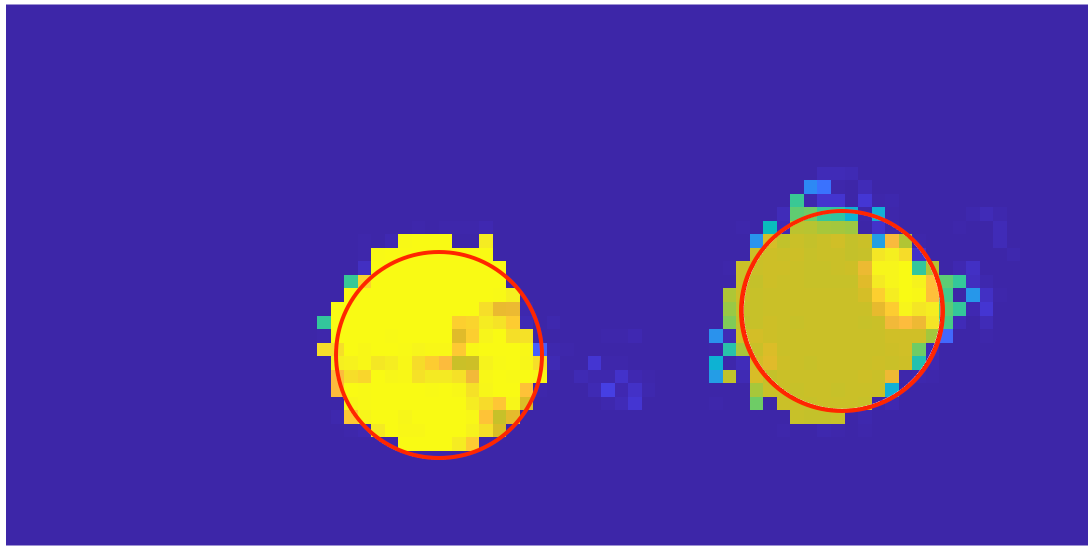}\hspace{0.025\textwidth}\includegraphics[width=\factor\textwidth]{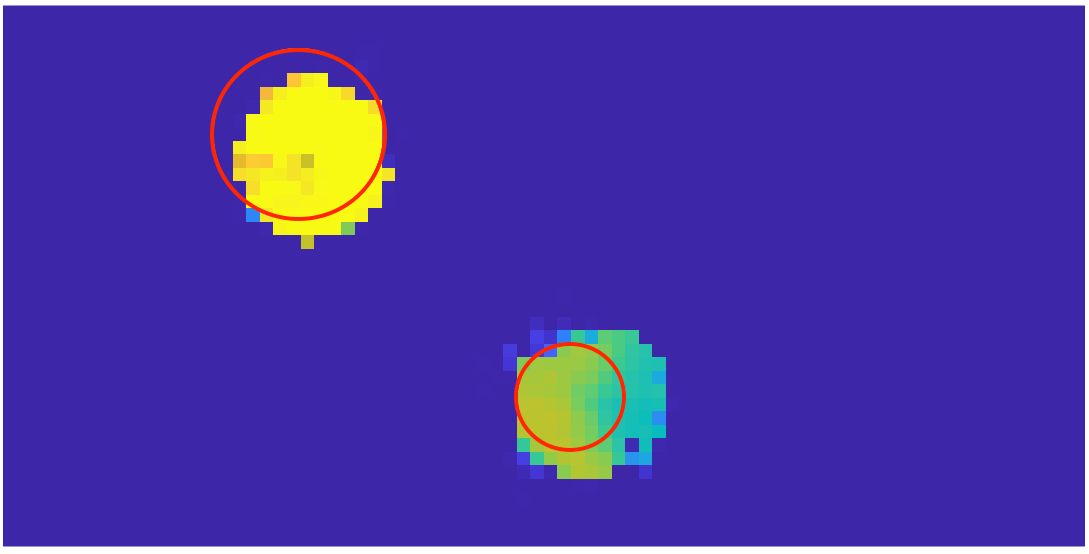}\hspace{0.025\textwidth}\includegraphics[width=\factor\textwidth]{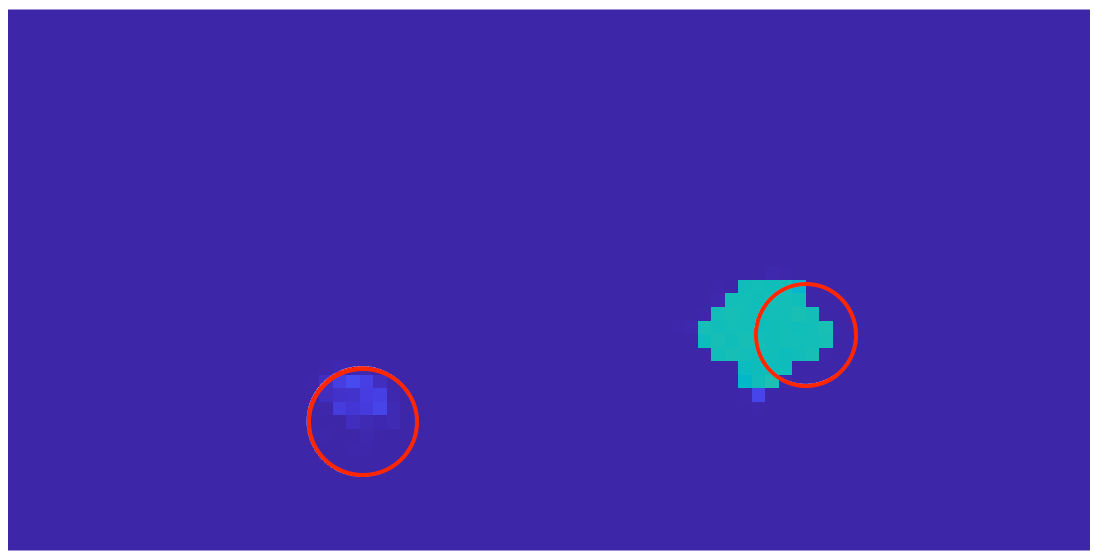}\hspace{0.025\textwidth}\includegraphics[width=\factor\textwidth]{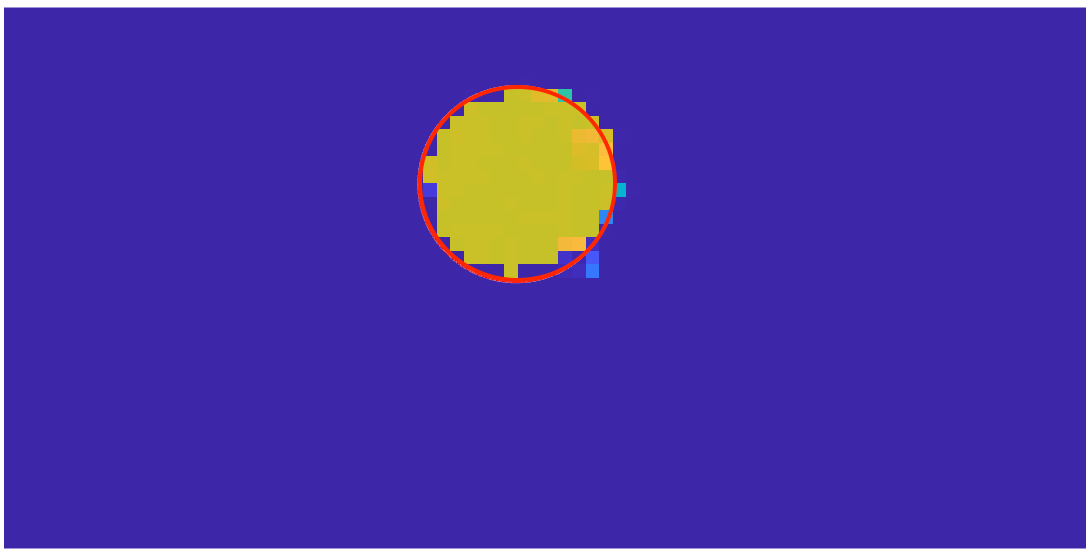}}

    \subfloat[\label{fig:G5_DOT}]{\includegraphics[width=\factor\textwidth]{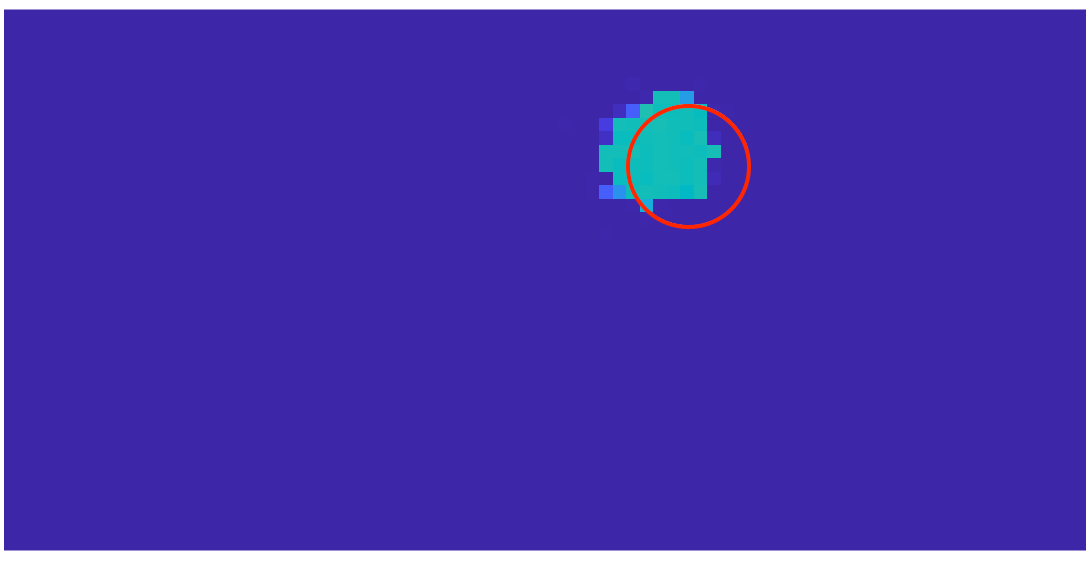}\hspace{0.025\textwidth}\includegraphics[width=\factor\textwidth]{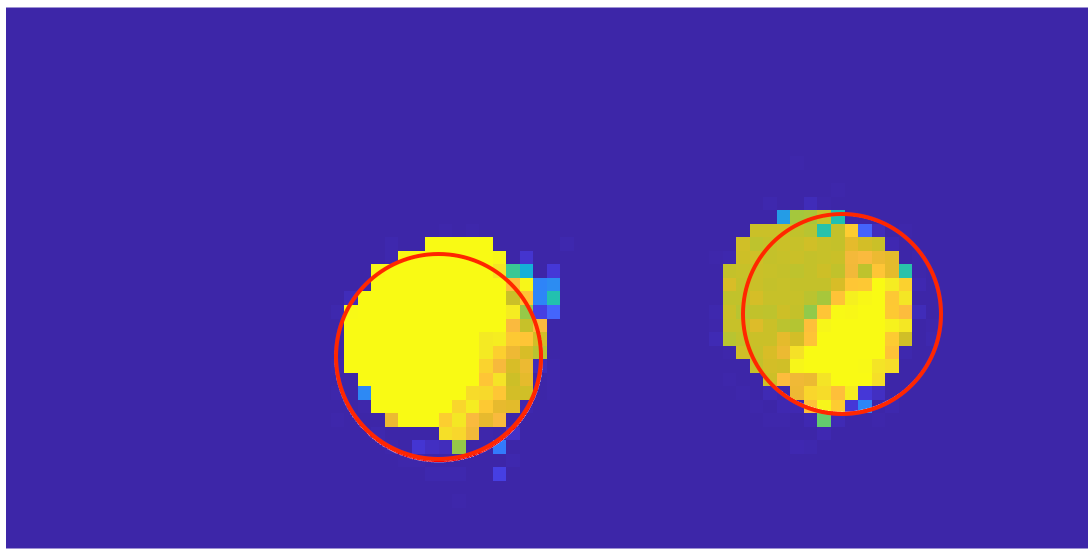}\hspace{0.025\textwidth}\includegraphics[width=\factor\textwidth]{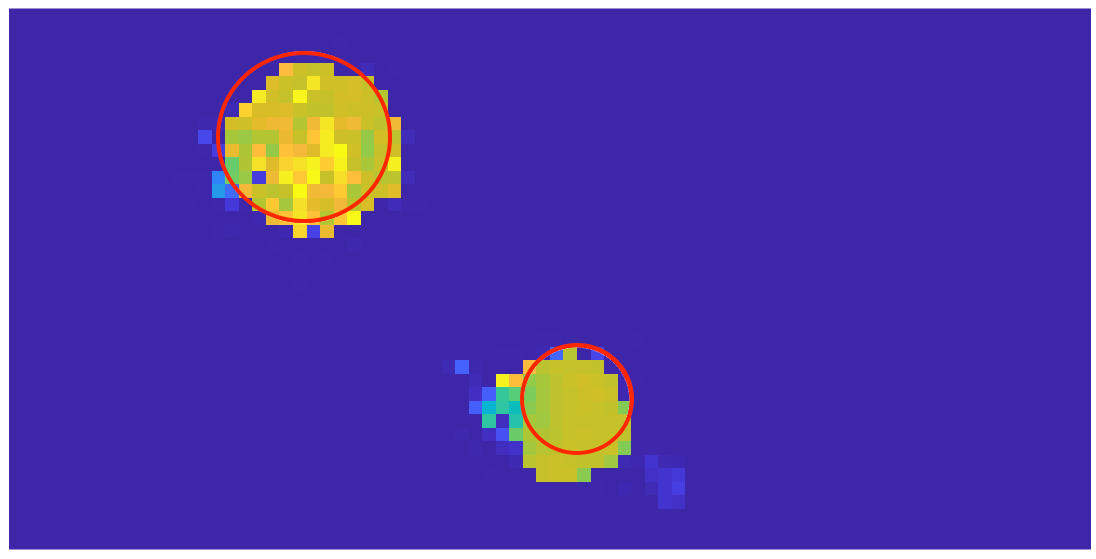}\hspace{0.025\textwidth}\includegraphics[width=\factor\textwidth]{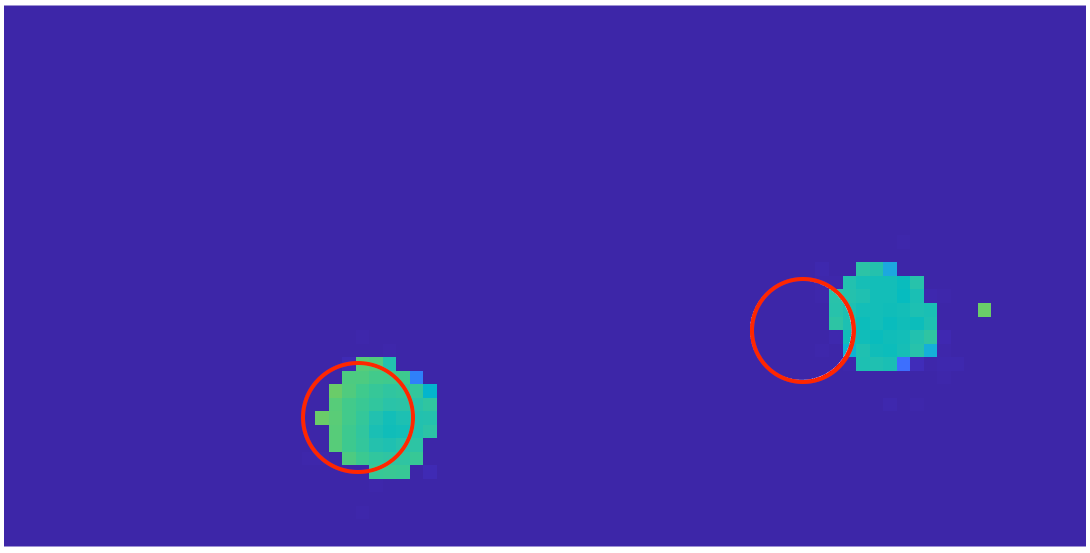}\hspace{0.025\textwidth}\includegraphics[width=\factor\textwidth]{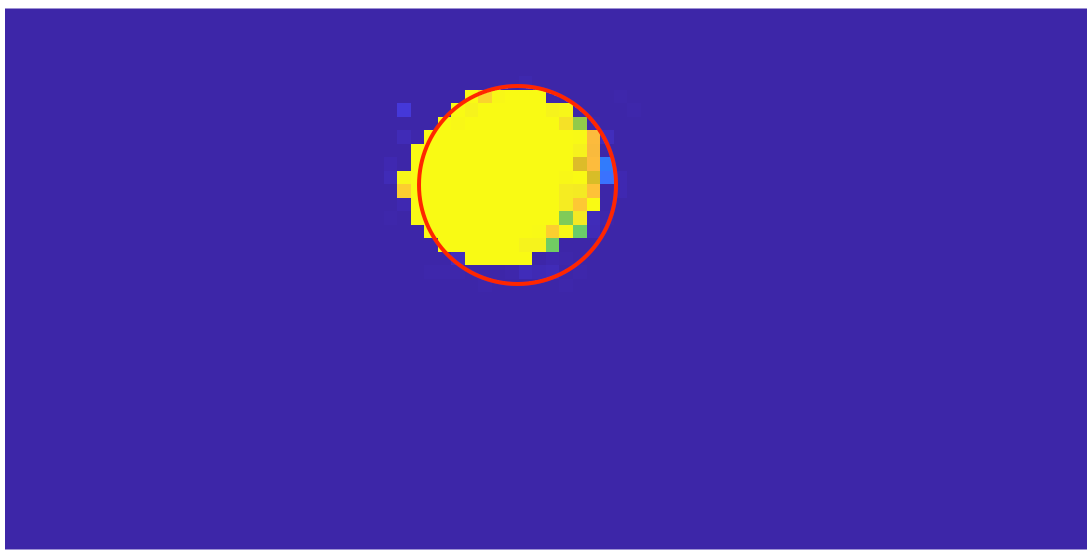}}  
    
\subfloat{\includegraphics[width=1\textwidth]{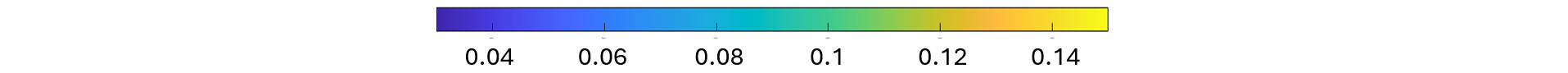}}
\end{center}
\caption{Reconstruction of the absorption coefficient distribution obtained with the Mod-DOT approach. (a): ground truth; (b): reconstruction without pretraining the AEs, before the application of the denoiser; (c): reconstruction with pretraining of the AEs 
before denoising; (d): reconstruction with pretraining of the AEs after denoising.(e): 1\% noise; (f): 3\% noise; (g): 5\% noise.
Rows (e) to (g) use pretrain and denoising.
Values are presented as $\mu_a/D$, the 
red circles denote the true contrast regions.}
\label{fig:fResults_DOT}
\end{figure*}

\subsubsection{Study of training strategies}
We investigate the effect of the training strategy
on the accuracy of the results. Namely, we compare the reconstructed images by the E2E 
(Fig.~\ref{fig:fResults_E2E}) and Mod-DOT architectures (Fig.~\ref{fig:fResults_DOT})
when training is performed monolithically on
the net or the two AEs are first pre-trained and
then the net is trained in a coupled fashion. 
It is apparent that for increasing 
noise, pretraining the AEs, and especially 
the signal-AE, yields superior results
for both architectures. 

\subsubsection{Study of out-of-distribution test samples}

\begin{figure}[htbp]
\centering
\includegraphics[width=.9\textwidth]{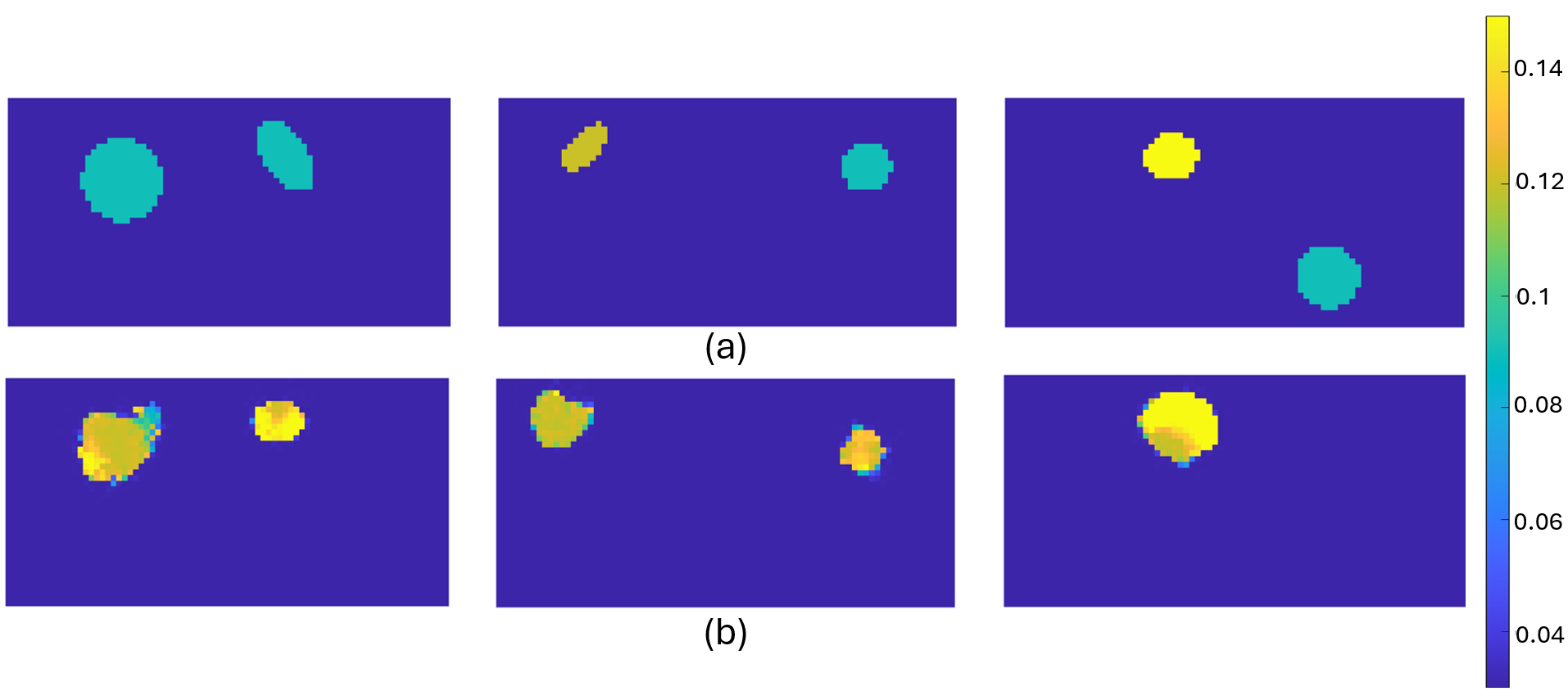}
\caption{Reconstruction of the absorption coefficient with out-of-distribution samples. Contrast regions have elliptical shapes with different locations, size and intensity. (a): ground-truth; (b): reconstruction obtained with the Mod-DOT approach trained on samples with
only circular contrast regions.
Values are presented as $\mu_a/D$.}
\label{fig:ellissi}
\end{figure}

We conclude our analysis by considering again the rectangular geometry 
and we try to reconstruct with the Mod-DOT
net 
absorption images with elliptical contrast regions. The net has been trained only on samples with circular contrast regions and we
present it with test samples having contrast regions made of ellipses with 
different location, size and orientation of the axes. 
We show in Fig.~\ref{fig:ellissi}
the reconstructed images of the
absorption coefficient for 3 different test samples. In general, the net is able to qualitatively reconstruct 
the correct location and intensity of the contrast regions. The shape of the region is not always correctly reconstructed. Observe, however, that
training only on circular regions
imposes a strong bias on the shape. In addition, the sinograms obtained with the elliptical contrast regions 
are minimally different from those
obtained from circular contrast regions (first two columns in Fig.~\ref{fig:ellissi}). There are  
samples where the reconstruction has a lower quality (third column in Fig.~\ref{fig:ellissi}). This may depend on the sinogram generated by a certain distribution of the coefficient,
even it is not easy to establish a general criterion. 

\section{Conclusion}
\label{sec:discu}
DOT reconstruction is a severely ill-conditioned problem which demands a careful treatment. 
Commonly used  strategies are based on physics-driven models
accompanied by priors which enforce constraints on the variance of the 
solution or in its sparsity.  These strategies often fail or are computationally
very intensive.  
In this work, we have investigated
a data-driven approach based on
a modular use of neural networks
which exploits the use of autoencoders. 
This architecture, called Mod-DOT, shares common elements with the Learned-SVD method originally proposed in~\cite{Boink19}
for general inverse problems, 
with the DL-ROM methods for parametrized
partial differential equations and
with the inversion strategy for EIT proposed in~\cite{seo2019learning}.
Specifically, our approach shares the idea to 
exploit the reduced dimension latent spaces obtained via an autoencoder
to improve the quality of the 
reconstruction. 
At the same time, the bridge-NN, which
connects the data and signal latent spaces, also acts here as a learned regularizator of the inversion.

\null

\noindent The Mod-DOT approach produces improved results with
respect to classic approaches and also with respect to
end-to-end neural network-based approaches.
This is more evident for higher levels of noise, the more realistic situation, where the quality of the reconstruction
remains acceptable.

\backmatter

\bmhead{Acknowledgments}
Alessandro Benfenati and Paola Causin are part of the group GNCS
(Gruppo Nazionale Calcolo Scientifico) of INdAM.

\bmhead{Author Contributions}
All the authors equally contributed in the development of the model,
MQ wrote the code and performed the numerical experiments. 
All the authors analyzed the results.    
PC wrote the paper, AB prepared Figs.8,9,10
and reviewed the manuscript. 

\bmhead{Funding}
We gratefully acknowledge the financial support of Italian MIUR--PRIN
Grant 20225STXSB {\em Sustainable Tomographic Imaging with Learning and
rEgularization (STILE)}

\bmhead{Availability of Data and Materials}
All the implemented NN architectures are available upon reasonable request.

%


\end{document}